\documentclass[11pt, a4paper]{article}
\usepackage{amsmath}
\usepackage{amsfonts}
\usepackage{amssymb}
\usepackage{ulem}

\usepackage{amsthm}
\usepackage{url}
\usepackage{dcolumn}

\usepackage{cite}
\usepackage{fancyhdr}

\usepackage{graphicx, color}
\usepackage{hyperref}
\usepackage{xcolor}

\newtheorem*{rem*}{Remark}

\begin{document}

\title{Pentagons and rhombuses that can form \\rotationally symmetric 
tilings}
\author{ Teruhisa SUGIMOTO$^{ 1), 2)}$ }
\date{}
\maketitle

{\footnotesize

\begin{center}
$^{1)}$ The Interdisciplinary Institute of Science, Technology and Art

$^{2)}$ Japan Tessellation Design Association

E-mail: ismsugi@gmail.com
\end{center}

}

{\small
\begin{abstract}
\noindent
In this study, various rotationally symmetric tilings that can be formed using 
pentagons that are related to rhombus are discussed. The pentagons can be 
convex or concave and can be degenerated into a trapezoid. If the pentagons 
are convex, they belong to the Type 2 family. Because the properties of pentagons 
correspond to those of rhombuses, the study also explains the correspondence 
between pentagons and various rhombic tilings.
\end{abstract}
}



\section{Introduction}
\label{section1}

In \cite{Sugimoto_2020_1} and \cite{Sugimoto_2020_2}, we introduced rotationally 
symmetric tilings and rotationally symmetric tilings (tiling-like patterns) with an 
equilateral convex polygonal hole at the center formed using convex pentagonal 
tiles.\footnote{
A \textit{tiling} (or \textit{tessellation}) 
of the plane is a collection of sets that are called tiles, which covers a plane without 
gaps and overlaps, except for the boundaries of the tiles. The term ``tile" refers to a 
topological disk, whose boundary is a simple closed curve. If all the tiles in a tiling are 
of the same size and shape, then the tiling is \textit{monohedral}~\cite{G_and_S_1987, 
wiki_pentagon_tiling}. In this study, a polygon that admits a monohedral tiling is called 
a \textit{polygonal tile}~\cite{Sugimoto_NoteTP, Sugimoto_2017_2}. Note that, in 
monohedral tiling, it admits the use of reflected tiles.
} 
These tilings have different connecting methods such as 
edge-to-edge\footnote{
A tiling by convex polygons is \textit{edge-to-edge} if any two convex polygons in a 
tiling are either disjoint or share one vertex or an entire edge in common. 
Then other case is \textit{non-edge-to-edge}~\cite{G_and_S_1987, Sugimoto_NoteTP, 
Sugimoto_2017_2}.
} 
and non-edge-to-edge. The convex pentagonal tiles forming the tilings belong to the 
Type 1 family.\footnote{
To date, fifteen families of convex pentagonal tiles, each of them referred to as a 
``Type," are known ~\cite{G_and_S_1987, Sugimoto_NoteTP, Sugimoto_2017_2, 
wiki_pentagon_tiling}. For example, if the sum of three consecutive angles in a convex 
pentagonal tile is  $360^ \circ $, the pentagonal tile belongs to the Type 1 family. 
Convex pentagonal tiles belonging to some families also exist. Known convex pentagonal 
tiles can form periodic tiling. In May 2017, Micha\"{e}l Rao declared that the complete 
list of Types of convex pentagonal tiles had been obtained (i.e., they have only 
the known 15 families), but it does not seem to be fixed as of March 
2020~\cite{wiki_pentagon_tiling}.
} 
Note that the convex pentagonal tiles in \cite{Sugimoto_2020_1} and \cite{Sugimoto_2020_2} 
are considered to be generated by bisecting equilateral concave octagons and equilateral 
convex hexagons, respectively.

\begin{figure}[htbp]
 \centering\includegraphics[width=14.5cm,clip]{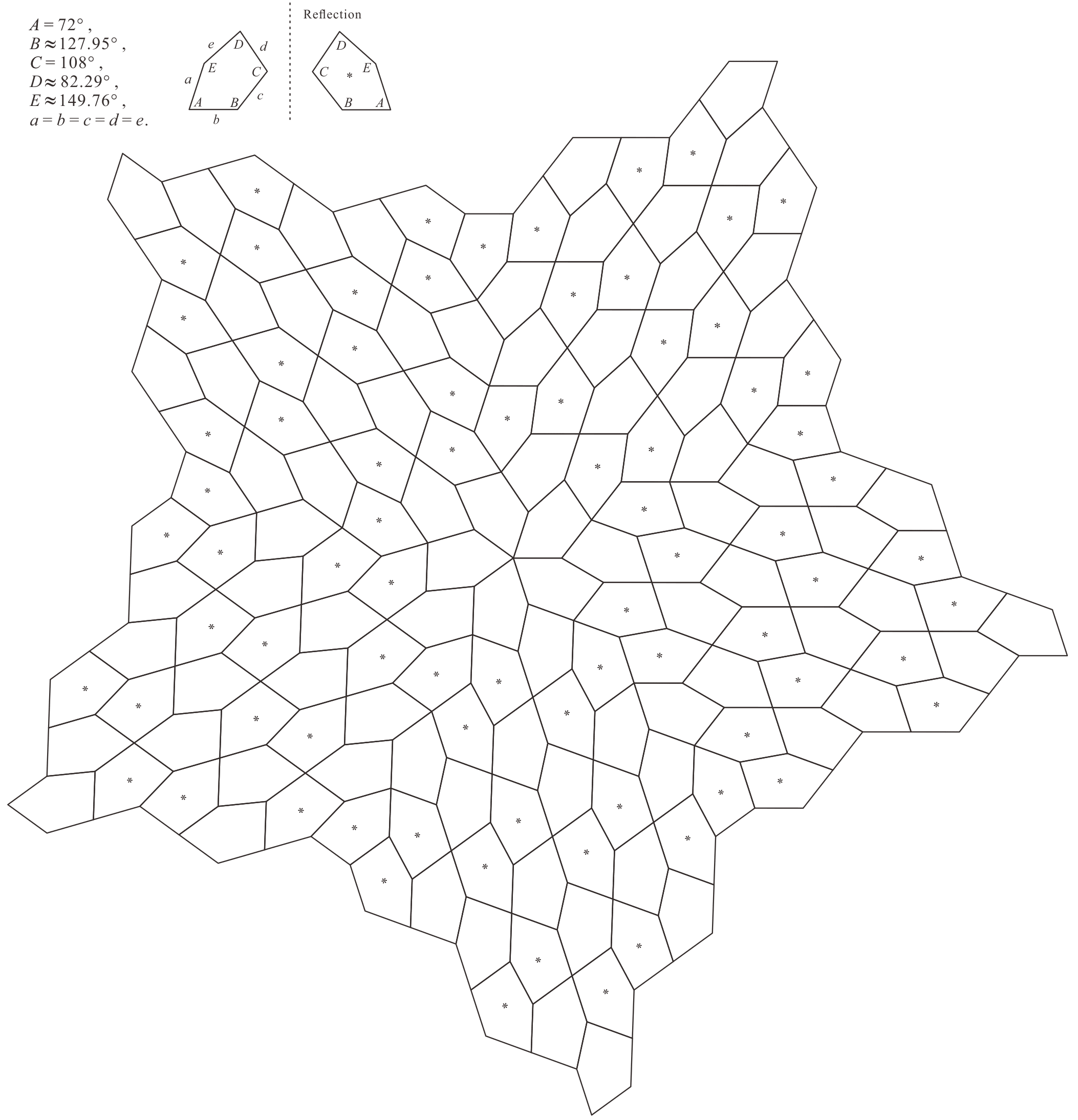} 
  \caption{{\small 
Five-fold rotationally symmetric tiling by an equilateral convex pentagonal 
tile that Hirschhorn, Hunt, and Zucca demonstrated (The figure is solely a 
depiction of  the area around the rotationally symmetric
 center, and the tiling can be spread in all directions as well)
 }
\label{fig01}
}
\end{figure}

\begin{figure}[htbp]
 \centering\includegraphics[width=14.5cm,clip]{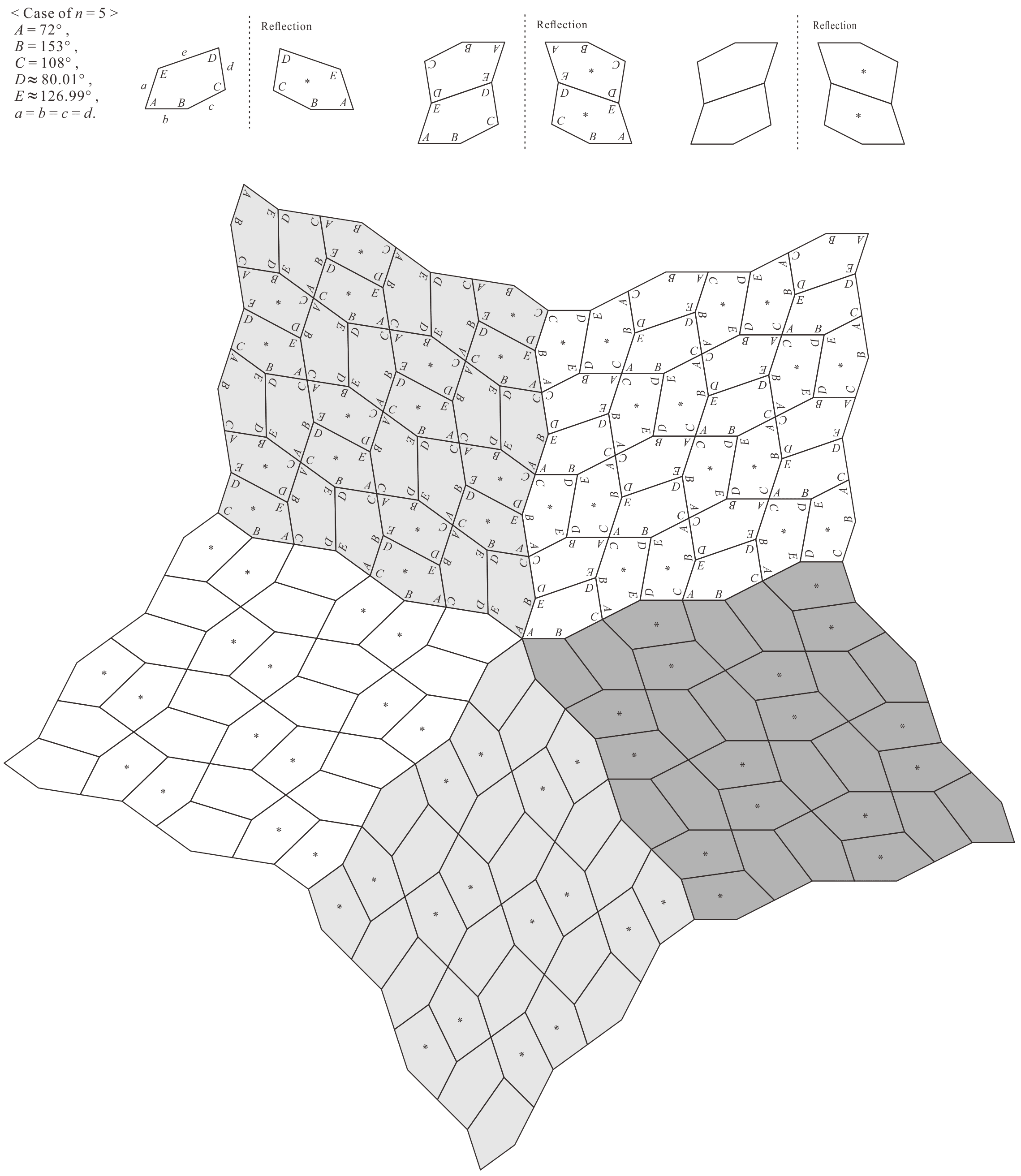} 
  \caption{{\small 
Five-fold rotationally symmetric tiling by a convex pentagonal tile 
with four equal-length edges 
 (The figure is solely a depiction of  the area around the rotationally symmetric
 center, and the tiling can be spread in all directions as well. 
Note that the gray area in the figure is used to clearly depict the structure)
} 
\label{fig02}
}
\end{figure}

Apart from the rotationally symmetric tilings with convex pentagonal tiles 
described above, Hirschhorn, Hunt, and Zucca demonstrated a five-fold rotationally 
symmetric tiling with equilateral convex pentagonal tiles belonging to the Type 2 
family, as shown in Figure~\ref{fig01} \cite{G_and_S_1987, H_and_H_1985, 
S_and_O_2009, wiki_pentagon_tiling, Zucca_2003}. In \cite{S_and_O_2009}, we considered 
edge-to-edge tilings with a convex pentagon having four equal-length edges 
and demonstrated that the convex pentagon in Figure~\ref{fig01} corresponds to 
a case of a convex pentagonal tile called ``C20-T2," which has five equal-length 
edges (i.e., equilateral edges) and an interior angle of $72^ \circ $. The 
results suggest that the five-fold rotationally symmetric tiling shown in 
Figure~\ref{fig01} can be formed using a convex pentagonal tile (C20-T2) with four 
equal-length edges, as shown in Figure~\ref{fig02}.

As in \cite{Sugimoto_2020_1} and \cite{Sugimoto_2020_2}, we expected that 
the convex pentagonal tile C20-T2 will be able to form not only five rotationally 
symmetric tilings, but also other rotationally symmetric tilings. We then confirm that 
C20-T2 is capable of forming such tilings. This study introduces the results 
obtained.

\renewcommand{\figurename}{{\small Figure.}}

\renewcommand{\figurename}{{\small Figure.}}

\renewcommand{\figurename}{{\small Figure.}}
\begin{figure}[htbp]
 \centering\includegraphics[width=12cm,clip]{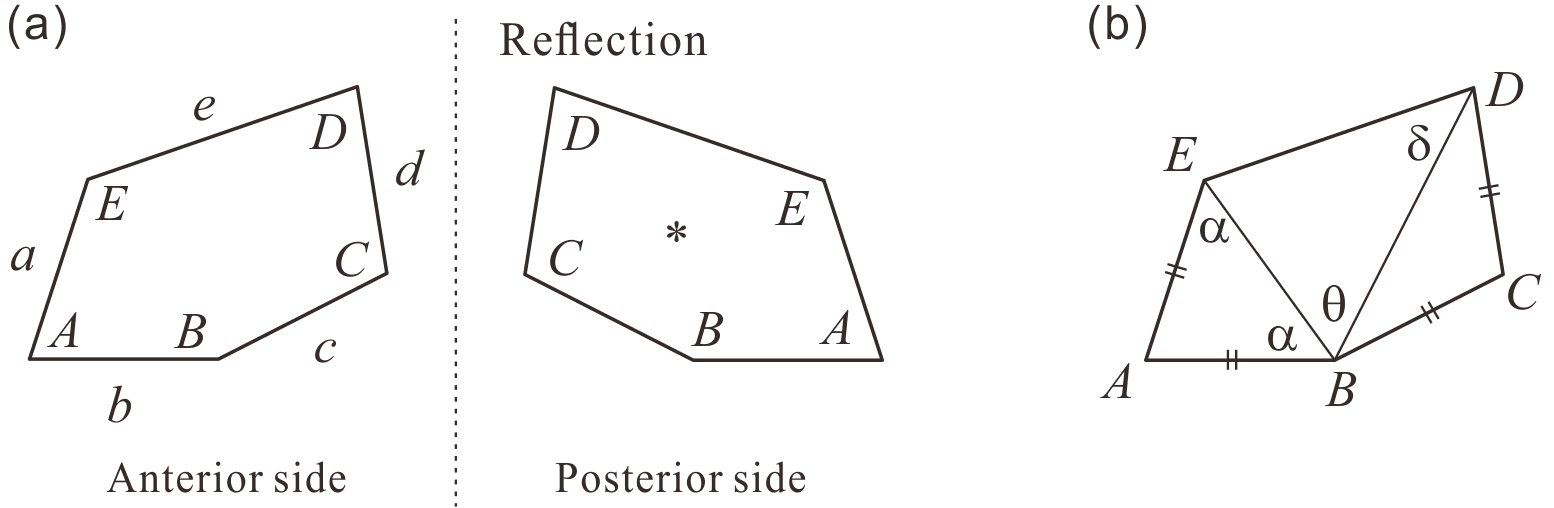} 
  \caption{{\small 
Nomenclature for vertices and edges of convex pentagon, and 
three triangles in the convex pentagonal tile C20-T2
} 
\label{fig03}
}
\end{figure}


\section{Conditions of pentagon that can form rotationally symmetric tilings}
\label{section2}

In this study, the vertices (interior angles) and edges of the pentagon will be 
referred to using the nomenclature shown in Figure~\ref{fig03}(a). C20-T2 shown 
in \cite{S_and_O_2009} is a convex pentagon that satisfies the conditions

\begin{equation}
\label{eq1}
\left\{ {\begin{array}{l}
 B + D + E = 360^ \circ , \\ 
 a = b = c = d, \\ 
 \end{array}} \right.
\end{equation}

\noindent
and can form the representative tiling (tiling of edge-to-edge version) of Type 2 
that has the relations ``$B + D + E = 360^ \circ , \;2A + 2C = 360^ \circ $." 
Because this convex pentagon has four equal-length edges, it can be divided into 
an isosceles triangle \textit{BCD}, an isosceles triangle \textit{ABE} with a base 
angle $\alpha $, and a triangle \textit{BDE} with $\angle DBE = \theta$ and 
$\angle BDE = \delta $, as shown in Figure~\ref{fig03}(b). 
Accordingly, using the relational expression for the interior angle of each 
vertex of C20-T2, the conditional expressions of (\ref{eq1}) can be rewritten 
as follows:

\begin{equation}
\label{eq2}
\left\{ {\begin{array}{l}
 A = 180^ \circ - 2\alpha , \\ 
 B = 90^ \circ + \theta , \\ 
 C = 2\alpha , \\ 
 D = 90^ \circ - \alpha + \delta , \\ 
 E = 180^ \circ + \alpha - \theta - \delta , \\ 
 a = b = c = d, \\ 
 \end{array}} \right.
\end{equation}

\noindent
where

\[
\delta = \tan ^{ - 1}\left( {\frac{\sin \theta }{\tan \alpha - \cos \theta 
}} \right)
\]

\noindent
and $0^ \circ < \alpha < 90^ \circ $ because $A > 0^ \circ $ and 
$C > 0^ \circ $~\cite{S_and_O_2009}. 
This pentagon has two degrees of freedom ($\alpha $ and $\theta $ 
parameters), besides its size. If the edge $e$ of this pentagon exists and the 
pentagon is convex, then $0^ \circ < \theta < 90^ \circ $. But depending 
on the value of $\alpha $, even if $\theta $ is selected in 
$\left( {0^ \circ ,\;90^ \circ } \right)$, the pentagon may not be convex 
or may be geometrically nonexistent. If $a = b = c = d = 1$, then the 
length of edge $e$ can be expressed as follows:

\[
e = 2\sqrt {1 - \sin (2\alpha )\cos \theta } .
\]

Let the interior angle of vertex $A$ be $\frac{360^ \circ }{n}$ (i.e., 
$\alpha = 90^  \circ-\frac{180^ \circ }{n}$) so that pentagons 
satisfying (\ref{eq2}) can form an $n$-fold rotationally 
symmetric tiling. (Remark: Due to the properties of the pentagons, 
the interior angle of vertex $C$, and not vertex $A$, will be 
$\frac{360^ \circ }{n}$.) Note that $n$ is an integer greater than or 
equal to three, because $C > 0^ \circ $. Therefore, the conditions 
of pentagonal tiles that can form $n$-fold rotationally symmetric 
tilings are expressed in (\ref{eq3}).
Note that the properties of the shape of pentagons that satisfy 
(\ref{eq3}) depending on the values of $n$ and $\theta $ are 
summarized in \ref{app}.

\begin{equation}
\label{eq3}
\left\{ {\begin{array}{l}
 A = \dfrac{360^ \circ }{n\strut}, \\ 
 B = 90^ \circ + \theta , \\ 
 C = 180^ \circ - \dfrac{360^ \circ }{n\strut}, \\ 
 D = \delta + \dfrac{180^ \circ }{n\strut}, \\ 
 E = 270^ \circ  - \theta - \delta - \dfrac{180^ \circ }{n\strut}, \\ 
 a = b = c = d. \\ 
 \end{array}} \right.
\end{equation}


\section{Relationships between pentagon and rhombus}
\label{section3}

\begin{figure}[htbp]
 \centering\includegraphics[width=15cm,clip]{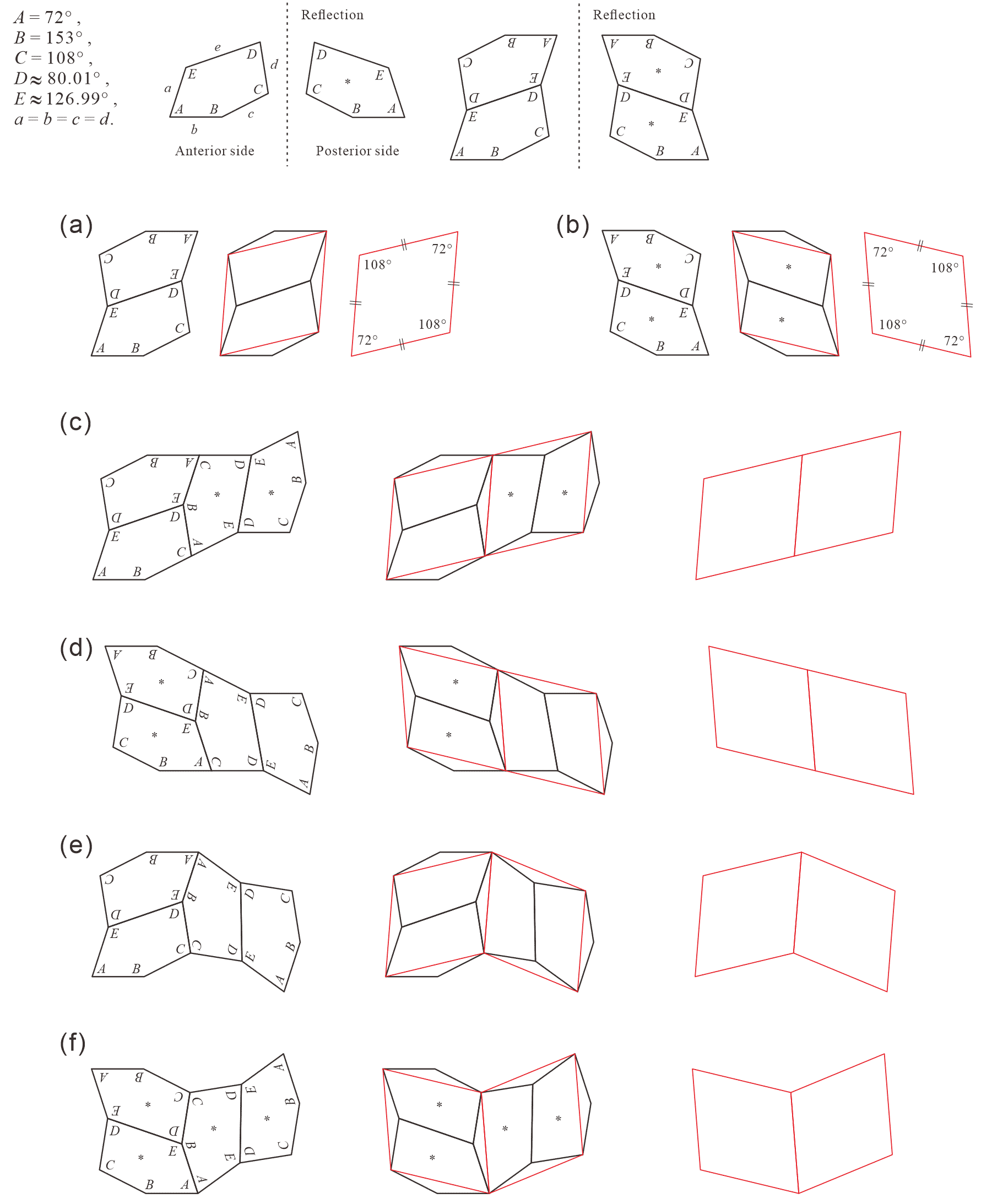} 
  \caption{{\small 
Relationships between pentagonal pair (basic unit) and rhombus
} 
\label{fig04}
}
\end{figure}

The convex pentagon shown in Figure~\ref{fig02} satisfies (\ref{eq3}), where 
$n = 5$ and $\theta = 63^ \circ $. Note that it is equivalent to the case where 
$\alpha = 54^ \circ$ and $\theta = 63^ \circ $, in (\ref{eq2}). By using this convex 
pentagon of Figure~\ref{fig02}, the method of forming tilings with pentagons 
satisfying the conditions of (\ref{eq2}) or (\ref{eq3}) is described below. 
In accordance with the relationship between the five interior angles of the 
pentagon, the vertex concentrations that can be always used in tilings are 
``$A + C = 180^ \circ ,\;B + D + E = 360^ \circ ,\;2A + 2C = 360^ \circ $." 
According to (\ref{eq2}) and (\ref{eq3}), the edge $e$ of the pentagon is the sole 
edge of different length. Therefore, the edge $e$ of one convex pentagon is always 
connected in an edge-to-edge manner with the edge $e$ of another convex 
pentagon. A pentagonal pair with their respective vertices $D$ and $E$ 
concentrated forms the basic unit of the tiling. This basic unit can be made of 
two types: a (anterior side) pentagonal pair as shown in Figure~\ref{fig04}(a) 
and a reflected (posterior side) pentagonal pair as shown in Figure~\ref{fig04}(b). 
Four different types of units, as shown in Figures~\ref{fig04}(c), \ref{fig04}(d), 
\ref{fig04}(e), and \ref{fig04}(f), are obtained by combining two pentagonal pairs 
shown in Figures~\ref{fig04}(a) and \ref{fig04}(b), so that $B + D + E = 360^ \circ $ 
can be assembled.

\renewcommand{\figurename}{{\small Figure.}}

\renewcommand{\figurename}{{\small Figure.}}
\begin{figure}[htbp]
 \centering\includegraphics[width=14.5cm,clip]{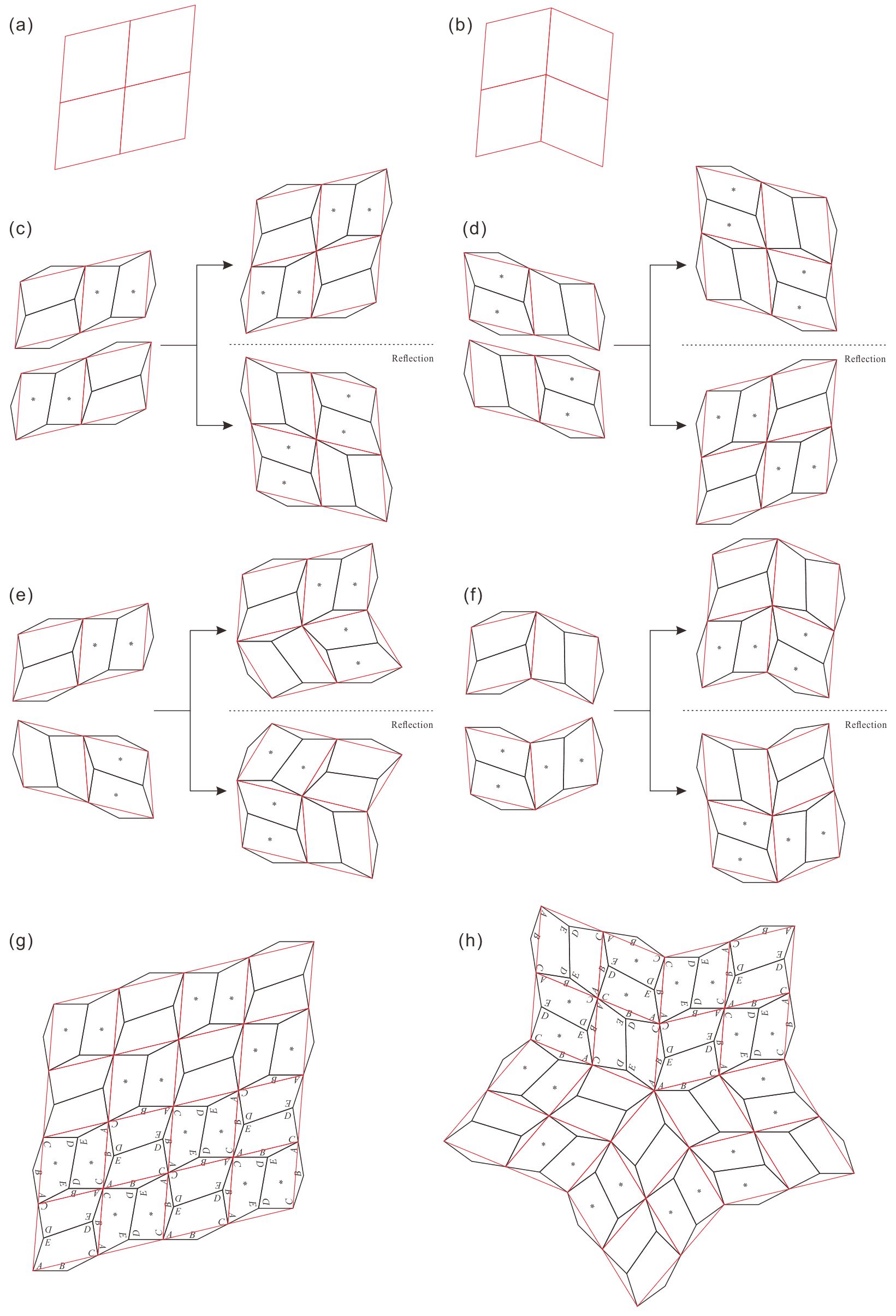} 
  \caption{{\small 
Combinations of rhombuses and pentagons
} 
\label{fig05}
}
\end{figure}

As shown in Figures~\ref{fig04}(a) and \ref{fig04}(b), a rhombus (red line), with an 
acute angle of $72^ \circ $, formed by connecting the vertices $A$ and $C$ of 
the pentagons, is applied to each basic unit of the pentagonal pair. (Remark: In 
this example, because the interior angle of the vertex $A$ is $72^ \circ $, the rhombus 
has an acute angle of $72^ \circ $. That is, the interior angles of the rhombus 
corresponding to the pair of pentagons in Figures~\ref{fig04}(a) and \ref{fig04}(b) 
are the same as the interior angles of vertices $A$ and $C$ in (\ref{eq2}) and 
(\ref{eq3}).) Consequently, the parts of pentagons that protrude from the rhombus 
match exactly with the parts that are more dented than the rhombus 
(refer to Figures~\ref{fig04}(c), \ref{fig04}(d), \ref{fig04}(e), and \ref{fig04}(f)). 
In fact, tilings in which ``$B + D + E = 360^ \circ ,\;2A + 2C = 360^ \circ $" 
using pentagons satisfying (\ref{eq2}) and (\ref{eq3}) are equivalent to 
 rhombic tilings (tilings formed by rhombuses). (Though a rhombus is a single 
entity, considering its internal pentagonal pattern, it will be considered as 
two entities.)

Rhombuses have two-fold rotational symmetry and two axes of reflection 
symmetry passing through the center of the rotational symmetry (hereafter, 
this property is described as $D_{2}$ symmetry\footnote{ 
``$D_{2}$" is based on the Schoenflies notation for symmetry in a two-dimensional point 
group~\cite{wiki_point_group, wiki_schoenflies_notation}. ``$D_{n}$" represents an 
$n$-fold rotation axis with $n$ reflection symmetry axes. The notation 
for symmetry is based on that presented in \cite{Klaassen_2016}.
}). 
Therefore, the rhombus and the reflected rhombus have identical outlines. 
Thus, the two methods of concentrating the four rhombic vertices at a 
point without gaps or overlaps are: Case (i) an arrangement by parallel 
translation as shown in Figure~\ref{fig05}(a); Case (ii) an arrangement by rotation 
(or reflection) as shown in Figure~\ref{fig05}(b). This concentration corresponds to 
forming a ``$2A + 2C = 360^ \circ $" at the center by four pentagons. In Case 
(i), because the pentagonal vertices circulate as ``$A \to C \to A \to C$" at the 
central ``$2A + 2C = 360^ \circ $," one combination (refer to Figure~\ref{fig05}(c)) is 
obtained by using two units of Figure~\ref{fig04}(c) and another combination (refer to  
Figure~\ref{fig05}(d)) is obtained by using two units of Figure~\ref{fig04}(d). In 
Case (ii), because the pentagonal vertices circulate as ``$A \to A \to C \to C$" at the 
central ``$2A + 2C = 360^ \circ $," one combination (refer to Figure~\ref{fig05}(e))
 is obtained by using units of Figures~\ref{fig04}(c) and \ref{fig04}(d), 
and another combination (refer to Figure~\ref{fig05}(f)) is obtained 
by using units of Figures~\ref{fig04}(e) and \ref{fig04}(f).

Only when the unit comprising eight pentagons in Figures~\ref{fig05}(c) or \ref{fig05}(d) 
are arranged in a parallel manner, a tiling, as shown in Figure~\ref{fig05}(g), is 
formed that represents a tiling of Type 2, in which  ``$B + D + E = 360^ \circ ,\; 
2A + 2C = 360^ \circ $." Because rhombuses can form rhombic belts by translation 
in the same direction vertically, rhombic tilings can also be formed by the belts 
that are freely connected horizontally by the connecting method shown in 
Figures~\ref{fig05}(a) and \ref{fig05}(b). Further, pentagonal tilings (tilings formed 
by pentagons) corresponding to those rhombic tilings can be formed.

When $n$ vertices, with interior angles of $\frac{360^ \circ}{n}$, of $n$ rhombuses are 
concentrated at a point, an $n$-fold rotationally symmetric arrangement is 
formed, with adjacent rhombuses connected as shown in Figure~\ref{fig05}(b). 
Therefore, an $n$-fold rotationally symmetric tiling with rhombuses can 
be formed by dividing each rhombus, in that arrangement, into similar shapes. 
By converting the rhombuses of such rhombic tiling into pentagons satisfying 
(\ref{eq3}), the rotationally symmetric tilings with convex pentagons can 
be obtained (refer to Figure~\ref{fig05}(h)). Therefore, when forming $n$-fold rotationally 
symmetric tilings from a pentagon satisfying (\ref{eq3}), the pentagonal arrangement 
can be known from the corresponding $n$-fold rotationally symmetric tilings with 
a rhombus.


\section{Rotationally symmetric tilings}
\label{section4}

Table~\ref{tab1} presents some of the relationships between the interior angles of 
convex pentagons satisfying (\ref{eq3}) that can form the $n$-fold rotationally 
symmetric edge-to-edge tilings. (For $n = 3\!-\!10, 16$, tilings with convex pentagonal 
tiles are drawn. For further details, Figures~\ref{fig02}, \ref{fig06}--\ref{fig13}.) Due to 
the presence of parameter $\theta $ in (\ref{eq3}), the shapes of convex pentagons 
that satisfy (\ref{eq3}) and can form an $n$-fold rotationally symmetric tiling are not 
fixed. Therefore, each example presented in Table~\ref{tab1} is a pentagon with a 
convex shape that can form an $n$-fold rotationally symmetric tiling. If the 
pentagons satisfying (\ref{eq3}) are convex, the $n$-fold rotationally symmetric 
tilings with the pentagonal tiles are connected in an edge-to-edge manner 
and have no axis of reflection symmetry.\footnote{
Hereafter, a figure with $n$-fold rotational symmetry without reflection is described 
as $C_{n}$ symmetry. ``$C_{n}$" is based on the Schoenflies notation for symmetry in 
a two-dimensional point group~\cite{wiki_point_group, wiki_schoenflies_notation}.
} 
The reason for this lack of symmetry is that the units comprising pentagons 
corresponding to that rhombus with $D_{2}$ symmetry (refer to Figures~\ref{fig04}(a), 
\ref{fig04}(b), \ref{fig05}(g), etc.) have $C_{2}$ symmetry.

\begin{table}[t]
 \begin{center}
{\small
\caption[Table 1]{
Example of interior angles of convex pentagons satisfying 
(\ref{eq3}) that can form the $n$-fold rotationally symmetric tilings
}
\label{tab1}
}
\begin{tabular}
{c| D{.}{.}{2} D{.}{.}{2} D{.}{.}{2} D{.}{.}{2} D{.}{.}{2} |c|c}
\hline
\raisebox{-1.50ex}[0cm][0cm]{$n$}& 
\multicolumn{5}{c|}{\shortstack{ Value of interior angle (degree) } } & 
\raisebox{-4.6ex}[0.7cm][0.5cm]{\shortstack{Edge \\length \\of $e$}} & 
\raisebox{-3.0ex}[0.7cm][0.5cm]{\shortstack{Figure \\number}} \\

 & 
\textit{A} & 
\textit{B}& 
\textit{C}& 
\textit{D}& 
\textit{E}& 
 & 
  \\
\hline
3& 
120& 
151& 
60& 
143.96 & 
65.04 & 
1.523 & 
\ref{fig06} \\
\hline
4& 
90& 
151& 
90& 
104.5 & 
104.5 & 
1.436 & 
\ref{fig07} \\
\hline
5& 
72& 
153& 
108& 
80.01 & 
126.99 & 
1.508 & 
\ref{fig02} \\
\hline
6& 
60& 
151& 
120& 
65.04 & 
143.96 & 
1.523 & 
\ref{fig08} \\
\hline
7& 
51.43 & 
146& 
128.57 & 
54.37 & 
159.63 & 
1.500 & 
\ref{fig09} \\
\hline
8& 
45& 
151& 
135& 
46.89 & 
162.11 & 
1.621 & 
\ref{fig10} \\
\hline
9& 
40& 
150& 
140& 
41.07 & 
168.93 & 
1.648 & 
\ref{fig11} \\
\hline
10& 
36& 
156& 
144& 
36.88 & 
167.12 & 
1.745 & 
\ref{fig12} \\
\hline
11& 
32.73 & 
156& 
147.27 & 
33.31 & 
170.69 & 
1.766 & 
 \\
\hline
12& 
30& 
160& 
150& 
30.49 & 
169.51 & 
1.821 & 
 \\
\hline
13& 
27.69 & 
160& 
152.31 & 
28.04 & 
171.96 & 
1.834 & 
 \\
\hline
14& 
25.71 & 
164& 
154.29 & 
26.03 & 
169.97 & 
1.877 & 
 \\
\hline
15& 
24& 
164& 
156& 
24.25 & 
171.75 & 
1.885 & 
 \\
\hline
16& 
22.5& 
170& 
157.5& 
22.72 & 
167.28 & 
1.932 & 
\ref{fig13} \\
\hline
17& 
21.18 & 
170& 
158.82 & 
21.36 & 
168.64 & 
1.936 & 
 \\
\hline
18& 
20& 
170& 
160& 
20.16 & 
169.84 & 
1.940 & 
 \\
\hline
...& 
...& 
...& 
...& 
...& 
...& 
...& 
 \\
\hline
\end{tabular}

\end{center}

\end{table}

The pentagons with $n = 3$ and $n = 6$ correspond to rhombuses with an acute 
angle of $60^ \circ $ (i.e., they correspond to tiling of an equilateral 
triangle), and these pentagons are opposite to each other. (In Table~\ref{tab1}, 
the interior angle of vertex $B$ is chosen to have the same value in both the 
cases. Note that the convex pentagonal tiles with $n = 3$ and $n = 6$ belong 
to both the Type 2 and Type 5 families~\cite{G_and_S_1987, Sugimoto_NoteTP, 
Sugimoto_2016, Sugimoto_2017_2, wiki_pentagon_tiling}.) According to this 
relationship, these tiles can form tilings with three-fold rotational symmetry 
that have six-fold rotational symmetry at the intersection of tilings, as shown 
in Figure~\ref{fig14}. Also, in addition to $2A + 2C = 360^ \circ $, 
``$3C = 360^ \circ ,\;4A + C = 360^ \circ ,\;6A = 360^ \circ $" 
are valid in these tilings (refer to Figure~\ref{fig15}). In particular, consider the 
unit comprising six pentagons, as shown in Figure~\ref{fig15}(a), that has 
outline shape with $D_{3}$ symmetry. The pentagons in this unit can be 
reversed freely (i.e., a unit comprising six anterior pentagons can be freely 
exchanged with a unit comprising six posterior pentagons). Therefore, 
various patterns, as shown in Figures~\ref{fig16} and \ref{fig17}, can be 
generated by the pentagon corresponding to the rhombus with an 
acute angle of $60^ \circ $.

In the case of $n = 4$, $A = C$ and $D = E$, and the pentagon has a line of 
symmetry connecting the vertex $B$ to the midpoint of the edge $e$ 
(refer to Figure~\ref{fig07}), i.e., there is no distinction between its anterior and 
posterior sides --- in the figures of this study, the posterior pentagons are 
marked with an asterisk mark. Accordingly, the rhombus corresponding to 
this case is a square. Therefore, the unit comprising eight pentagons 
corresponding to Figure~\ref{fig05}(a) has $C_{4}$ symmetry. The convex 
pentagonal tiling of this case is called Cairo tiling, and the convex pentagonal 
tiles belong to both the Type 2 and Type 4 families~\cite{G_and_S_1987, 
Sugimoto_NoteTP, Sugimoto_2016, Sugimoto_2017_2, wiki_pentagon_tiling}.

\begin{table}[htb]
 \begin{center}
{\small
\caption[Table 2]{
Trapezoids based on pentagons satisfying (\ref{eq3}) 
that can form the $n$-fold rotationally symmetric tilings
}
\label{tab2}
}
\begin{tabular}
{c| D{.}{.}{2} D{.}{.}{2} D{.}{.}{2} D{.}{.}{2} D{.}{.}{2} |c |c}
\hline
\raisebox{-1.50ex}[0cm][0cm]{$n$}& 
\multicolumn{5}{c|}{\shortstack{ Value of interior angle (degree) } } & 
\raisebox{-4.6ex}[0.7cm][0.5cm]{\shortstack{Edge \\length \\of $e$}} & 
\raisebox{-3.0ex}[0.7cm][0.5cm]{\shortstack{Figure \\number}} \\

 & 
\textit{A} & 
\textit{B}& 
\textit{C}& 
\textit{D}& 
\textit{E}& 
 & 
  \\
\hline
5& 
72& 
108& 
108& 
72& 
180& 
0.618 & 
\ref{fig18} \\
\hline
6& 
60& 
120& 
120& 
60& 
180& 
1 & 
\ref{fig19} \\
\hline
7& 
51.43 & 
128.57 & 
128.57 & 
51.43 & 
180& 
1.247 & 
 \\
\hline
8& 
45& 
135& 
135& 
45& 
180& 
1.414 & 
\ref{fig20} \\
\hline
9& 
40& 
140& 
140& 
40& 
180& 
1.532 & 
 \\
\hline
10& 
36& 
144& 
144& 
36& 
180& 
1.618 & 
 \\
\hline
11& 
32.73 & 
147.27& 
147.27& 
32.73 & 
180& 
1.683 & 
 \\
\hline
12& 
30& 
150& 
150& 
30& 
180& 
1.732 & 
 \\
\hline
13& 
27.69 & 
152.31 & 
152.31 & 
27.69 & 
180& 
1.771 & 
 \\
\hline
14& 
25.71 & 
154.29 & 
154.29 & 
25.71 & 
180& 
1.802 & 
 \\
\hline
15& 
24 & 
156 & 
156 & 
24 & 
180& 
1.827 & 
 \\
\hline
16& 
22.5 & 
157.5 & 
157.5 & 
22.5 & 
180& 
1.848 & 
 \\
\hline
17& 
21.18 & 
158.82 & 
158.82 & 
21.18 & 
180& 
1.865 & 
 \\
\hline
18& 
20 & 
160 & 
160 & 
20 & 
180& 
1.879 & 
 \\
\hline
...& 
...& 
...& 
...& 
...& 
...& 
...& 
 \\
\hline
\end{tabular}

\end{center}

\end{table}

Equilateral pentagons that satisfy (\ref{eq3}) exist, provided $n = 4, 5, 6, 7$. 
The pentagons that are convex and have equilateral edges are the cases with 
$n = 4$ ($B \approx 131.41^ \circ$) and $n = 5$ ($B \approx 127.95^ \circ$). 
Figure~\ref{fig01} shows the five-fold rotationally symmetric tiling with 
equilateral convex pentagons with $n = 5$. In the case of $n = 6$, the shape 
changes to a trapezoid, and the trapezoid can form three-fold or six-fold 
rotationally symmetric tiling (refer to Figures~\ref{fig19} and \ref{fig21}. Note 
that the line corresponding to edge $e$ in the figures is shown as a blue line). 
In the case of $n = 7$, the pentagon is concave and can form a seven-fold 
rotationally symmetric tiling (refer to Figure~\ref{fig25})~\cite{Zucca_2003}

Here, let us introduce some $n$-fold rotationally symmetric tilings with 
$C_{n}$ symmetry formed of trapezoids based on pentagons that satisfies 
(\ref{eq3}), similar to the equilateral pentagonal case with $n = 6$. If the 
pentagons satisfying (\ref{eq3}) with $n \ge 5$ have $\theta = 90^ \circ - A$, 
then ``$A = D,\;B = C,\;E = 180^ \circ $." Therefore, they are trapezoids 
with a line of symmetry. Table~\ref{tab2} presents some of these trapezoids. 
(For $n = 5, 6, 8$, tilings with trapezoidal tiles are drawn. For further details, 
Figures~\ref{fig18}--\ref{fig20}. Note that the line corresponding to edge $e$ 
in the figures is shown as a blue line.) Because the pentagons with $n = 3$ and 
$n = 6$ are opposite to each other, the trapezoid for the case of $n = 6$ 
can form three-fold or six-fold rotational symmetry tilings and mixed tilings, 
as shown in Figure~\ref{fig21}. 
The trapezoid for the case of $n = 6$ corresponds to a rhombus with an 
acute angle of $60^ \circ $, and the shape shown in Figure~\ref{fig15}(a) 
corresponds to Figure~\ref{fig22}(a). Therefore, similar to the case of a convex 
pentagon, various patterns can be generated by this trapezoid. At the same 
time, focusing on the equilateral triangles that appear in the tiling, as 
shown in Figure~\ref{fig22}(b), various patterns can be generated by 
replacing the trapezoids without using the shape shown in Figure~\ref{fig22}(a).

\begin{table}[ht]
 \begin{center}
{\small
\caption[Table 3]{
Example of interior angles of concave pentagons satisfying 
(\ref{eq3}) that can form the $n$-fold rotationally symmetric tilings
}
\label{tab3}
}
\begin{tabular}
{c| D{.}{.}{2} D{.}{.}{2} D{.}{.}{2} D{.}{.}{2} D{.}{.}{2} |c |c}
\hline
\raisebox{-1.50ex}[0cm][0cm]{$n$}& 
\multicolumn{5}{c|}{\shortstack{ Value of interior angle (degree) } } & 
\raisebox{-4.6ex}[0.7cm][0.5cm]{\shortstack{Edge \\length \\of $e$}} & 
\raisebox{-3.0ex}[0.7cm][0.5cm]{\shortstack{Figure \\number}} \\

 & 
\textit{A} & 
\textit{B}& 
\textit{C}& 
\textit{D}& 
\textit{E}& 
 & 
  \\
\hline
5& 
72& 
98& 
108& 
55.82 & 
206.18 & 
0.482 & 
\ref{fig23} \\
\hline
6& 
60& 
98& 
120& 
40.63 & 
221.37 & 
0.755 & 
\ref{fig24} \\
\hline
7& 
51.43 & 
106.41 & 
128.57 & 
39.90 & 
213.69 & 
1 & 
\ref{fig25} \\
\hline
8& 
45& 
112& 
135& 
36.64 & 
211.36 & 
1.174 & 
\ref{fig26} \\
\hline
9& 
40& 
112& 
140& 
31.63 & 
216.37 & 
1.271 & 
 \\
\hline
10& 
36& 
112& 
144& 
27.88 & 
220.12 & 
1.349 & 
\ref{fig27} \\
\hline
11& 
32.73 & 
112& 
147.27 & 
24.96 & 
223.04 & 
1.412 & 
 \\
\hline
12& 
30& 
135& 
150& 
28.16 & 
196.84 & 
1.608 & 
\ref{fig28} \\
\hline
13& 
27.69 & 
112& 
152.31 & 
20.67 & 
227.33 & 
1.509 & 
 \\
\hline
14& 
25.71 & 
112& 
154.29 & 
19.05 & 
228.95 & 
1.546 & 
 \\
\hline
15& 
24& 
112& 
156& 
17.66 & 
230.34 & 
1.578 & 
 \\
\hline
16& 
22.5& 
112& 
157.5& 
16.47 & 
231.53 & 
1.606 & 
 \\
\hline
17& 
21.18 & 
112& 
158.82 & 
15.43 & 
232.57 & 
1.631 & 
 \\
\hline
18& 
20& 
112& 
160& 
14.51 & 
233.49 & 
1.653 & 
 \\
\hline
...& 
...& 
...& 
...& 
...& 
...& 
...& 
 \\
\hline
\end{tabular}

\end{center}

\end{table}

Next, let us introduce some $n$-fold rotationally symmetric tilings with 
$C_{n}$ symmetry formed of concave pentagons that satisfies (\ref{eq3}) 
similar to the equilateral pentagonal case with $n = 7$. There are two cases 
of concave pentagons that satisfy (\ref{eq3}), $B > 180^ \circ$ and $E > 180^ \circ$.
A concave pentagon with $E > 180^ \circ$ is geometrically nonexistent if $n < 5$. 
Due to the presence of parameter $\theta $ in (\ref{eq3}), the shapes of concave 
pentagons that satisfy (\ref{eq3}) and can form an $n$-fold rotationally symmetric 
tiling are not fixed. Each example presented in Table~\ref{tab3} is such a concave 
pentagon with $E > 180^ \circ$. (For $n = 5\!-\!8, 10, 12$, 
tilings with concave pentagonal tiles are drawn. For further details, 
Figures~\ref{fig23}--\ref{fig28}.) Because the concave pentagon for the case 
of $n = 6$ also corresponds to a rhombus with an acute angle of $60^ \circ $, 
similar to that of a convex pentagon, various patterns can be generated by 
this concave pentagon (refer to Figure~\ref{fig29}).


\section{Rotationally symmetric tilings (tiling-like patterns) with an 
equilateral concave polygonal hole at the center}
\label{section5}

The rhombus can form various tilings, one of which is a rotationally symmetric 
tiling-like pattern with a regular polygonal hole at the center~\cite{H_and_M_2015}. 
Note that the tiling-like patterns are not considered tilings due to the presence 
of a gap, but are simply called tilings in this study. According to the properties 
deduced from \cite{H_and_M_2015}, pentagons satisfying (\ref{eq3}) can form 
rotationally symmetric tilings with a polygonal hole at the center, as shown 
in \cite{Sugimoto_2020_1} and \cite{Sugimoto_2020_2}. Though the rhombus 
has $D_{2}$ symmetry, the units comprising pentagons satisfying (\ref{eq3}) 
corresponding to the rhombus have $C_{2}$ symmetry. Therefore, pentagons 
satisfying (\ref{eq3}) can form rotationally symmetric tilings with an equilateral 
polygonal hole at the center, provided $n$ in (\ref{eq3}) is an even number. 
The hole formed at the center is an equilateral concave $2n$-gon with 
$D_{\frac{n}{2}} $ symmetry, and the tiling with hole has $C_{\frac{n}{2}} $ symmetry.

Let us introduce figures of these tilings. Figure~\ref{fig30} shows a rotationally 
symmetric tiling with $C_{4}$ symmetry, with an equilateral concave 16-gonal 
hole with $D_{4}$ symmetry at the center, using a convex pentagon with $n = 8$, 
as presented in Table~\ref{tab1}. Figure~\ref{fig31} shows a rotationally symmetric 
tiling with $C_{5}$ symmetry, with an equilateral concave 20-gonal hole with 
$D_{5}$ symmetry at the center, using a convex pentagon with $n = 10$, as 
presented in Table~\ref{tab1}. Figure~\ref{fig32} shows a rotationally symmetric 
tiling with $C_{8}$ symmetry, with an equilateral concave 32-gonal hole with $D_{8}$ 
symmetry at the center, using a convex pentagon with $n = 16$, as presented 
in Table~\ref{tab1}. As shown in these figures, the two types of rhombuses generated 
by pentagons (with and without gray color) are reflections of each other. 
That is, because these tilings are formed by alternately connecting the two types of 
rhombuses, they have $C_{\frac{n}{2}} $ symmetry and form an equilateral concave 
$2n$-gonal hole with $D_{\frac{n}{2}} $ symmetry with iterating concave and convex 
edges. According to the above properties, if $n$ in (\ref{eq3}) is an odd number, the 
polygonal holes cannot close. As described in Section~\ref{section4}, due to the presence 
of parameter $\theta $ in (\ref{eq3}), for pentagons satisfying the condition (\ref{eq3}), 
their shape is not fixed, and they need not be convex. Figure~\ref{fig33} shows a 
rotationally symmetric tiling with $C_{4}$ symmetry, with an equilateral 
concave 16-gonal hole with $D_{4}$ symmetry at the center, using a trapezoid 
with $n = 8$, as presented in Table~\ref{tab2}. Figure~\ref{fig34} shows a rotationally 
symmetric tiling with $C_{4}$ symmetry, with an equilateral concave 16-gonal 
hole with $D_{4}$ symmetry at the center, using a concave pentagon with $n = 8$, 
as presented in Table~\ref{tab3}. Figure~\ref{fig35} shows a rotationally symmetric 
tiling with $C_{5}$ symmetry, with an equilateral concave 20-gonal hole with 
$D_{5}$ symmetry at the center, using a concave pentagon with $n = 10$, as 
presented in Table~\ref{tab3}. Figure~\ref{fig36} shows a rotationally symmetric 
tiling with $C_{6}$ symmetry, with an equilateral concave 24-gonal hole with $D_{6}$ 
symmetry at the center, using a concave pentagon with $n = 12$, as presented 
in Table~\ref{tab3}.

Using pentagons satisfying (\ref{eq3}) with $n = 4$, similar to those shown in 
Figures~\ref{fig30}--\ref{fig36}, a rotationally symmetric tiling with $C_{2}$ 
symmetry, with an equilateral concave octagonal hole with $D_{2}$ symmetry 
at the center, can be formed. Because the concave octagonal hole corresponds to 
the shape of the pentagonal pair of Figure~\ref{fig04}(a), it can be filled with 
two pentagons.

Using pentagons satisfying (\ref{eq3}) with $n = 6$, corresponding to the 
rhombus with an acute angle of $60^ \circ $, similar to those shown in 
Figures~\ref{fig30}--\ref{fig36}, a rotationally symmetric tiling with $C_{3}$ 
symmetry, with an equilateral concave 12-gonal hole with $D_{3}$ symmetry 
at the center, can be formed. Because the concave 12-gonal hole corresponds to 
the shape shown in Figure~\ref{fig15}(a), it can be filled with six pentagons. 
Furthermore, this pentagon can form a three-fold rotational symmetric 
tiling as shown in Figure~\ref{fig06}. The outline of six pentagons at the center 
of such a tiling corresponds to an equilateral concave 12-gon shown in 
Figure~\ref{fig15}(a). Therefore, if the six pentagons at the center of such a 
tiling are removed, it appears as a three-fold rotationally symmetric tiling, with an 
equilateral concave 12-gonal hole with $D_{3}$ symmetry at the center. In the 
case of $n = 6$, as explained in Section~\ref{section4}, because the arrangement of 
pentagons inside the tilings can be replaced as shown in Figures~\ref{fig16} and 
\ref{fig17}, it can form different patterns with three-fold or six-fold rotational symmetry, 
or patterns without rotational symmetry. The above patterns of tilings with 
an equilateral concave 12-gonal hole at the center by the pentagons with 
$n = 6$ are one such variation.

The above-mentioned rotationally symmetric tiling with a regular polygonal 
hole at the center, using rhombuses, is formed by the following method. 
Because one inner angle (interior angle) of a regular $m$-gon is 
``$180^ \circ -\frac{360^ \circ }{m}$," the outer angle of one vertex of a regular 
$m$-gon is ``$180^ \circ + \frac{360^ \circ }{m}$," and that can be achieved by 
a combination of the acute and obtuse angles of the rhombus. For example, in 
the case of a regular octagon ($m = 8$), the interior angle of one vertex is 
$135^ \circ $, so the value of ``$360^ \circ - 135^ \circ = 225^ \circ $" will 
be shared by one obtuse and multiple acute angles. This sharing can be done in 
rhombuses with an acute angle of $\frac{360^ \circ }{8k}$, where $k$ is 
an integer greater than or equal to one. For a rhombus with an acute angle 
of $45^ \circ $ (when $k = 1$), sharing will be ``$2\times 45^ \circ + 135^ \circ 
= 225^ \circ $"; for a rhombus with an acute angle of $22.5^ \circ $ (when $k = 2$), 
sharing will be ``$3\times 22.5^ \circ + 157.5{ }^ \circ = 225^ \circ $" and so on. 
In fact, a rhombus with an acute angle of $\frac{360^ \circ }{m \cdot k}$ 
can form a rotationally symmetric tiling with a regular $m$-gonal hole at the 
center \cite{H_and_M_2015}. Therefore, a pentagon satisfying (\ref{eq3}) with 
$n = m \cdot k$ is a candidate for forming a rotationally symmetric tiling with 
$C_{\frac{m}{2}} $ symmetry, with an equilateral concave $2m$-gon hole with 
$D_{\frac{m}{2}} $ symmetry at the center. This may be geometrically established 
depending on how parameter $\theta $ is selected, and it is possible provided 
$n$ is an even number, as described above.

For example, if the pentagons can form rotationally symmetric tilings with 
$C_{4}$ symmetry, with an equilateral concave 16-gonal hole with $D_{4}$ 
symmetry at the center, they correspond to pentagons satisfying (\ref{eq3}) 
whose $n$ is a multiple of eight. Figure~\ref{fig30} is a case of $k = 1$ and $n = 8$. 
Figure~\ref{fig37} is a case of $k = 2$, i.e., a rotationally symmetric tiling with 
$C_{4}$ symmetry, with an equilateral concave 16-gonal hole with $D_{4}$ 
symmetry at the center, by a convex pentagon with $n = 16$, as presented in 
Table~\ref{tab1}. Similarly, the concave pentagon with $n = 12$, as presented in 
Table~\ref{tab3}, can form a rotationally symmetric tiling with $C_{3}$ symmetry, 
with an equilateral concave 12-gonal hole with $D_{3}$ symmetry at the center, as 
shown in Figure~\ref{fig38}.


\section{Tilings with multiple pentagons of different shapes}
\label{section6}

Rhombuses can form edge-to-edge tilings by using different shapes of 
rhombuses with differing interior angles when the lengths of edges are same. 
Tilings using two or more types of pentagons that satisfy (\ref{eq2}), in which 
$\theta $ has same value and $\alpha $ has different values, correspond to 
the tilings with two or more types of rhombuses. That is, the tiling is not 
monohedral. Figure~\ref{fig39}(a) is an example of tiling using convex pentagons 
with $n = 3, 4, 8$, as presented in Table~\ref{tab1}, and Figure~\ref{fig39}(b) 
is an example of tiling using concave pentagons with $n = 8, 10$, as presented 
in Table~\ref{tab3}. In these examples, pentagons satisfying (\ref{eq3}) 
are used, but we note that tilings can be formed by pentagons satisfying 
(\ref{eq2}), whose interior angle $A$ is not $\frac{360^ \circ }{n}$. In the case of 
trapezoids, presented in Table~\ref{tab2}, it is clear that tiling is formed by 
multiple different trapezoids.

Furthermore, pentagons satisfying (\ref{eq2}) with the same value of $\theta $ 
can be used in a tiling, whether convex, concave, or trapezoidal. 
For $\theta = 45^ \circ $, Figure~\ref{fig40} shows an example of tiling by 
convex pentagons satisfying (\ref{eq2}) with $\alpha = 54^ \circ $, trapezoids 
(pentagons) satisfying (\ref{eq2}) with $\alpha = 67.5^ \circ $, and concave 
pentagons satisfying (\ref{eq2}) with $\alpha = 75^ \circ $.

The pentagonal tilings in Figures~\ref{fig39} and \ref{fig40} satisfying 
``$B + D + E = 360^ \circ ,\;2A + 2C = 360^ \circ $" are rhombic tilings that are 
formed from rhombic belts made by translation in the same direction. By adjusting 
the combination of rhombuses used, tilings other than the combination of above 
belts can be formed. (They correspond to pentagonal tilings that admit the 
vertex  concentrations ``$B + D + E = 360^ \circ ,\;2A + 2C = 360^ \circ $" and 
also vertex concentrations other than ``$B + D + E = 360^ \circ ,\;2A + 2C = 360^ \circ $.") 
For example, because squares and rhombuses with an acute angle of $45^ \circ $ 
can form an eight-fold rotationally symmetric tiling \cite{H_and_M_2015}, a pentagonal 
tiling, as shown in Figure~\ref{fig41}, corresponding to it can be formed by 
convex pentagons with $n = 4, 8$, as presented in Table~\ref{tab1}. In addition, 
because rhombuses with acute angles of $72^ \circ $ and $36^ \circ $ can form 
a five-fold rotationally symmetric tiling, a pentagonal tiling \cite{H_and_M_2015}, as shown in 
Figure~\ref{fig42}, corresponding to it can be formed by convex pentagons 
satisfying (\ref{eq3}) with $n = 5$ and $\theta = 45^ \circ $, and concave pentagons 
satisfying (\ref{eq3}) with $n = 10$ and $\theta = 45^ \circ $. Note that the number 
of pentagons satisfying (\ref{eq3}) included in the corresponding rhombuses can be 
changed. (In Figure~\ref{fig41}, one rhombus includes 32 pentagons, and in 
Figure~\ref{fig42}, one rhombus includes eight pentagons. It is possible to 
have a case where one rhombus includes two pentagons or $8\times u^2$ 
pentagons where $u = 1, 2, 3, \ldots $.) Similarly, tilings with three or more 
types of rhombuses can be converted into tilings with pentagons.


\section{Conclusions}
\label{section7}

In \cite{Sugimoto_2020_1} and \cite{Sugimoto_2020_2}, we introduced convex 
pentagonal tiles, belonging to the Type 1 family, that can generate countless 
rotationally symmetric tilings. 
In this study, we have shown that convex pentagonal tiles belonging 
to the Type 2 family can generate countless rotational symmetric tilings. In 
addition, because the pentagons have two parameters, the study discussed 
that the tilings can be generated by shapes other than convex. 

Because the properties of pentagons dealt with in this study correspond to 
those of rhombuses, it also explained the correspondence between pentagons 
and various rhombic tilings. Not all rhombic tilings (including tilings with 
holes as introduced in Section~\ref{section5}) can be converted into pentagonal 
tilings by the method discussed in this study. But, various knowledge of rhombic 
tilings can be used to generate various pentagonal tilings.

Livio Zucca presented interesting tilings using equilateral pentagons in 
\cite{Zucca_2003}. However, that study does not consider pentagons with
four equal-length edges and the relationship between pentagons and rhombuses.

\bigskip
\noindent
\textbf{Acknowledgments.} 
The author would like to thank Yoshiaki ARAKI of Japan Tessellation Design 
Association, for discussions and comments.



\bigskip
\appendix
\def\thesection{Appendix}
\section{}
\label{app}

We organized the properties of the shape of pentagons that satisfy (\ref{eq3}) 
depending on the values of $n$ and $\theta $. If the pentagons that satisfy (\ref{eq3}) 
exist geometrically, then $n > 2$, because $C > 0^ \circ $. It is $0^ \circ < 
\theta < 180^ \circ $ because the edge $e$ of the pentagon that satisfies (\ref{eq3}) 
exists and does not intersect the edge $b$. As mentioned in Section~\ref{section4}, the 
pentagons that satisfy (\ref{eq3}) with $n = 3$ and $n = 6$ correspond to rhombuses 
with an acute angle of $60^ \circ $, and these pentagons are opposite to 
each other. The shapes of the pentagons that satisfy (\ref{eq3}) change depending on 
the values of  $n$ and $\theta $ are as follows:

\begin{enumerate}
\setlength{\itemindent}{0pt}
\item[]

\textbf{Case where the pentagons that satisfy (\ref{eq3}) admit $n = 4$}

\begin{itemize}
	\item $0^ \circ < \theta < 90^ \circ$ : Convex pentagons
	\item $\theta = 90^ \circ$ : Parallelograms ($B = 180^ \circ$)
	\item $90^ \circ < \theta < 180^ \circ$ : Concave pentagons ($B > 180^ \circ$)
\end{itemize}

\end{enumerate}

\vspace{1pt}

\begin{enumerate}
\setlength{\itemindent}{0pt}
\item[]

\textbf{Cases where the pentagons that satisfy (\ref{eq3}) admit $n \ge 5$}

\begin{itemize}
	\item $0^ \circ < \theta < 90^ \circ - \dfrac{360^ \circ }{n}$ : Concave pentagons 
($E > 180^ \circ$)

	\item $\theta = 90^ \circ - \dfrac{360^ \circ }{n}$ : Trapezoids ($E = 180^ \circ$)

	\item $90^ \circ - \dfrac{360^ \circ }{n} < \theta < 90^ \circ$ : Convex pentagons

	\item $\theta = 90^ \circ$ : Parallelograms ($B = 180^ \circ$)

	\item $90^ \circ < \theta < 180^ \circ$ : Concave pentagons ($B > 180^ \circ$)
\end{itemize}

\end{enumerate}

\vspace{1pt}

\begin{enumerate}
\setlength{\itemindent}{0pt}
\item[]

\textbf{Cases where the pentagons that satisfy (\ref{eq3}) become the polygons with 
$a = b = c = d = e$}

\begin{itemize}
	\item Case where $n = 4$ and $\theta \approx 41.41^ \circ $ (Convex pentagon with 
$E \approx 114.30^ \circ$)

	\item Case where $n = 5$ and $\theta \approx 37.95^ \circ $ (Convex pentagon with 
$E \approx 149.76^ \circ$)

	\item Case where $n = 6$ and $\theta = 30^ \circ $ (Trapezoid of $E = 180^ \circ$)

	\item Case where $n = 7$ and $\theta \approx 16.41^ \circ $ (Concave pentagon with 
$E \approx 213.69^ \circ$)
\end{itemize}

\end{enumerate}

If the pentagons that satisfy (\ref{eq3}) have $\theta = 90^ \circ $ (i.e., 
$B = 180^ \circ$), then they are parallelograms (the case of $n = 4$ is a 
rectangle). Such parallelograms correspond to half of the rhombus 
corresponding to the basic unit of Figures~\ref{fig04}(a) or \ref{fig04}(b).

If the pentagons that satisfy (\ref{eq3}) have $E = 180^ \circ $, then 
``$A = D,\;B = C$.'' Thus, they are trapezoids with a line of symmetry and 
$\theta = 90^ \circ - \frac{360^ \circ }{n}$ holds. It is when $n \ge 5$ that the pentagons 
that satisfy (\ref{eq3}) become trapezoids, because $\theta > 0^ \circ $.

If the pentagons that satisfy (\ref{eq3}) have $a = b = c = d = e = 1$, then 
$1 = 2\sqrt {1 - \sin (2\alpha )\cos \theta } $ where $\alpha = 90^  \circ-\frac{180^ \circ }{n}$ 
(refer to Section~\ref{section2}). In this case, $n \le 7$ for $\theta $ to exist.

Finally, example figures using a concave pentagon with $B > 180^ \circ$ are shown.
Figure~\ref{fig43} shows an eight-fold rotationally symmetric tiling by a concave 
pentagon satisfying (\ref{eq3}) with $n=8$ and $B = 224^ \circ$. Figure~\ref{fig44} shows 
a rotationally symmetric tiling with $C_{4}$ symmetry, with an equilateral concave 
16-gonal hole with $D_{4}$ symmetry at the center, using a convex pentagon 
satisfying (\ref{eq3}) with $n = 8$ and $B =224^ \circ$.


\renewcommand{\figurename}{{\small Figure.}}
\begin{figure}[htbp]
 \centering\includegraphics[width=15cm,clip]{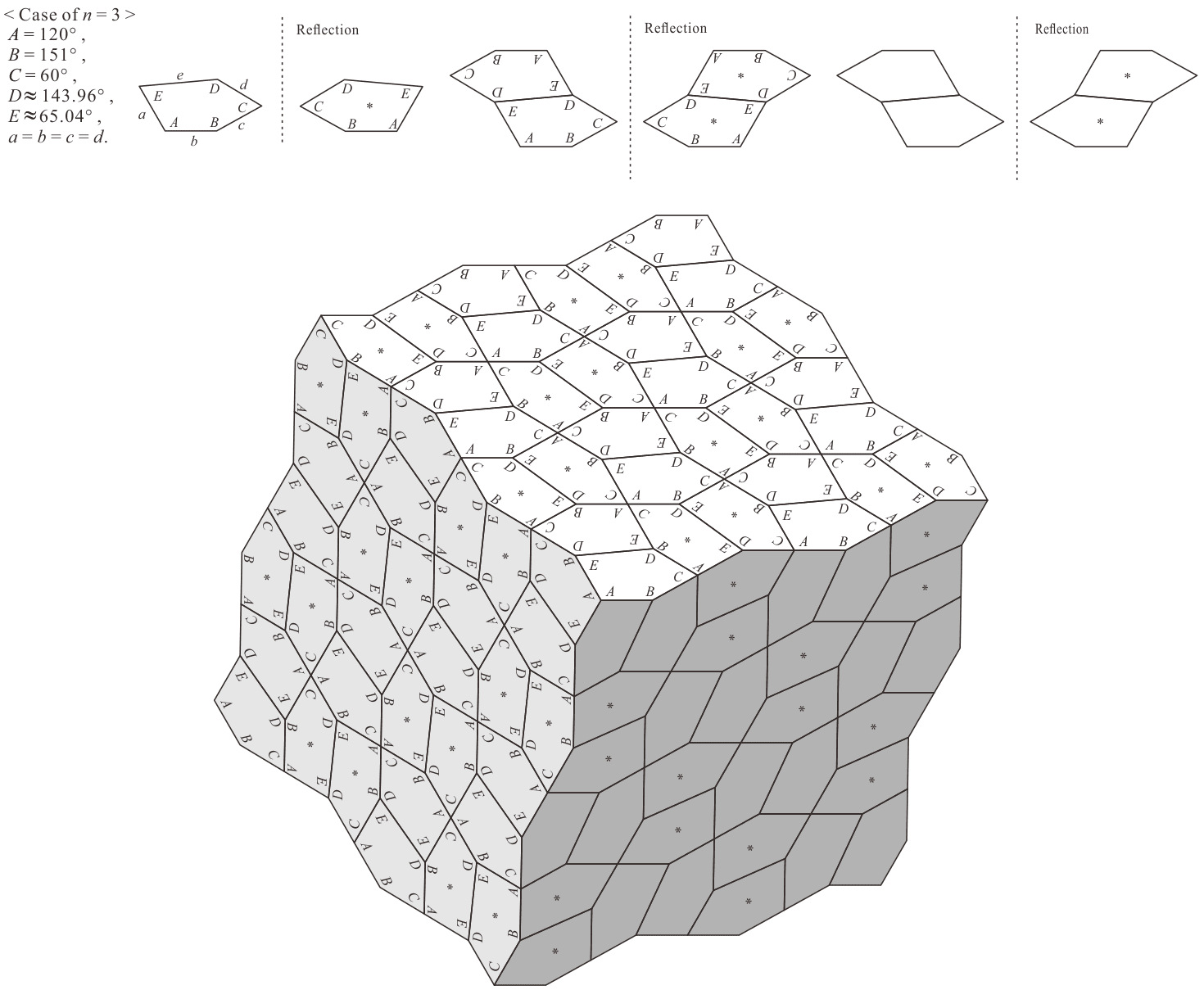} 
  \caption{{\small 
Three-fold rotationally symmetric tiling by a convex 
pentagon of $n=3$ in Table~\ref{tab1} 
 (The figure is solely a depiction of  the area around the rotationally symmetric
 center, and the tiling can be spread in all directions as well)
} 
\label{fig06}
}
\end{figure}

\renewcommand{\figurename}{{\small Figure.}}
\begin{figure}[htbp]
 \centering\includegraphics[width=15cm,clip]{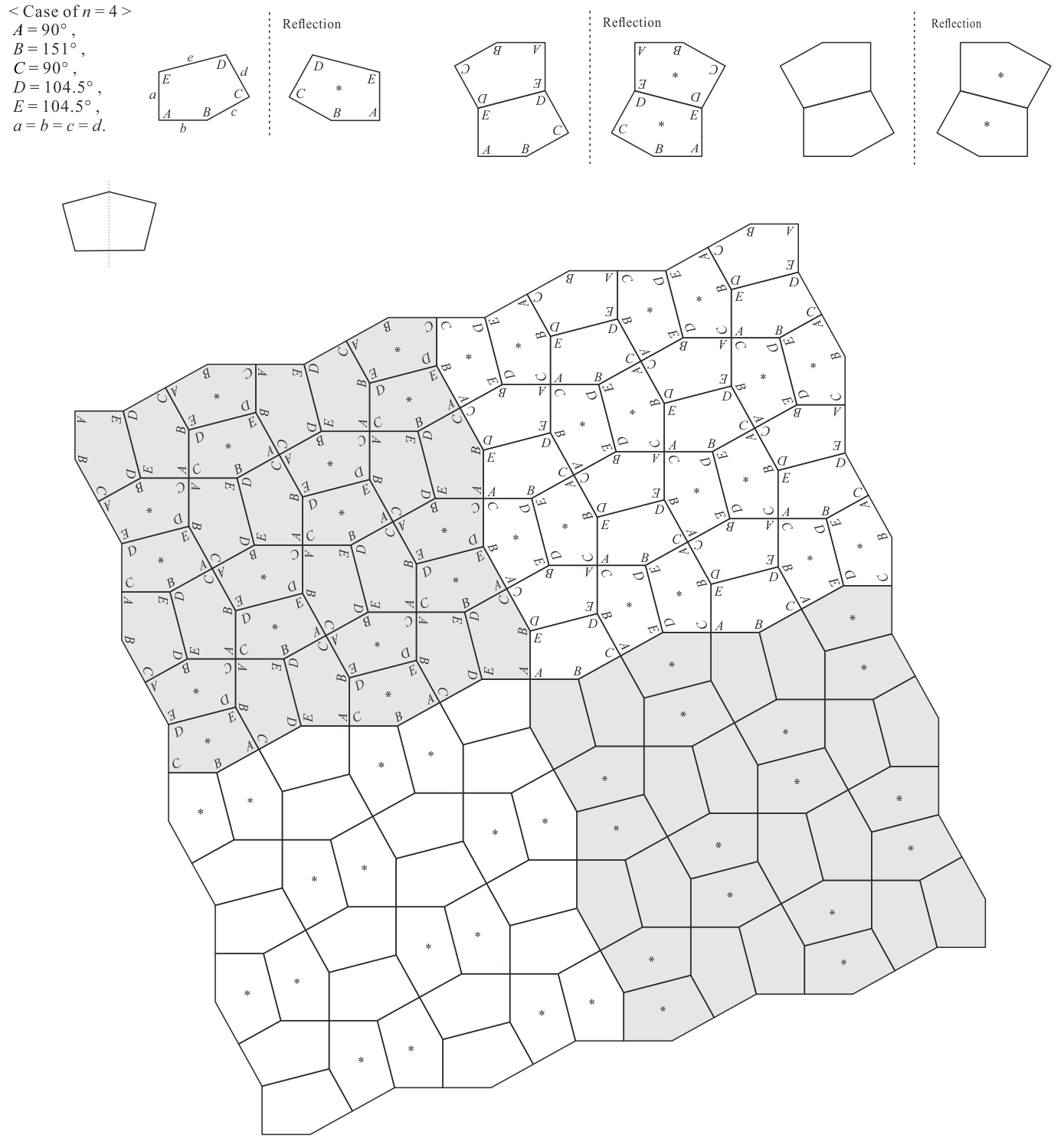} 
  \caption{{\small 
Four-fold rotationally symmetric tiling by a convex 
pentagon of $n=4$ in Table~\ref{tab1} 
 (The figure is solely a depiction of  the area around the rotationally symmetric
 center, and the tiling can be spread in all directions as well)
} 
\label{fig07}
}
\end{figure}

\renewcommand{\figurename}{{\small Figure.}}
\begin{figure}[htbp]
 \centering\includegraphics[width=15cm,clip]{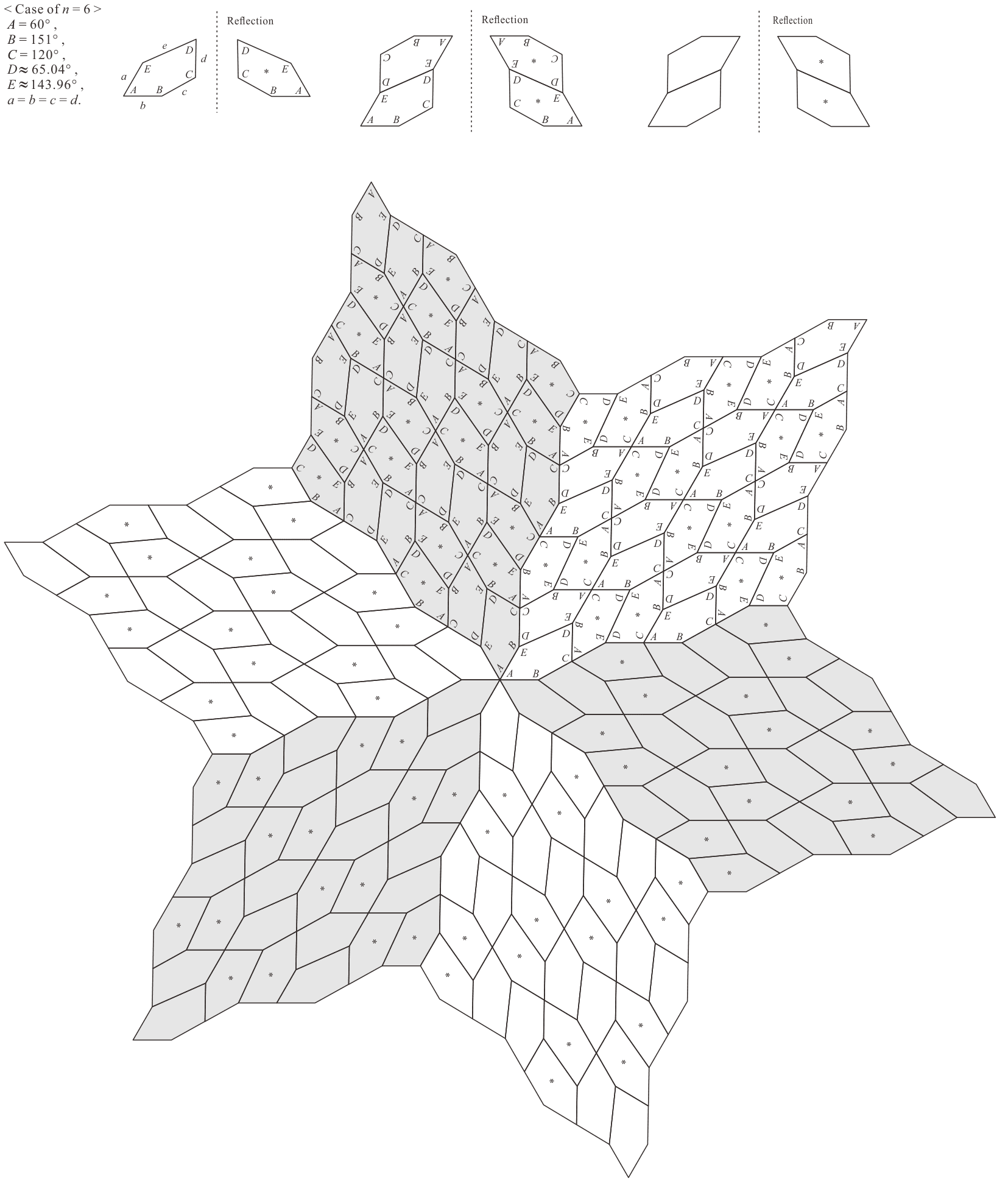} 
  \caption{{\small 
Six-fold rotationally symmetric tiling by a convex 
pentagon of $n=6$ in Table~\ref{tab1} 
 (The figure is solely a depiction of  the area around the rotationally symmetric
 center, and the tiling can be spread in all directions as well)
} 
\label{fig08}
}
\end{figure}

\renewcommand{\figurename}{{\small Figure.}}
\begin{figure}[htbp]
 \centering\includegraphics[width=15cm,clip]{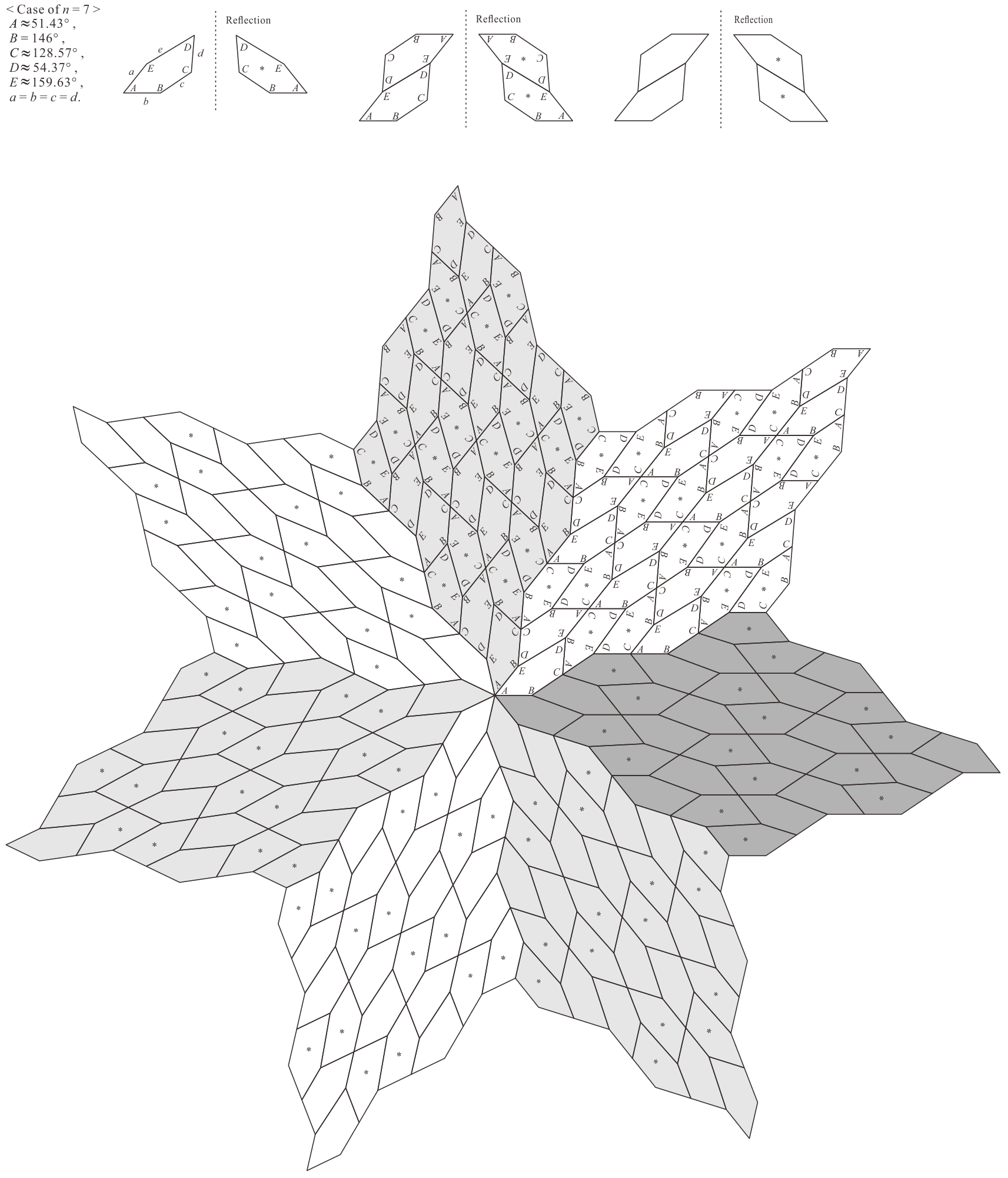} 
  \caption{{\small 
Seven-fold rotationally symmetric tiling by a convex 
pentagon of $n=7$ in Table~\ref{tab1} 
 (The figure is solely a depiction of  the area around the rotationally symmetric
 center, and the tiling can be spread in all directions as well)
} 
\label{fig09}
}
\end{figure}

\renewcommand{\figurename}{{\small Figure.}}
\begin{figure}[htbp]
 \centering\includegraphics[width=15cm,clip]{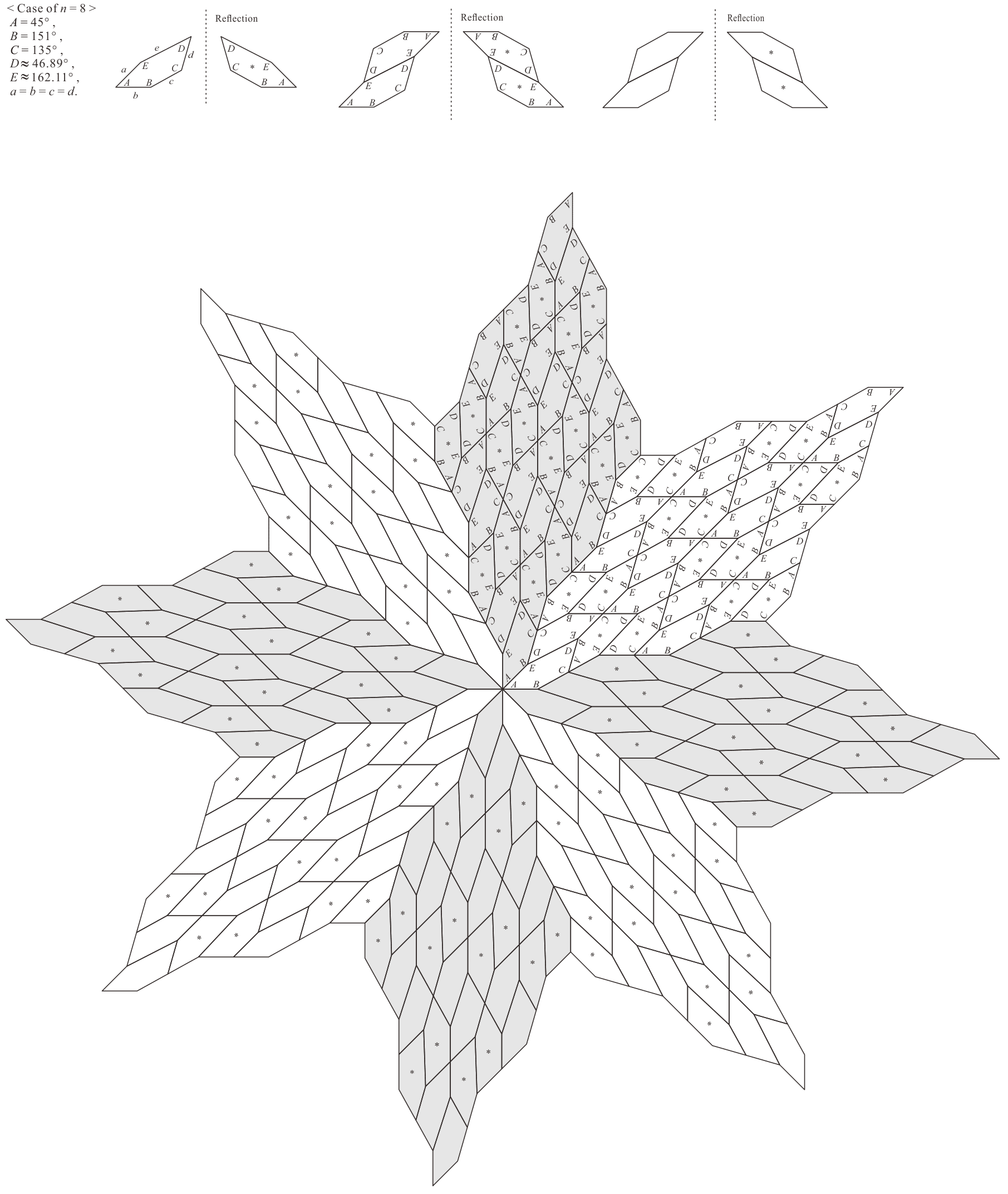} 
  \caption{{\small 
Eight-fold rotationally symmetric tiling by a convex 
pentagon of $n=8$ in Table~\ref{tab1} 
 (The figure is solely a depiction of  the area around the rotationally symmetric
 center, and the tiling can be spread in all directions as well)
}  
\label{fig10}
}
\end{figure}

\renewcommand{\figurename}{{\small Figure.}}
\begin{figure}[htbp]
 \centering\includegraphics[width=15cm,clip]{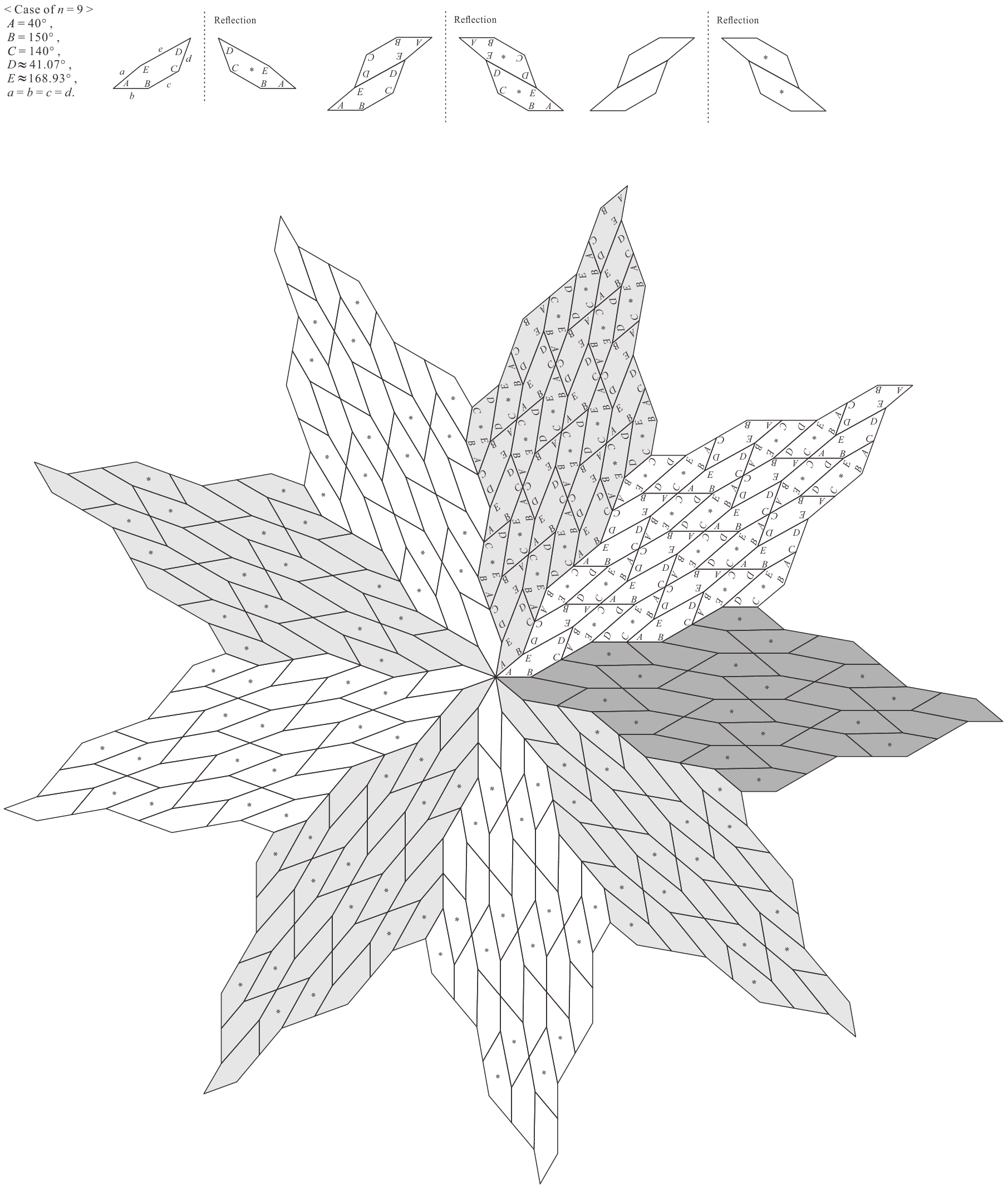} 
  \caption{{\small 
Nine-fold rotationally symmetric tiling by a convex 
pentagon of $n=9$ in Table~\ref{tab1} 
 (The figure is solely a depiction of  the area around the rotationally symmetric
 center, and the tiling can be spread in all directions as well)
} 
\label{fig11}
}
\end{figure}

\renewcommand{\figurename}{{\small Figure.}}
\begin{figure}[htbp]
 \centering\includegraphics[width=15cm,clip]{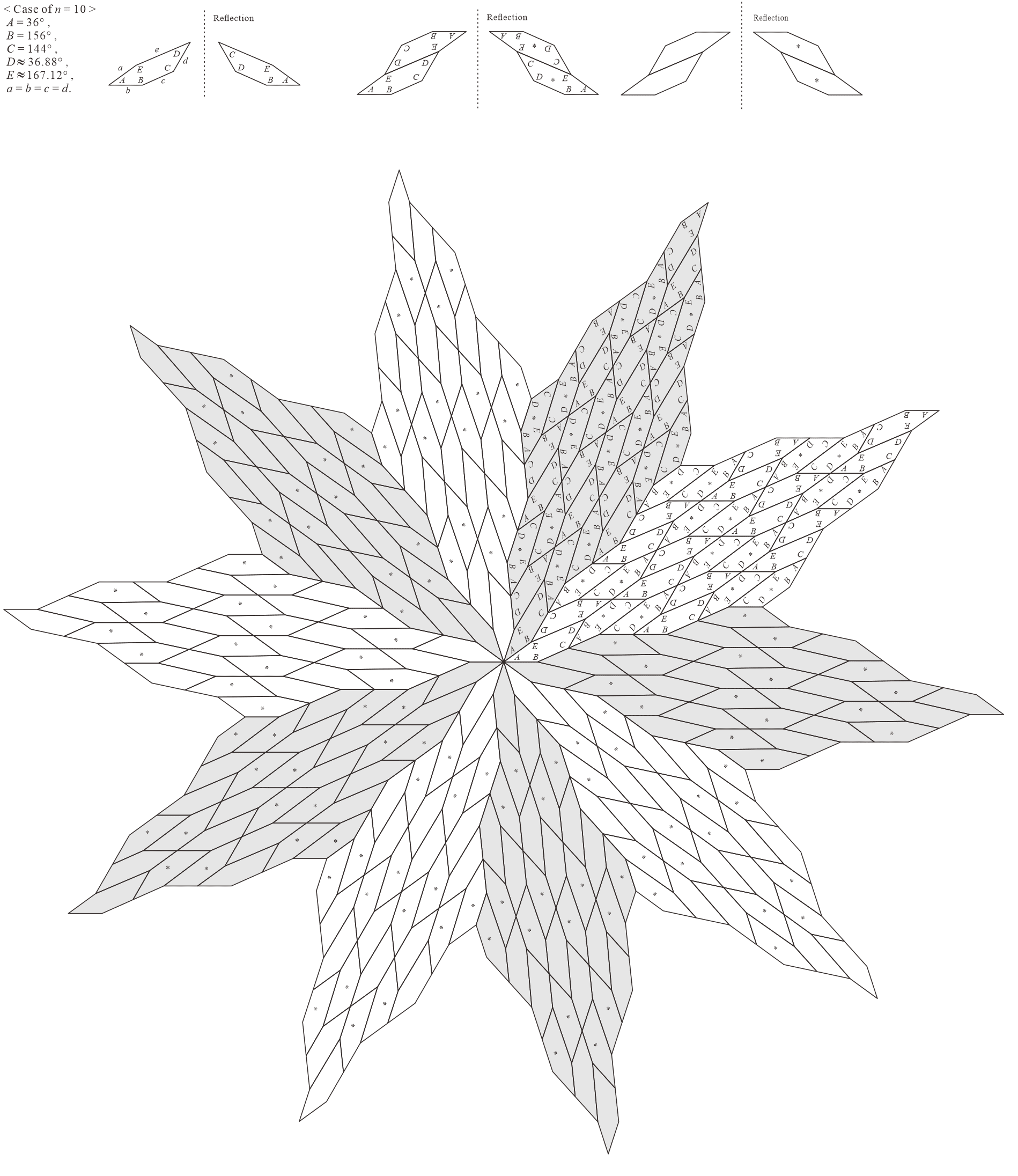} 
  \caption{{\small 
10-fold rotationally symmetric tiling by a convex 
pentagon of $n=10$ in Table~\ref{tab1} 
 (The figure is solely a depiction of  the area around the rotationally symmetric
 center, and the tiling can be spread in all directions as well)
}
\label{fig12}
}
\end{figure}

\renewcommand{\figurename}{{\small Figure.}}
\begin{figure}[htbp]
 \centering\includegraphics[width=15cm,clip]{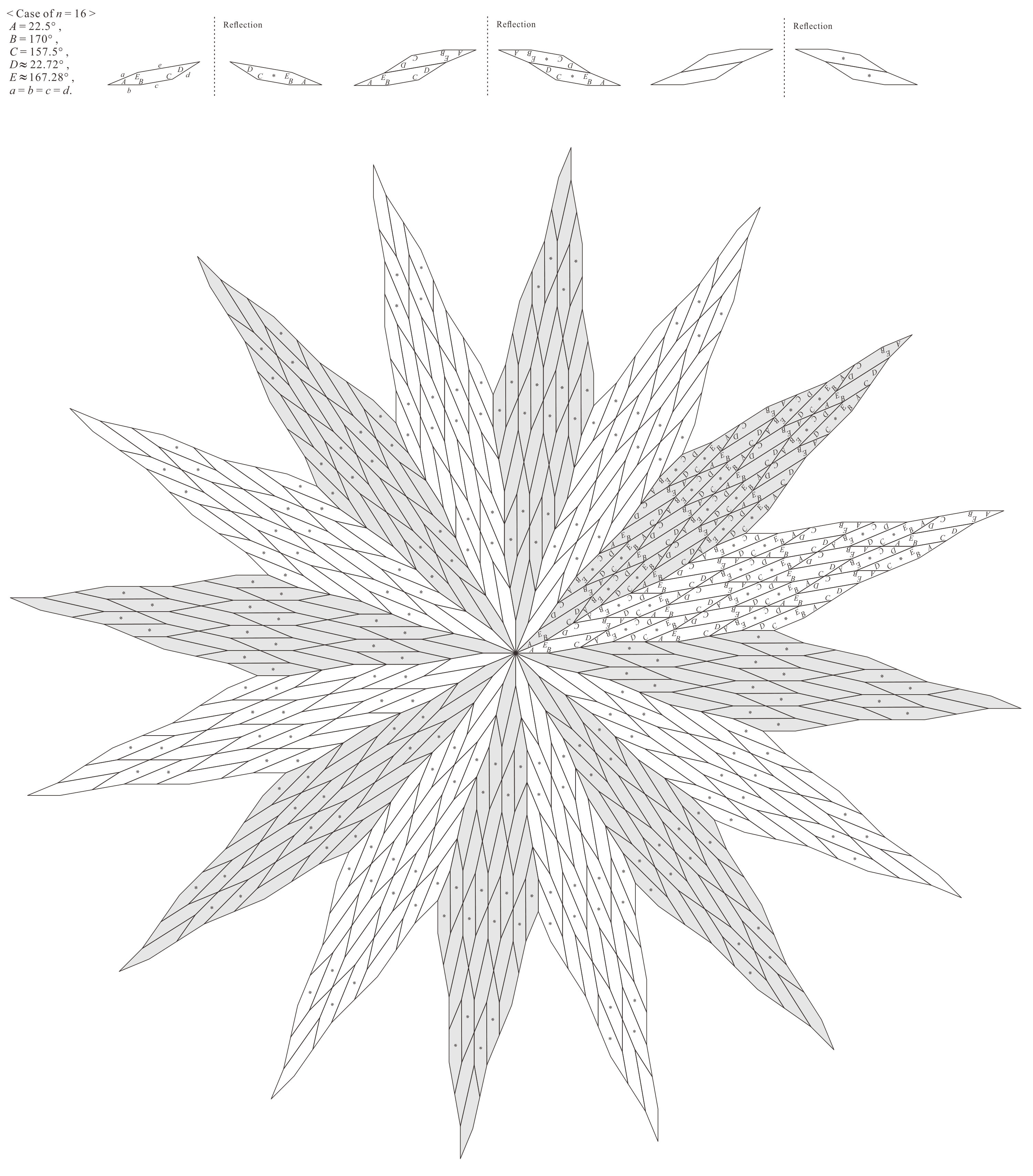} 
  \caption{{\small 
16-fold rotationally symmetric tiling by a convex 
pentagon of $n=16$ in Table~\ref{tab1} 
 (The figure is solely a depiction of  the area around the rotationally symmetric
 center, and the tiling can be spread in all directions as well)
} 
\label{fig13}
}
\end{figure}

\renewcommand{\figurename}{{\small Figure.}}
\begin{figure}[htbp]
 \centering\includegraphics[width=15cm,clip]{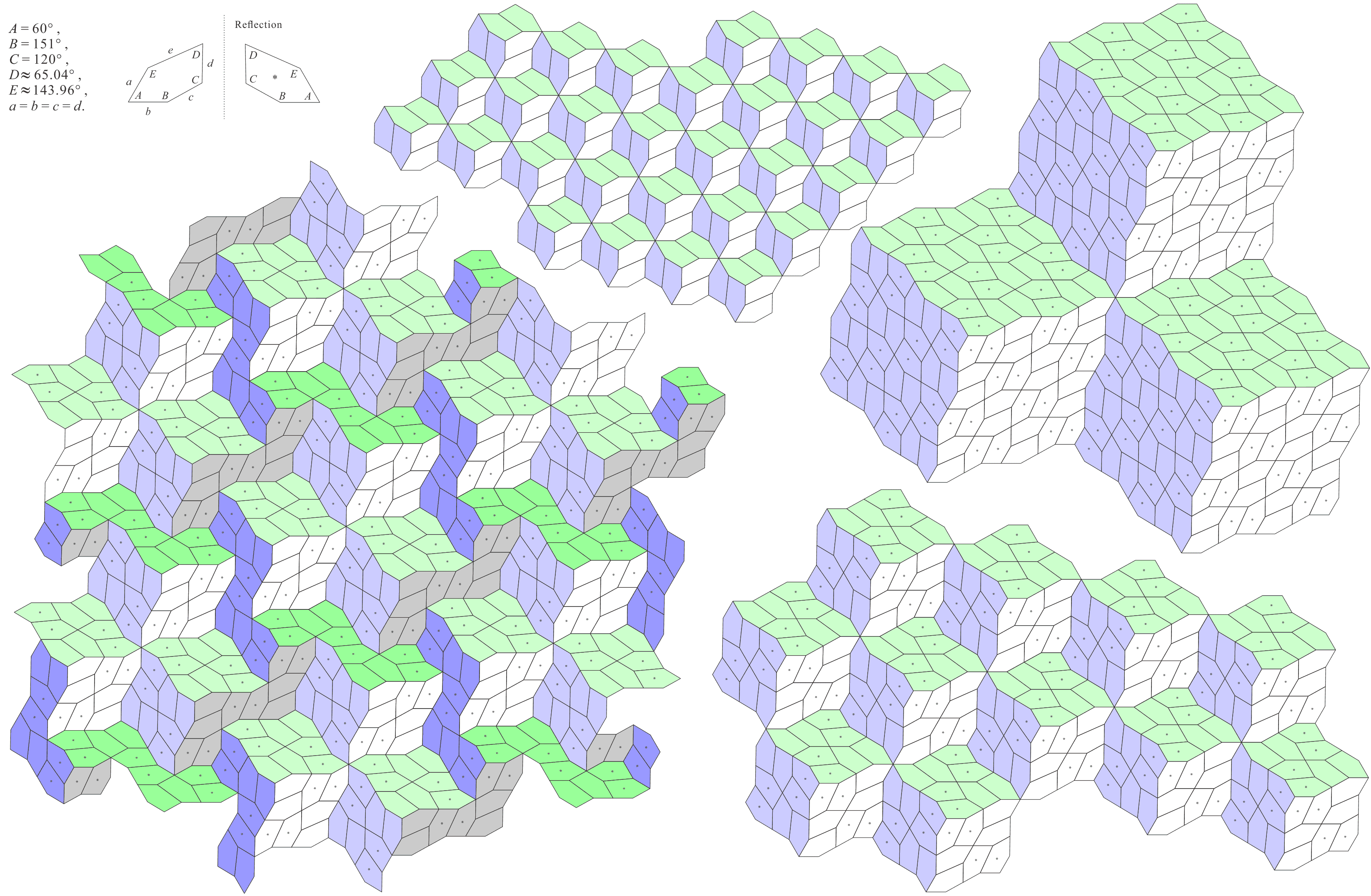} 
  \caption{{\small 
Examples of tilings with three-fold and six-fold rotational symmetry by 
a pentagon that corresponds to rhombus with an acute angle 
of $60^ \circ$ } 
\label{fig14}
}
\end{figure}

\renewcommand{\figurename}{{\small Figure.}}
\begin{figure}[htbp]
 \centering\includegraphics[width=14cm,clip]{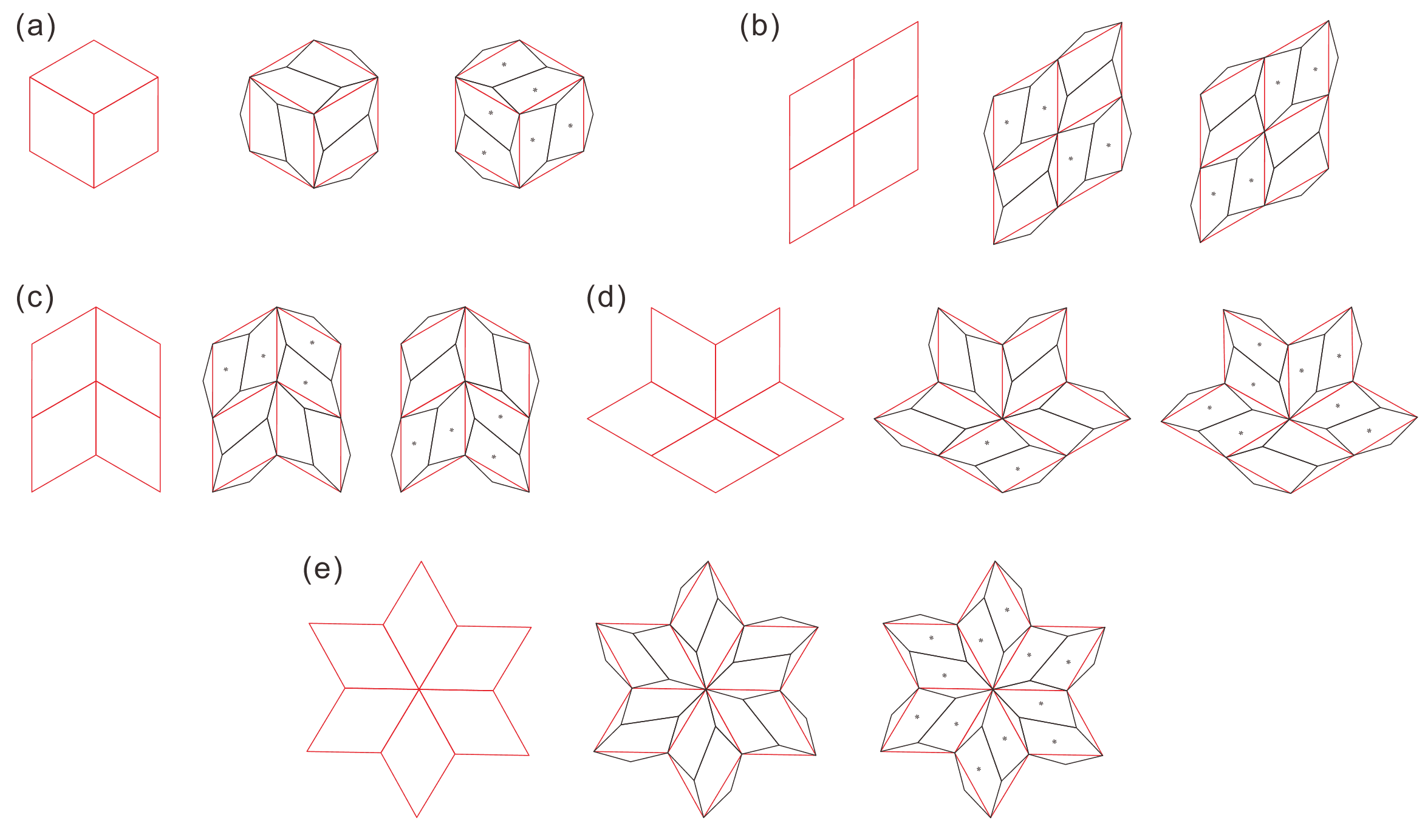} 
  \caption{{\small 
Combinations of vertices $A $ and $C$ of convex pentagons that 
correspond to rhombuses with an acute angle of $60^ \circ$} 
\label{fig15}
}
\end{figure}

\renewcommand{\figurename}{{\small Figure.}}
\begin{figure}[htbp]
 \centering\includegraphics[width=15cm,clip]{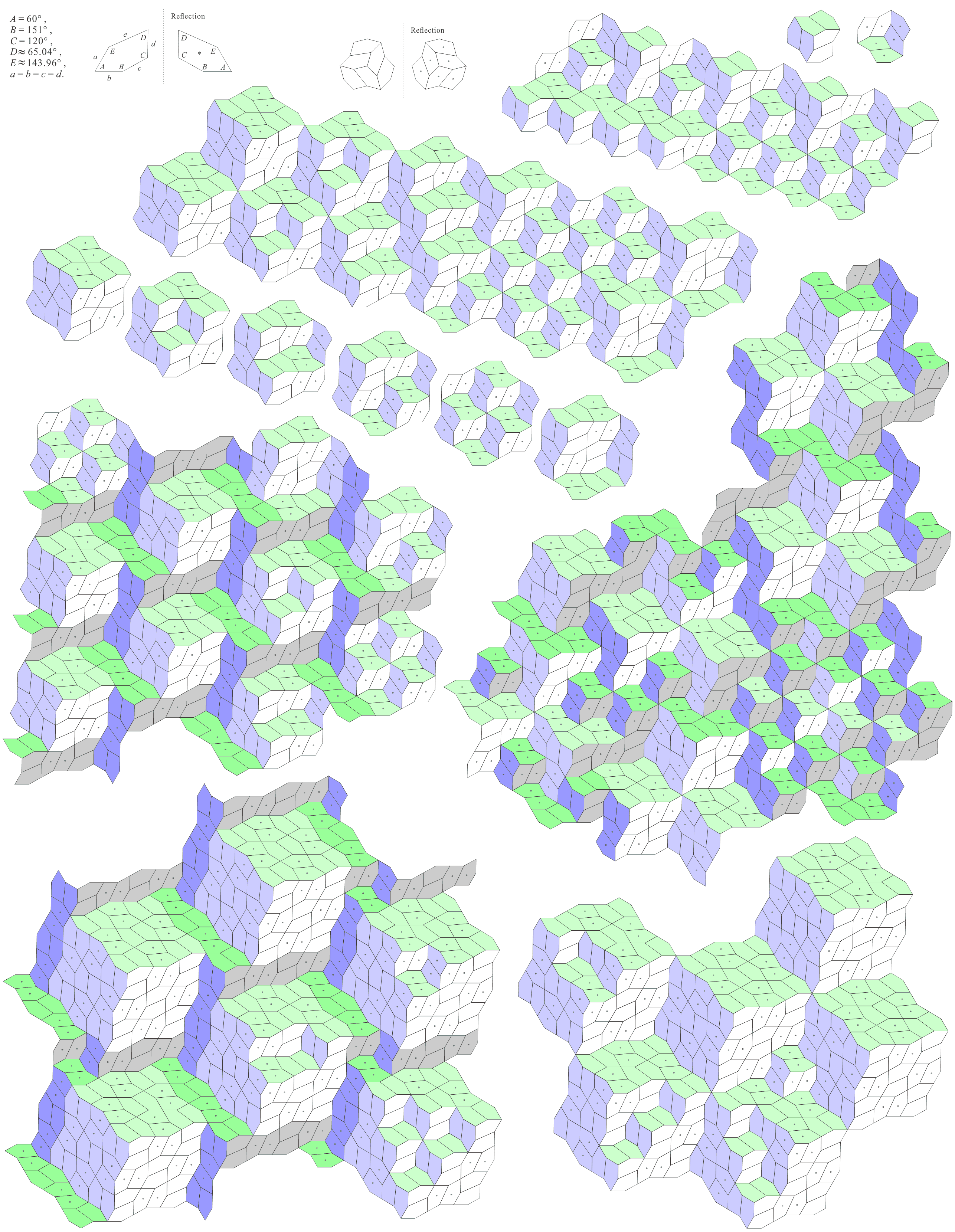} 
  \caption{{\small 
Examples of tilings by a pentagon that corresponds to rhombus 
with an acute angle of $60^ \circ$, Part 1 } 
\label{fig16}
}
\end{figure}

\renewcommand{\figurename}{{\small Figure.}}
\begin{figure}[htbp]
 \centering\includegraphics[width=15cm,clip]{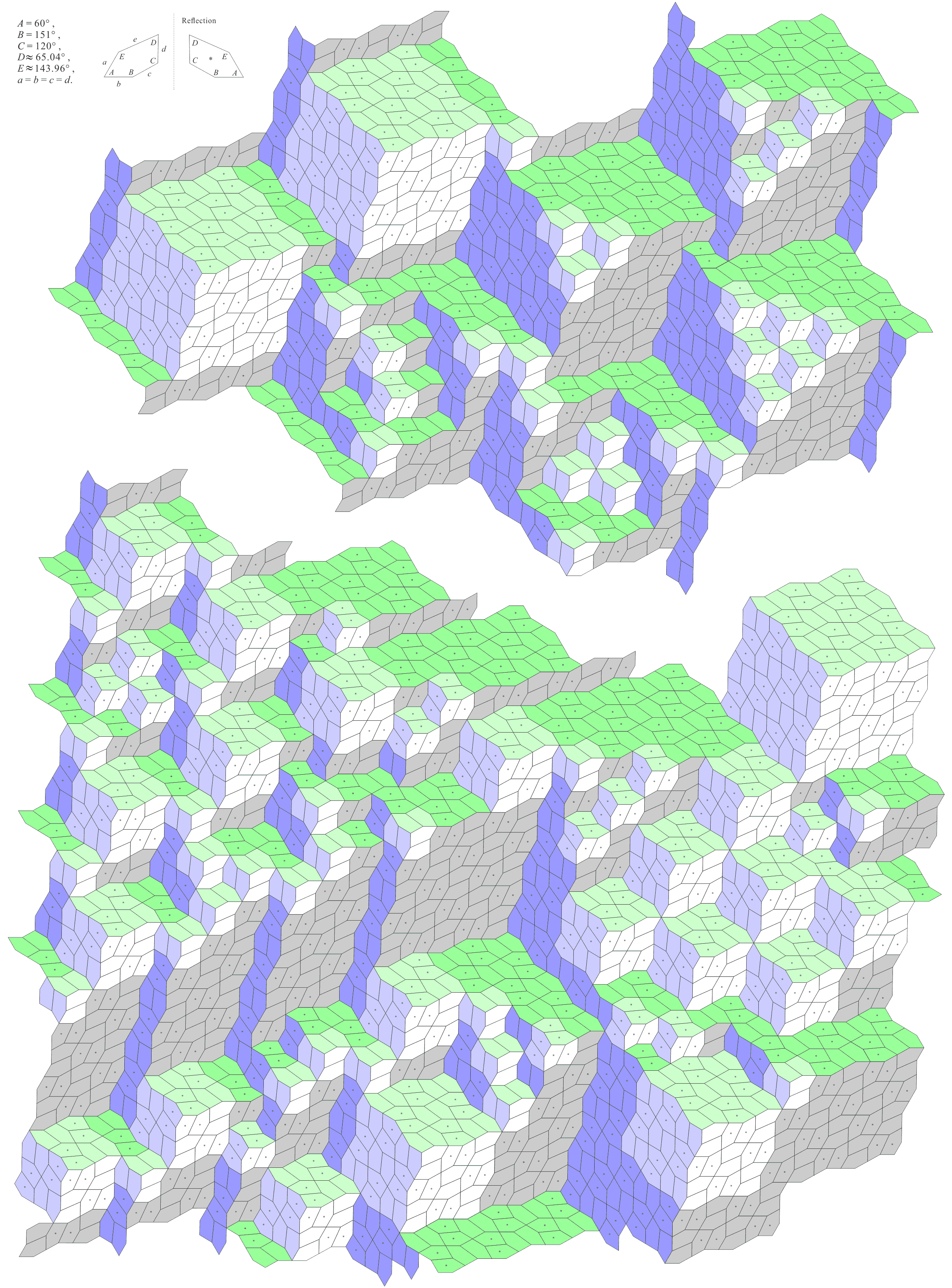} 
  \caption{{\small 
Examples of tilings by a pentagon that corresponds to rhombus 
with an acute angle of $60^ \circ$, Part 2 } 
\label{fig17}
}
\end{figure}

\renewcommand{\figurename}{{\small Figure.}}
\begin{figure}[htbp]
 \centering\includegraphics[width=15cm,clip]{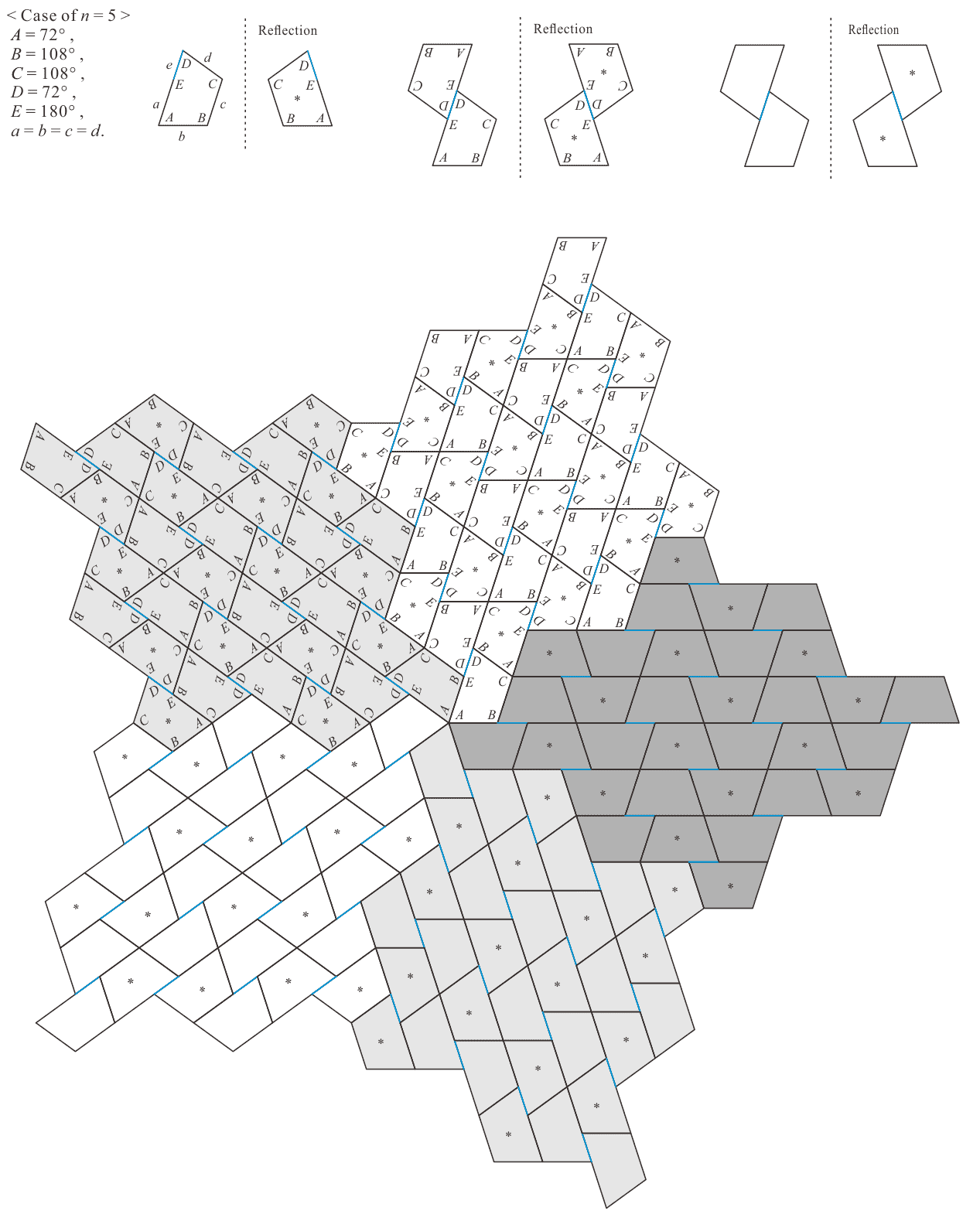} 
  \caption{{\small 
Five-fold rotationally symmetric tiling by a trapezoid of $n=5$ 
in Table~\ref{tab2} 
 (The figure is solely a depiction of  the area around the rotationally symmetric
 center, and the tiling can be spread in all directions as well)
} 
\label{fig18}
}
\end{figure}

\renewcommand{\figurename}{{\small Figure.}}
\begin{figure}[htbp]
 \centering\includegraphics[width=15cm,clip]{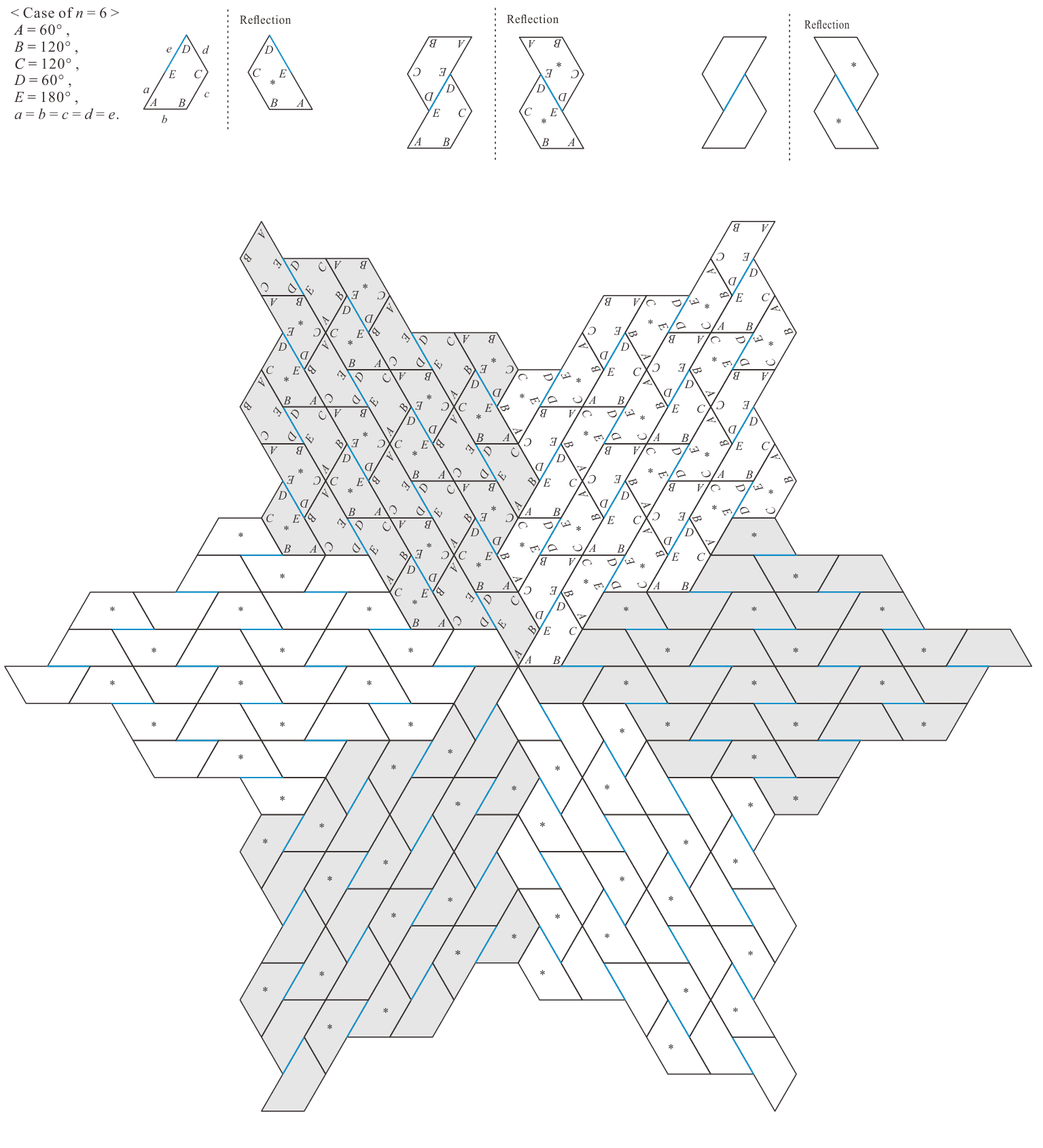} 
  \caption{{\small 
Six-fold rotationally symmetric tiling by a trapezoid of $n=6$ 
in Table~\ref{tab2} 
 (The figure is solely a depiction of  the area around the rotationally symmetric
 center, and the tiling can be spread in all directions as well)
} 
\label{fig19}
}
\end{figure}

\renewcommand{\figurename}{{\small Figure.}}
\begin{figure}[htbp]
 \centering\includegraphics[width=15cm,clip]{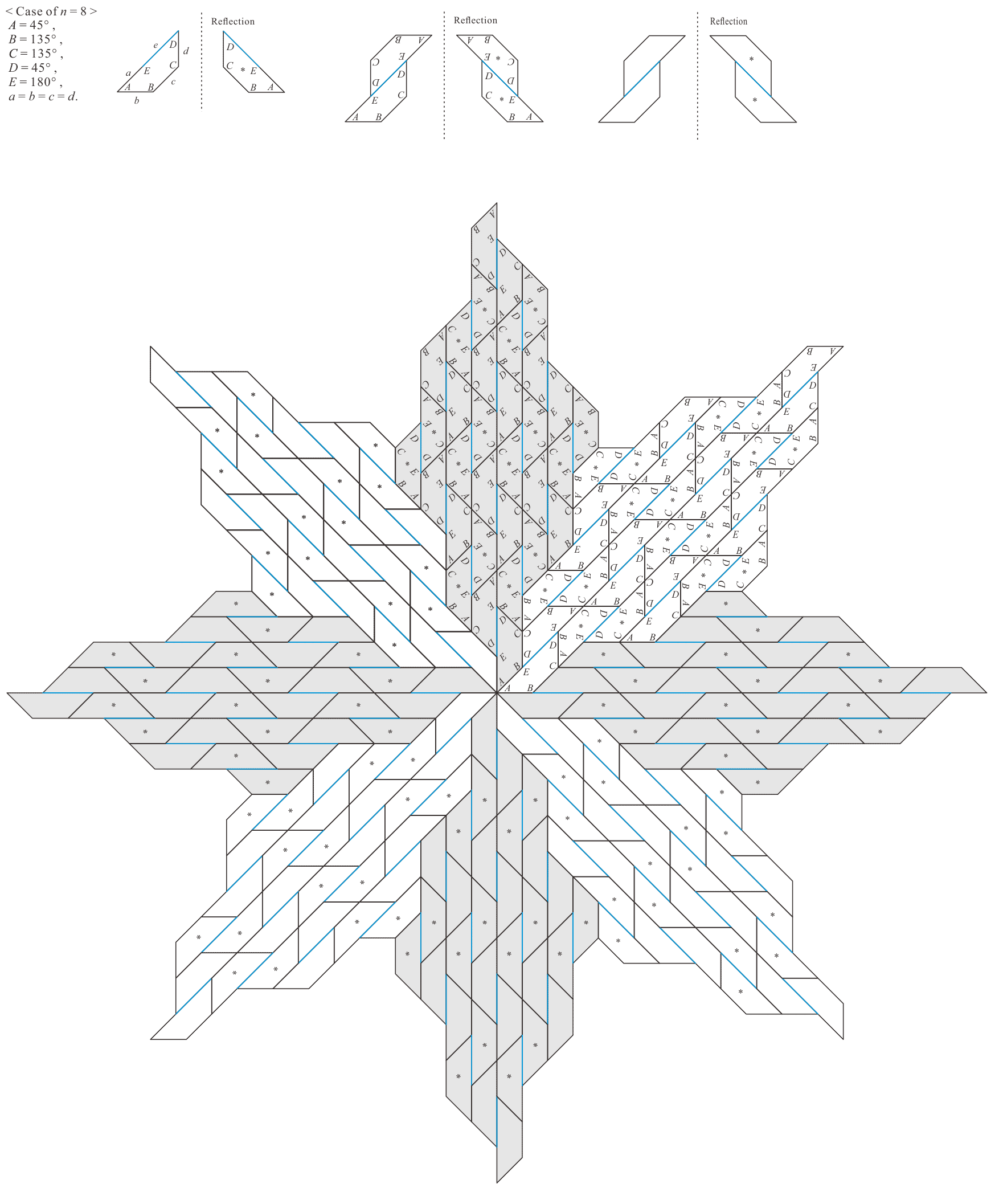} 
  \caption{{\small 
Eight-fold rotationally symmetric tiling by a trapezoid of $n=8$ 
in Table~\ref{tab2} 
 (The figure is solely a depiction of  the area around the rotationally symmetric
 center, and the tiling can be spread in all directions as well)
} 
\label{fig20}
}
\end{figure}

\renewcommand{\figurename}{{\small Figure.}}
\begin{figure}[htbp]
 \centering\includegraphics[width=15cm,clip]{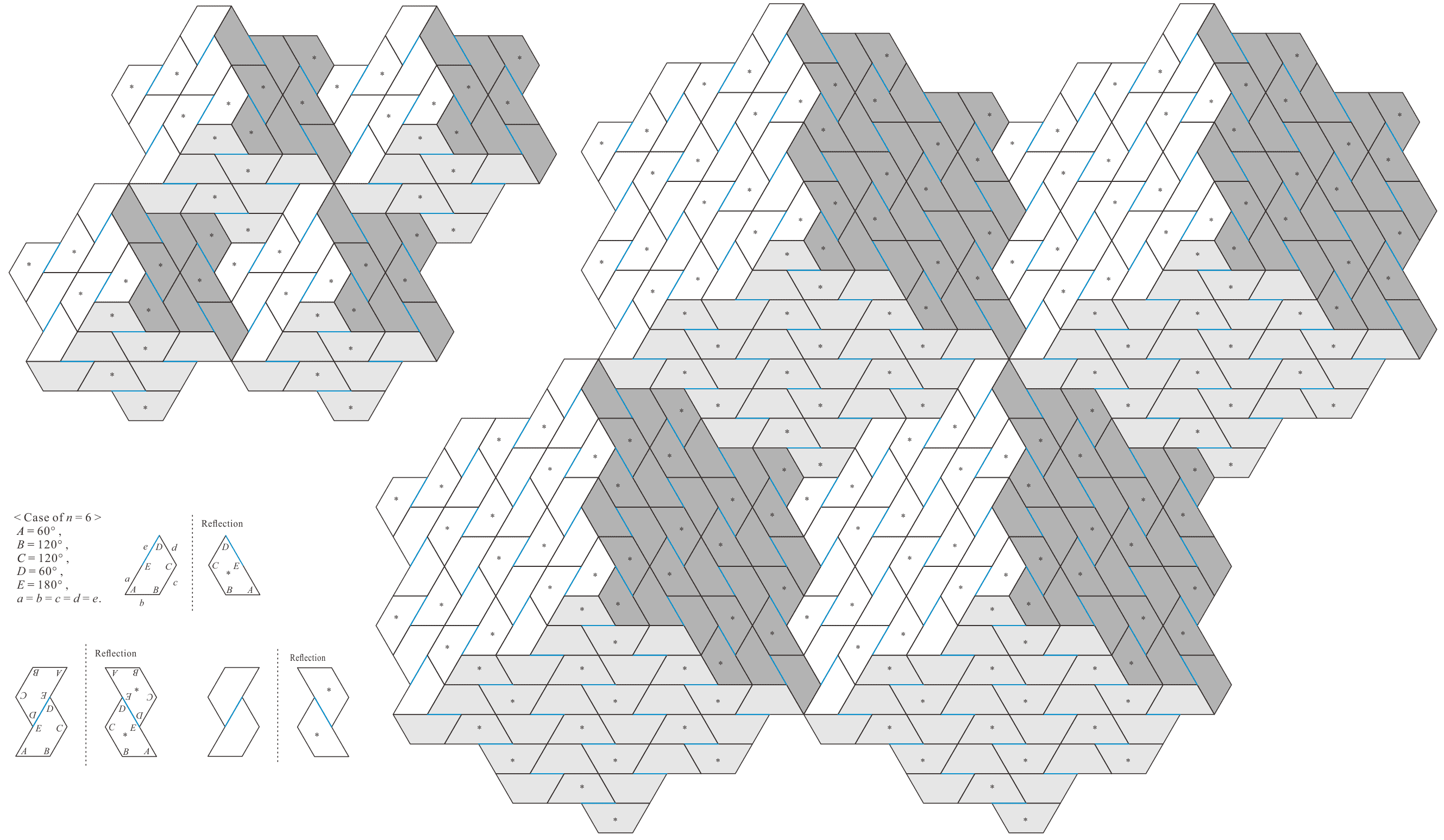} 
  \caption{{\small 
Examples of tilings with three-fold and six-fold rotational symmetry by 
a trapezoid that corresponds to rhombus with an acute angle 
of $60^ \circ$ } 
\label{fig21}
}
\end{figure}

\renewcommand{\figurename}{{\small Figure.}}
\begin{figure}[htbp]
 \centering\includegraphics[width=15cm,clip]{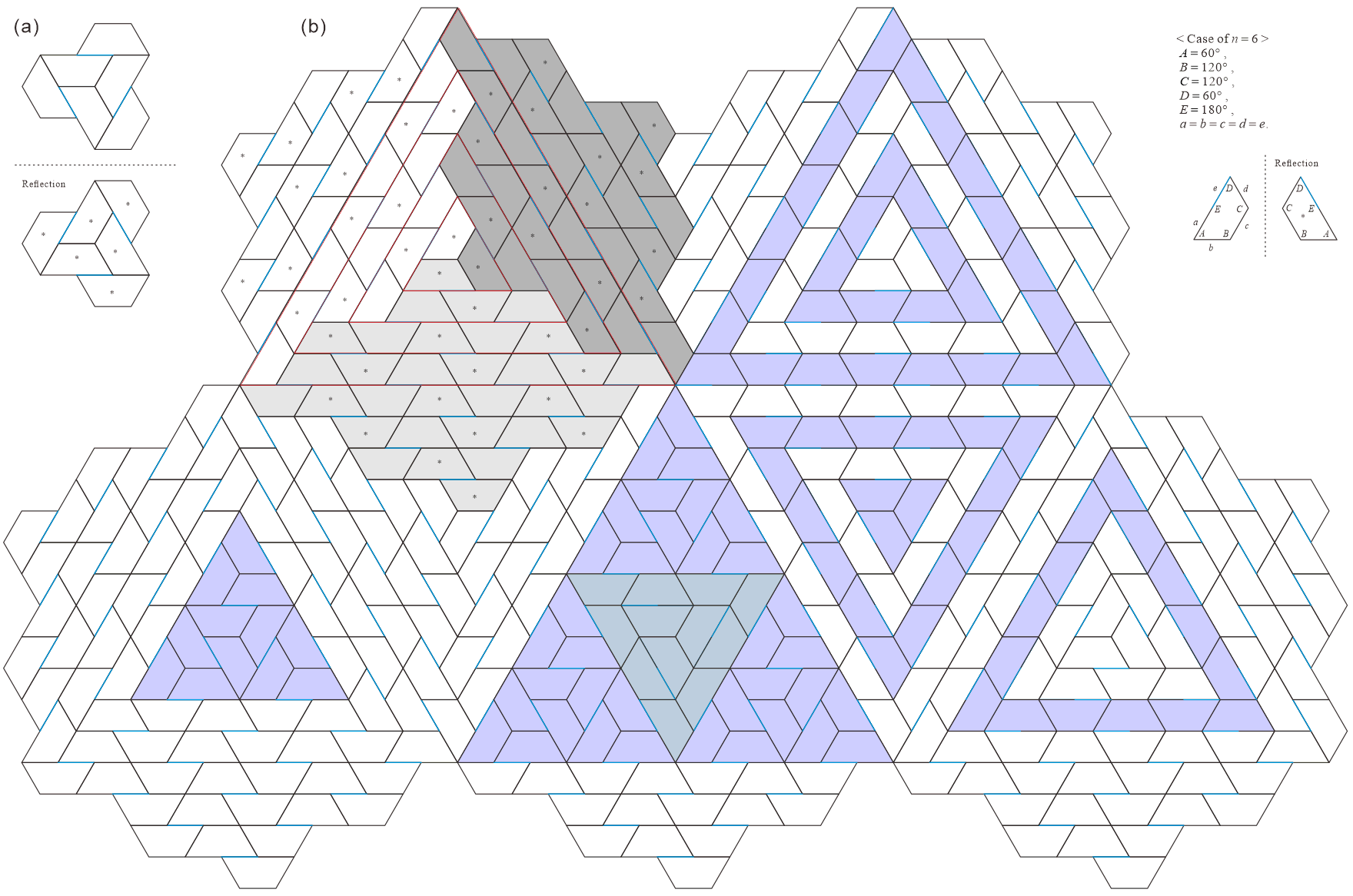} 
  \caption{{\small 
Examples of tilings by a trapezoid that corresponds to rhombus 
with an acute angle of $60^ \circ$ } 
\label{fig22}
}
\end{figure}

\renewcommand{\figurename}{{\small Figure.}}
\begin{figure}[htbp]
 \centering\includegraphics[width=15cm,clip]{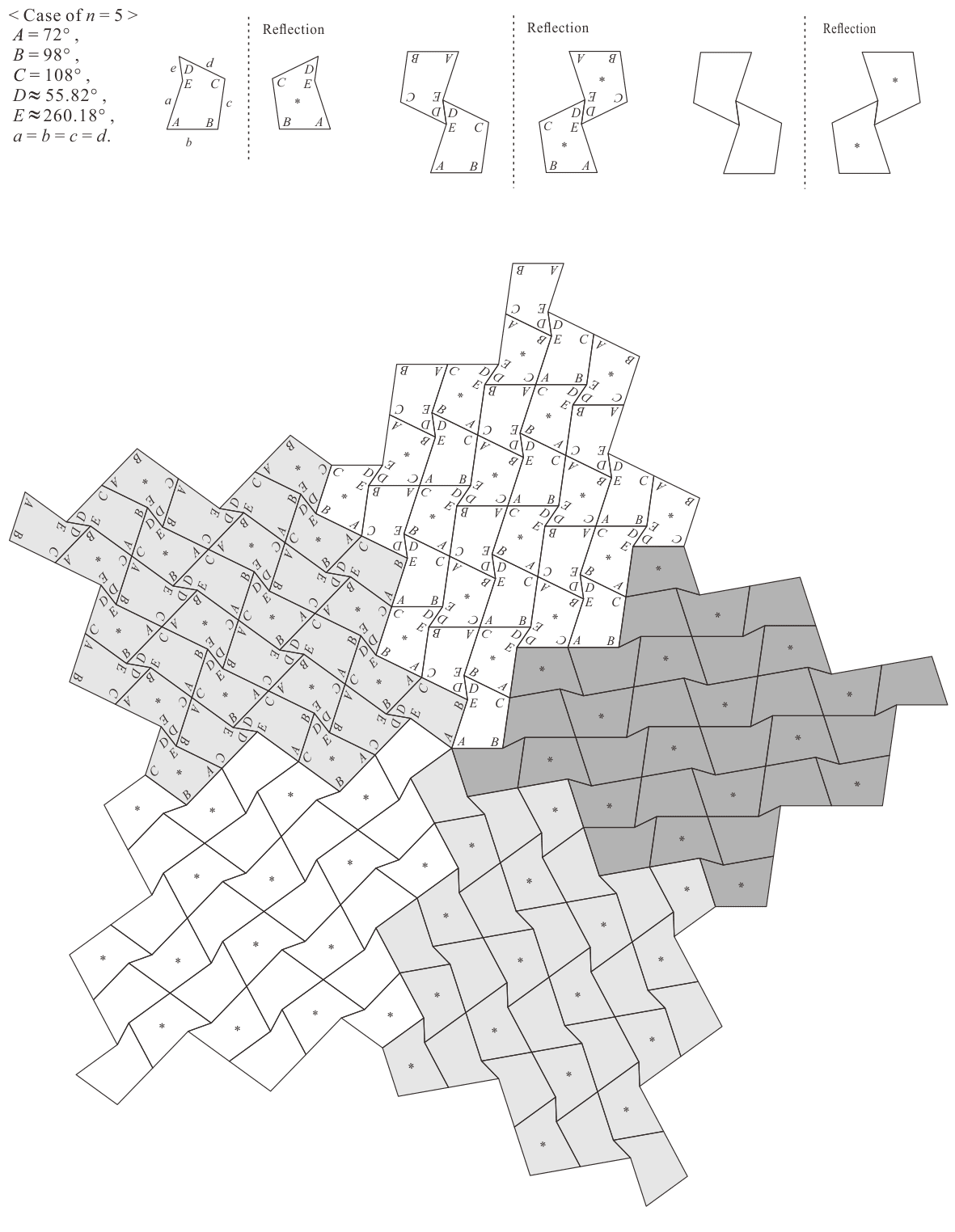} 
  \caption{{\small 
Five-fold rotationally symmetric tiling by a concave 
pentagon of $n=5$ in Table~\ref{tab3} 
 (The figure is solely a depiction of  the area around the rotationally symmetric
 center, and the tiling can be spread in all directions as well)
}
\label{fig23}
}
\end{figure}

\renewcommand{\figurename}{{\small Figure.}}
\begin{figure}[htbp]
 \centering\includegraphics[width=15cm,clip]{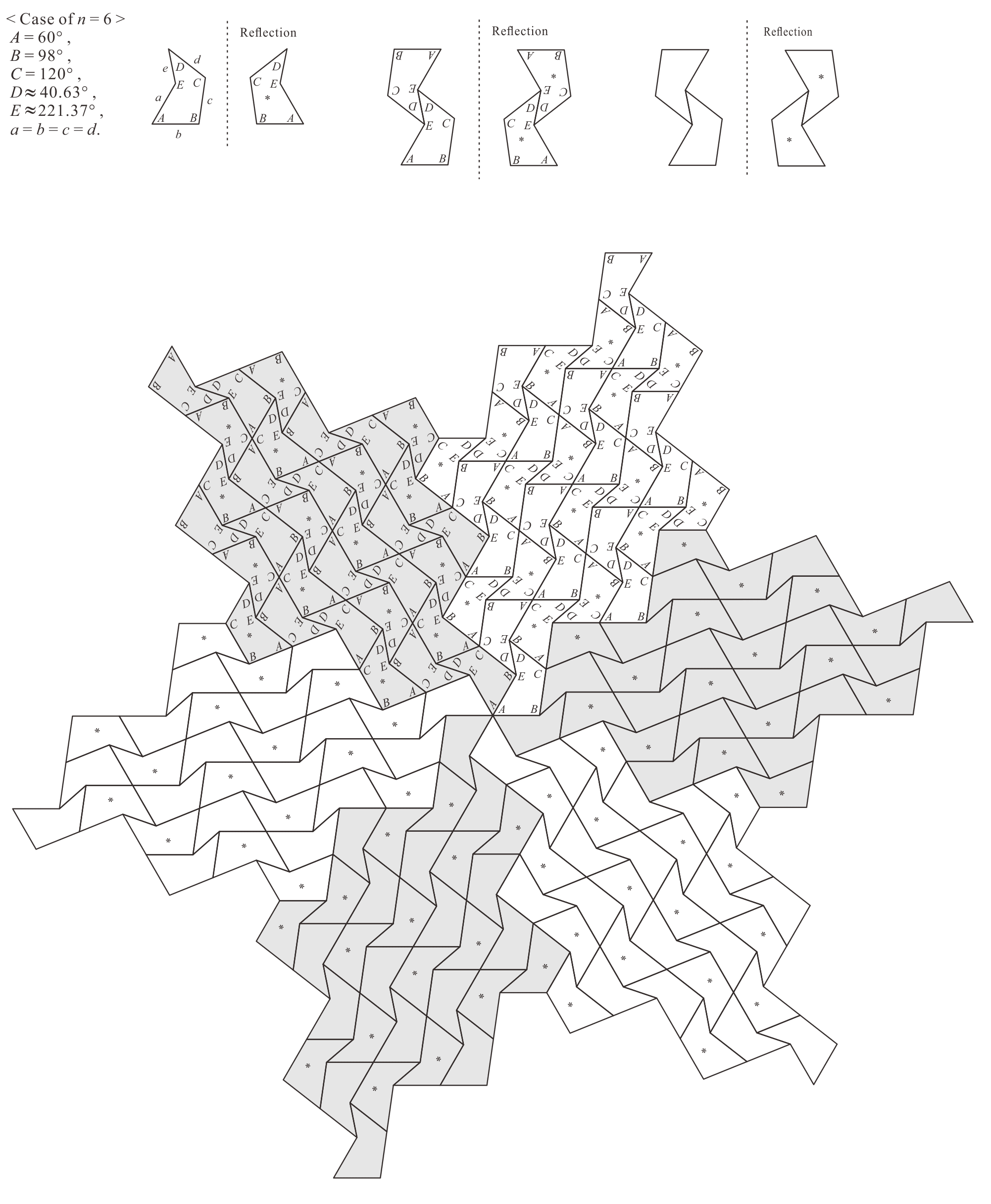} 
  \caption{{\small 
Six-fold rotationally symmetric tiling by a concave 
pentagon of $n=6$ in Table~\ref{tab3} 
 (The figure is solely a depiction of  the area around the rotationally symmetric
 center, and the tiling can be spread in all directions as well)
}
\label{fig24}
}
\end{figure}

\renewcommand{\figurename}{{\small Figure.}}
\begin{figure}[htbp]
 \centering\includegraphics[width=15cm,clip]{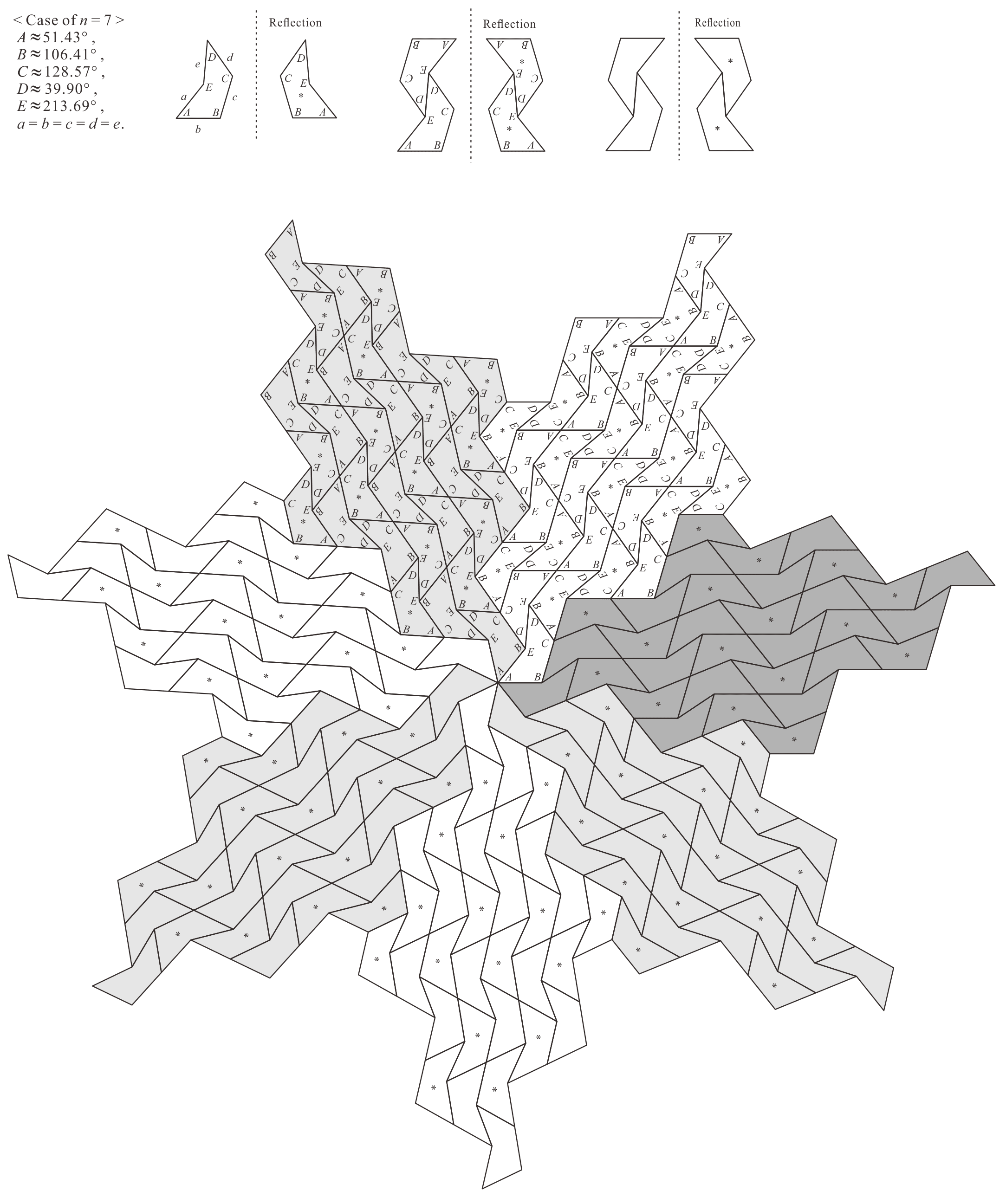} 
  \caption{{\small 
Seven-fold rotationally symmetric tiling by a concave 
pentagon of $n=7$ in Table~\ref{tab3} 
 (The figure is solely a depiction of  the area around the rotationally symmetric
 center, and the tiling can be spread in all directions as well)
}
\label{fig25}
}
\end{figure}

\renewcommand{\figurename}{{\small Figure.}}
\begin{figure}[htbp]
 \centering\includegraphics[width=15cm,clip]{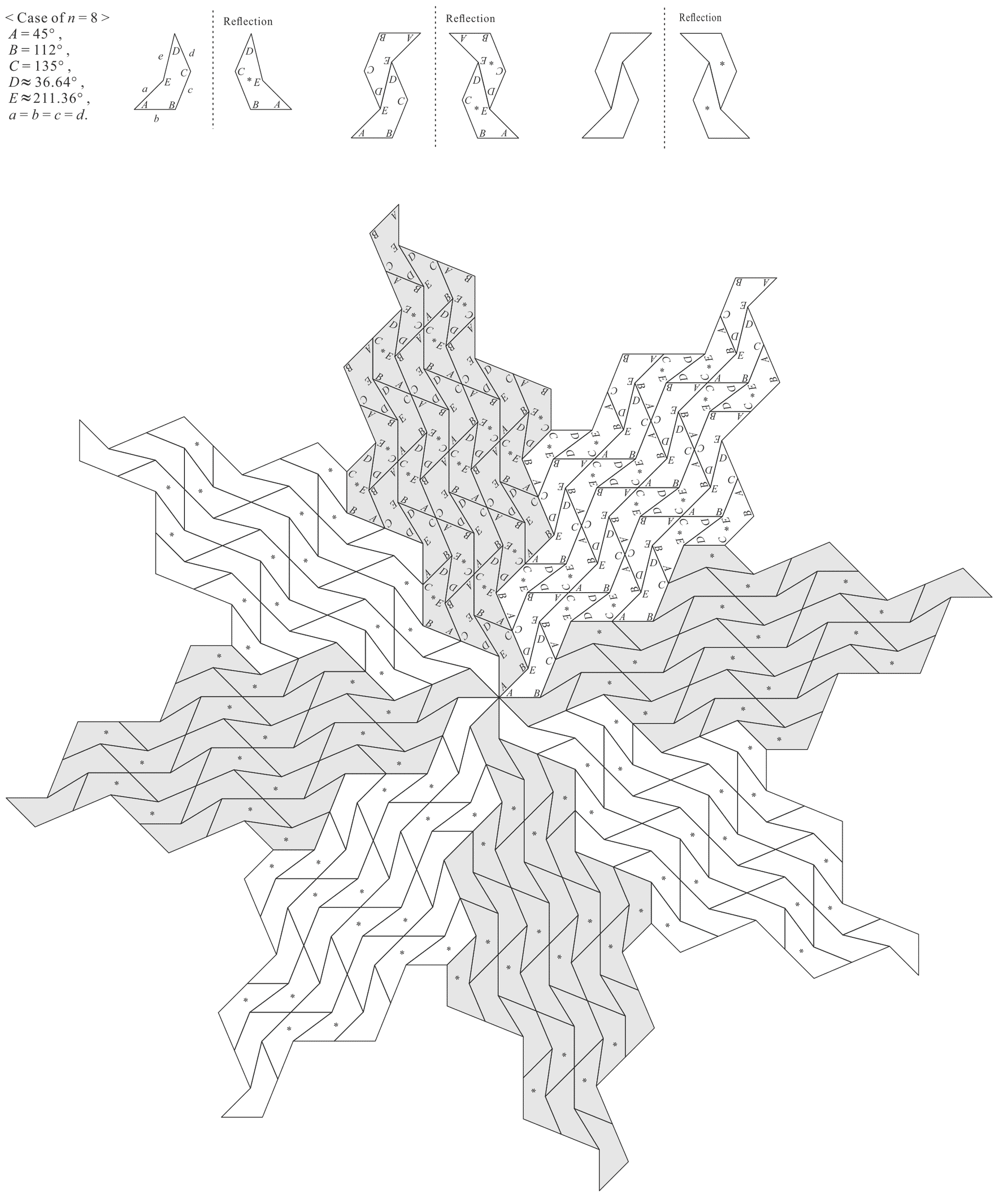} 
  \caption{{\small 
Eight-fold rotationally symmetric tiling by a concave 
pentagon of $n=8$ in Table~\ref{tab3} 
 (The figure is solely a depiction of  the area around the rotationally symmetric
 center, and the tiling can be spread in all directions as well)
}
\label{fig26}
}
\end{figure}

\renewcommand{\figurename}{{\small Figure.}}
\begin{figure}[htbp]
 \centering\includegraphics[width=15cm,clip]{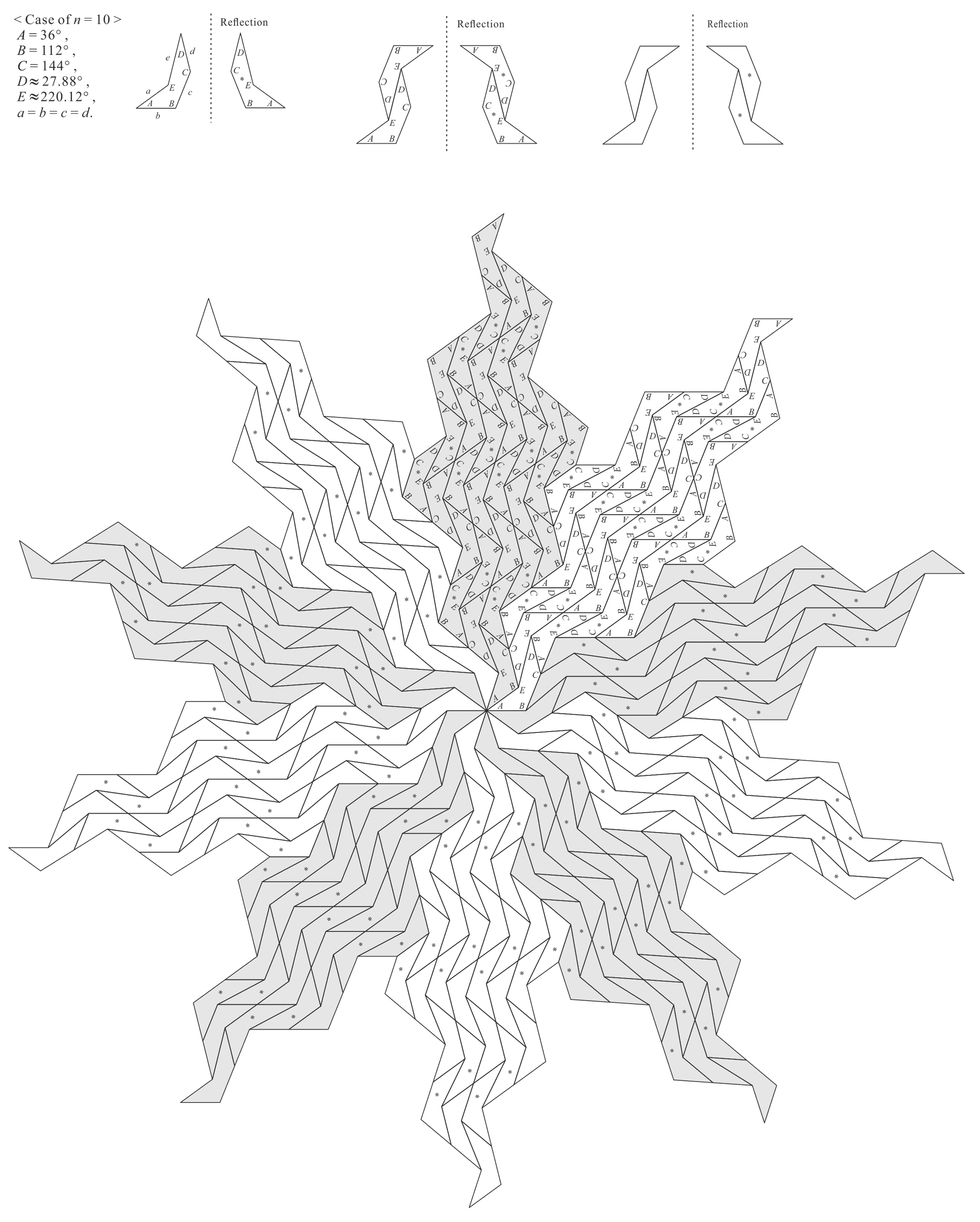} 
  \caption{{\small 
10-fold rotationally symmetric tiling by a concave 
pentagon of $n=10$ in Table~\ref{tab3} 
 (The figure is solely a depiction of  the area around the rotationally symmetric
 center, and the tiling can be spread in all directions as well)
}
\label{fig27}
}
\end{figure}

\renewcommand{\figurename}{{\small Figure.}}
\begin{figure}[htbp]
 \centering\includegraphics[width=15cm,clip]{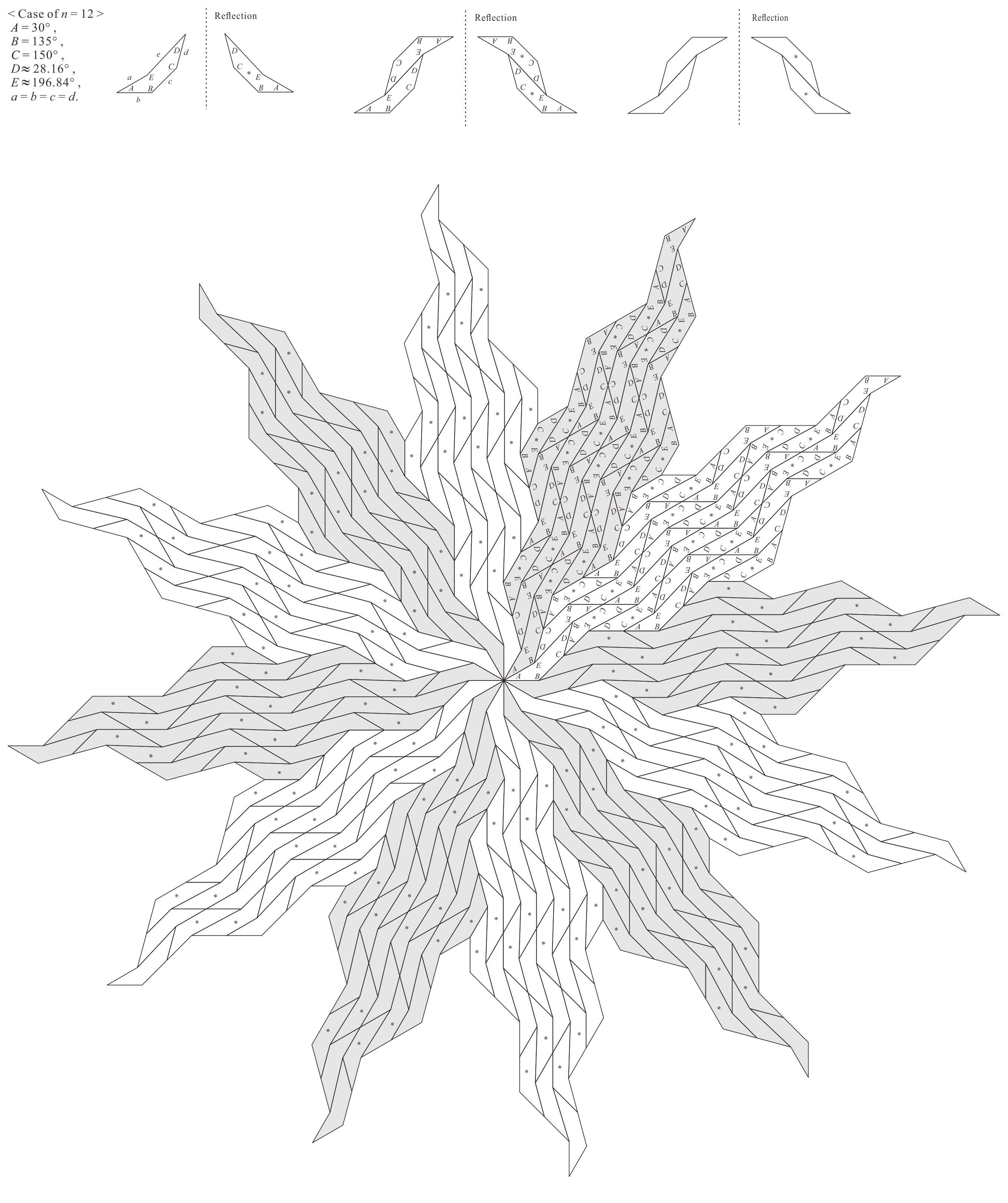} 
  \caption{{\small 
12-fold rotationally symmetric tiling by a concave 
pentagon of $n=12$ in Table~\ref{tab3} 
 (The figure is solely a depiction of  the area around the rotationally symmetric
 center, and the tiling can be spread in all directions as well)
}
\label{fig28}
}
\end{figure}

\renewcommand{\figurename}{{\small Figure.}}
\begin{figure}[htbp]
 \centering\includegraphics[width=15cm,clip]{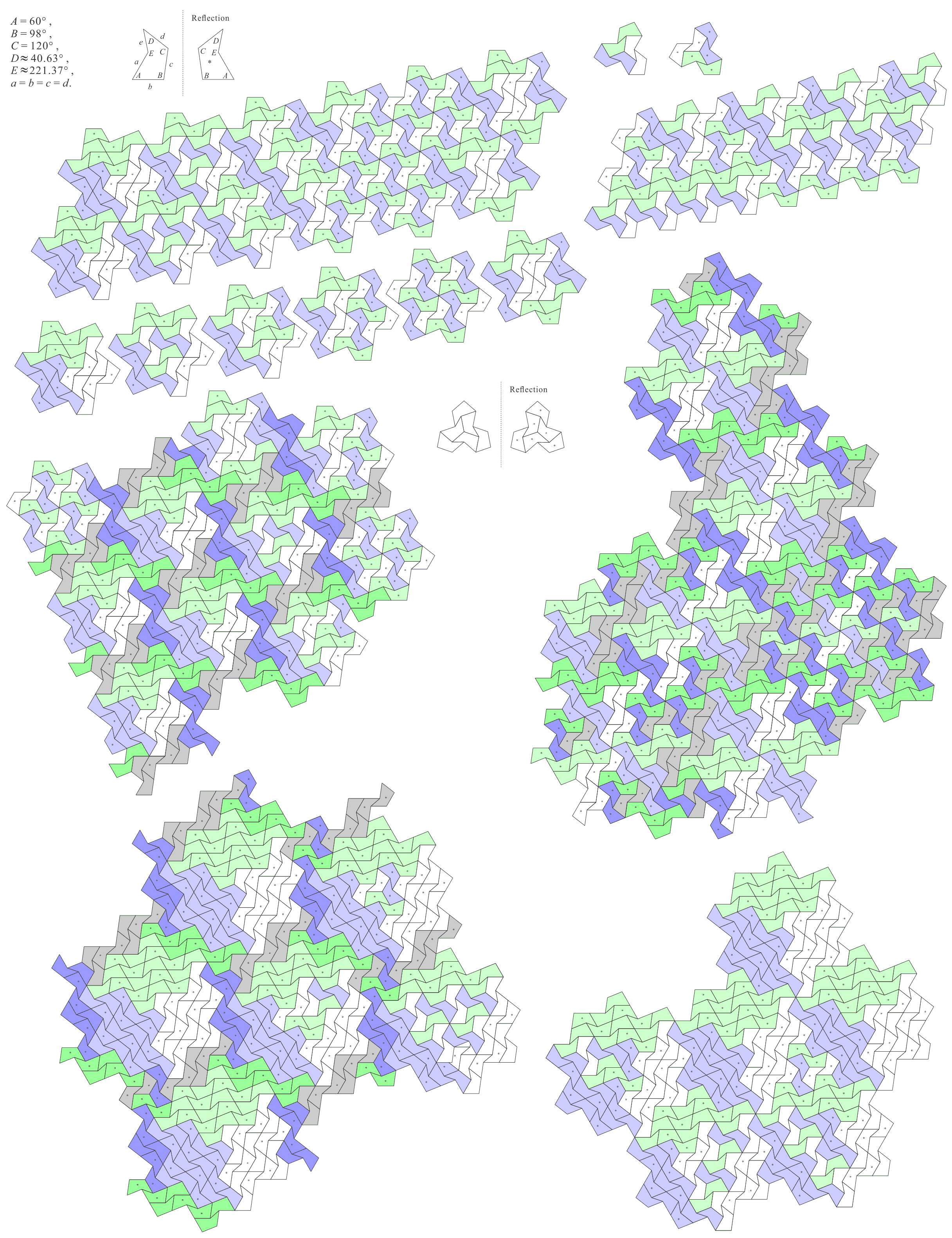} 
  \caption{{\small 
Examples of tilings by a concave pentagon that corresponds to rhombus 
with an acute angle of $60^ \circ$ }
\label{fig29}
}
\end{figure}

\renewcommand{\figurename}{{\small Figure.}}
\begin{figure}[htbp]
 \centering\includegraphics[width=15cm,clip]{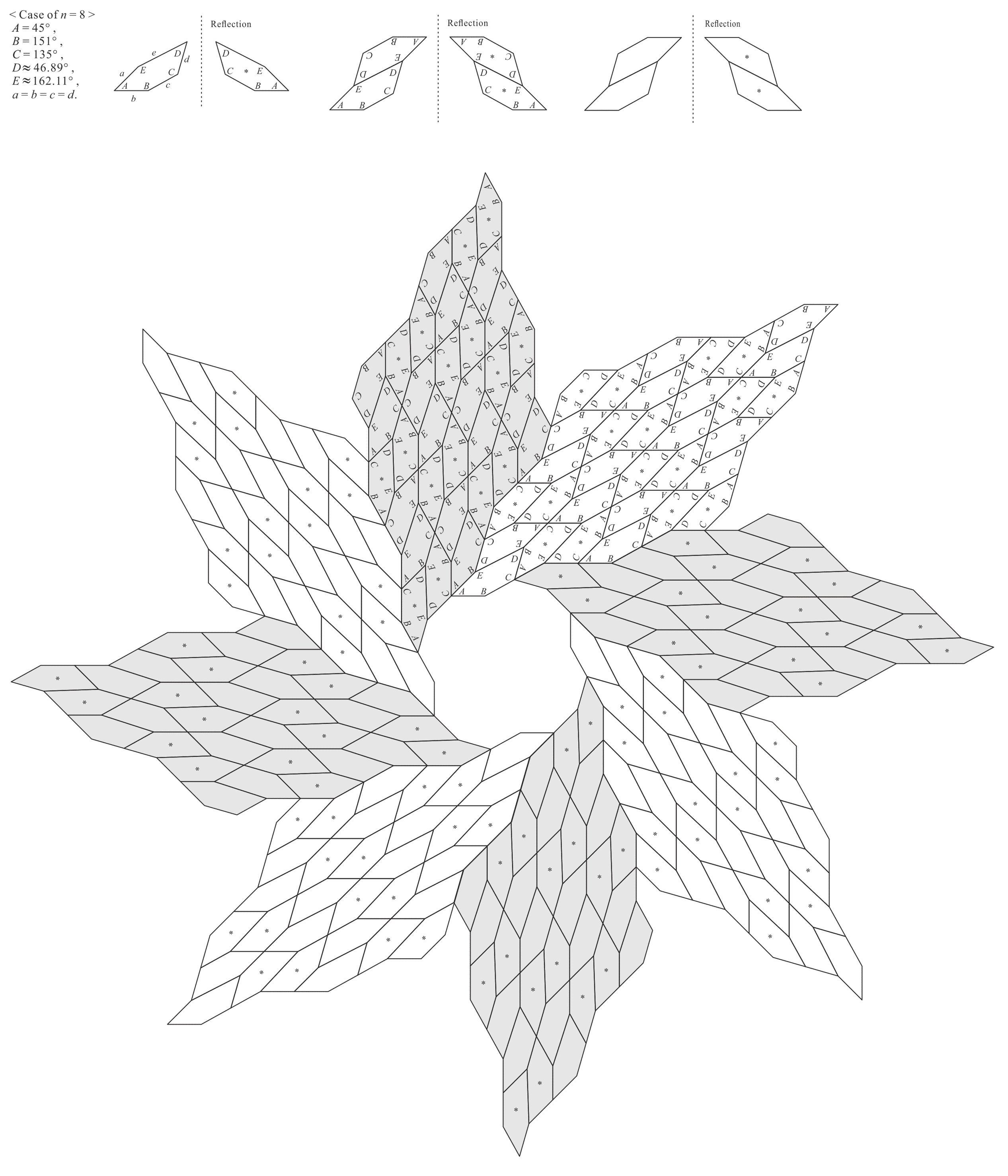} 
  \caption{{\small 
Rotationally symmetric tiling with $C_{4}$ symmetry, with an equilateral 
concave 16-gonal hole with $D_{4}$ symmetry at the center, by a convex 
pentagon of $n = 8$ in Table~\ref{tab1} 
 (The figure is solely a depiction of  the area around the rotationally symmetric
 center, and the tiling can be spread in all directions as well)
} 
\label{fig30}
}
\end{figure}

\renewcommand{\figurename}{{\small Figure.}}
\begin{figure}[htbp]
 \centering\includegraphics[width=15cm,clip]{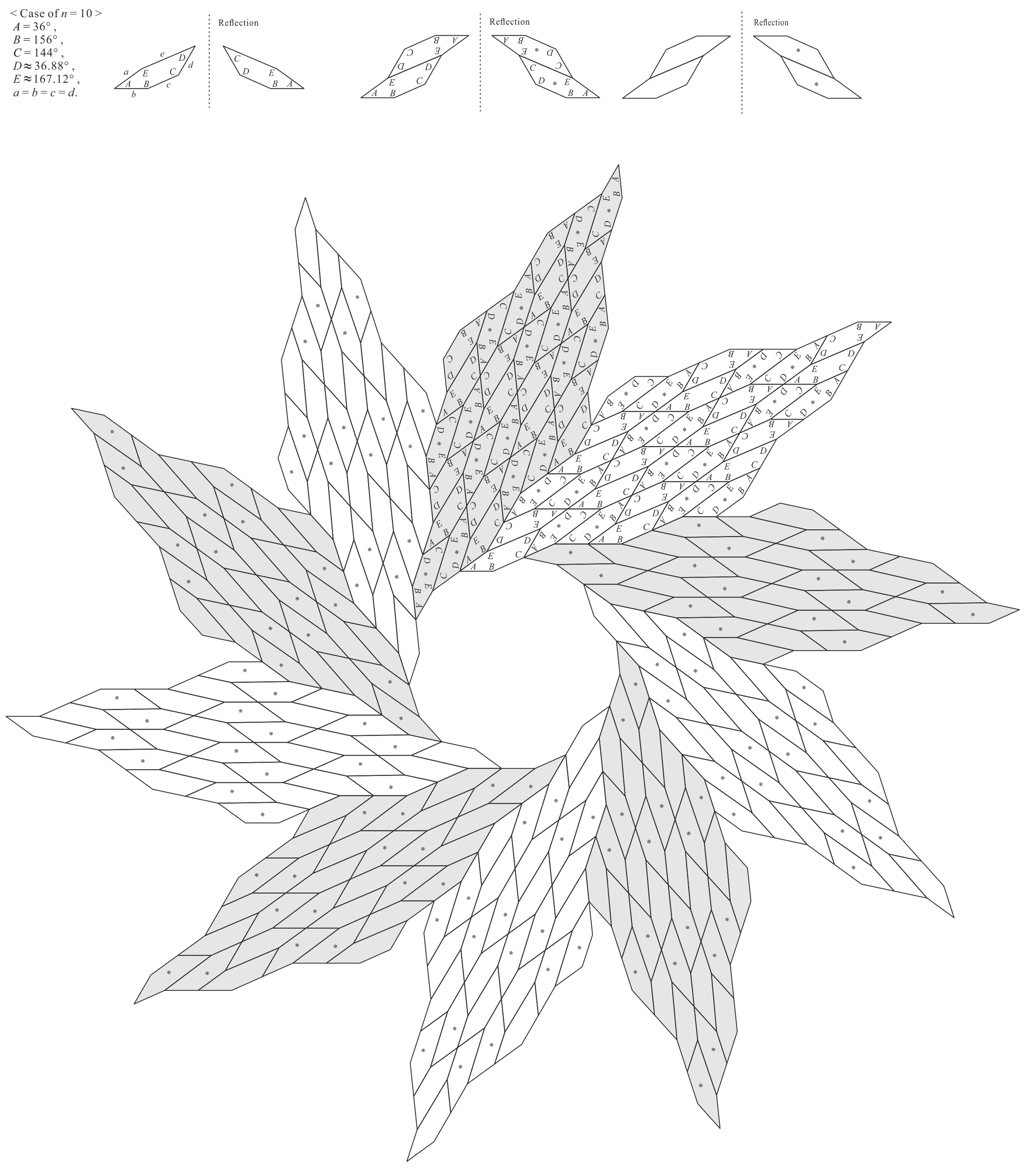} 
  \caption{{\small 
Rotationally symmetric tiling with $C_{5}$ symmetry, with an equilateral 
concave 20-gonal hole with $D_{5}$ symmetry at the center, by a convex 
pentagon of $n = 10$ in Table~\ref{tab1} 
 (The figure is solely a depiction of  the area around the rotationally symmetric
 center, and the tiling can be spread in all directions as well)
} 
\label{fig31}
}
\end{figure}

\renewcommand{\figurename}{{\small Figure.}}
\begin{figure}[htbp]
 \centering\includegraphics[width=15cm,clip]{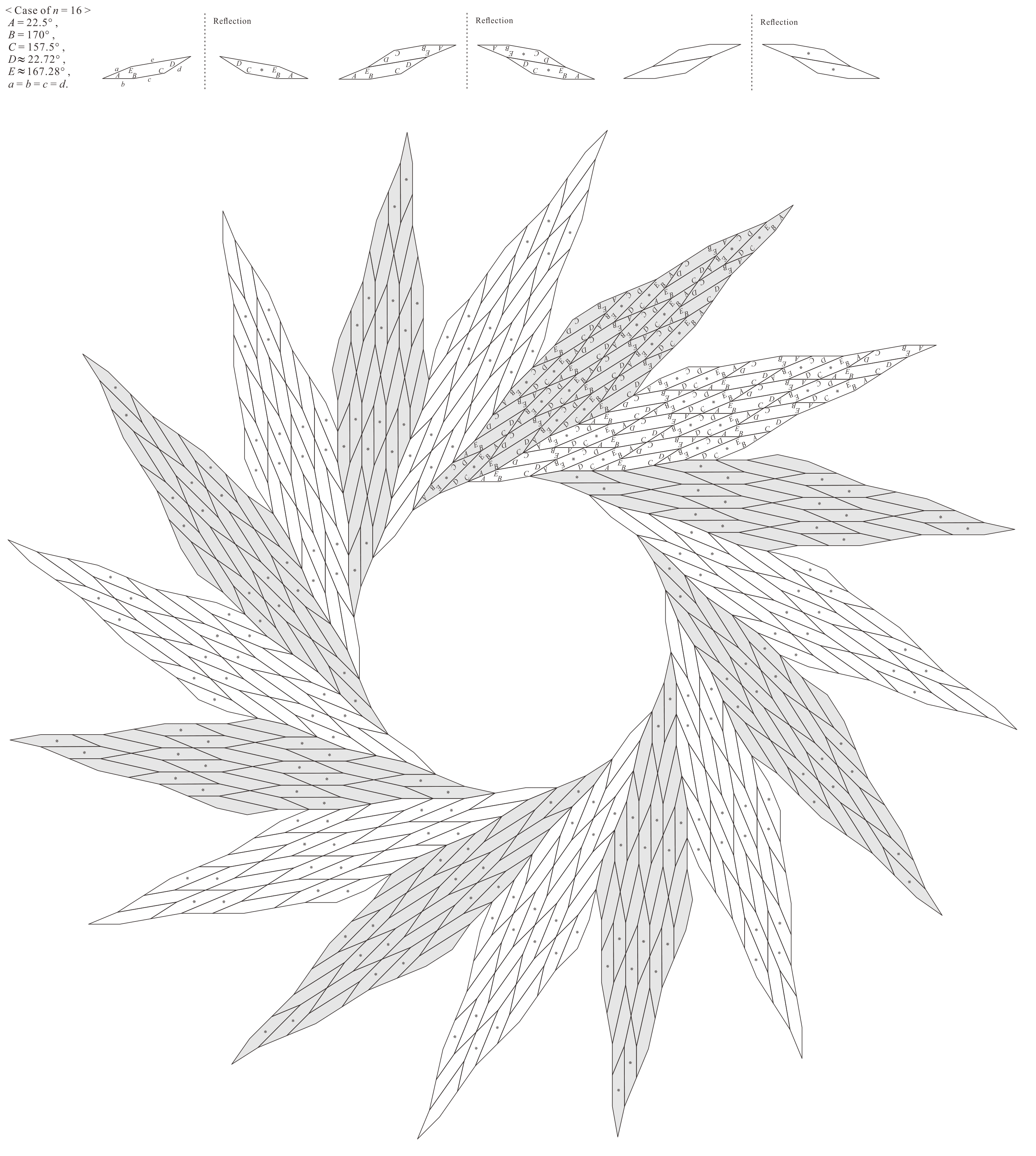} 
  \caption{{\small 
Rotationally symmetric tiling with $C_{8}$ symmetry, with an equilateral 
concave 32-gonal hole with $D_{8}$ symmetry at the center, by a convex 
pentagon of $n = 16$ in Table~\ref{tab1} 
 (The figure is solely a depiction of  the area around the rotationally symmetric
 center, and the tiling can be spread in all directions as well)
} 
\label{fig32}
}
\end{figure}

\renewcommand{\figurename}{{\small Figure.}}
\begin{figure}[htbp]
 \centering\includegraphics[width=15cm,clip]{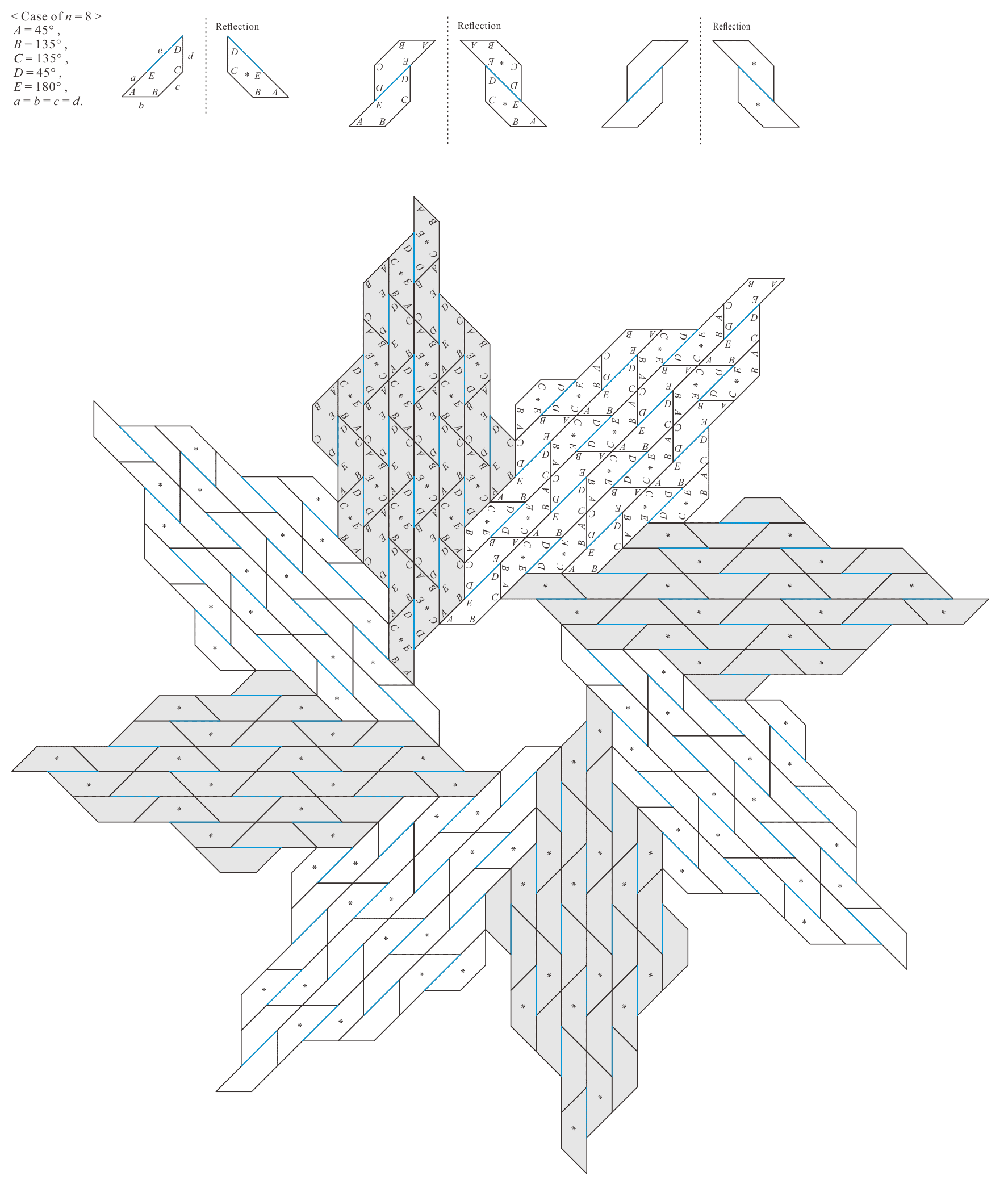} 
  \caption{{\small 
Rotationally symmetric tiling with $C_{4}$ symmetry, with an equilateral 
concave 16-gonal hole with $D_{4}$ symmetry at the center, by a
trapezoid of $n = 8$ in Table~\ref{tab2} 
 (The figure is solely a depiction of  the area around the rotationally symmetric
 center, and the tiling can be spread in all directions as well)
} 
\label{fig33}
}
\end{figure}

\renewcommand{\figurename}{{\small Figure.}}
\begin{figure}[htbp]
 \centering\includegraphics[width=15cm,clip]{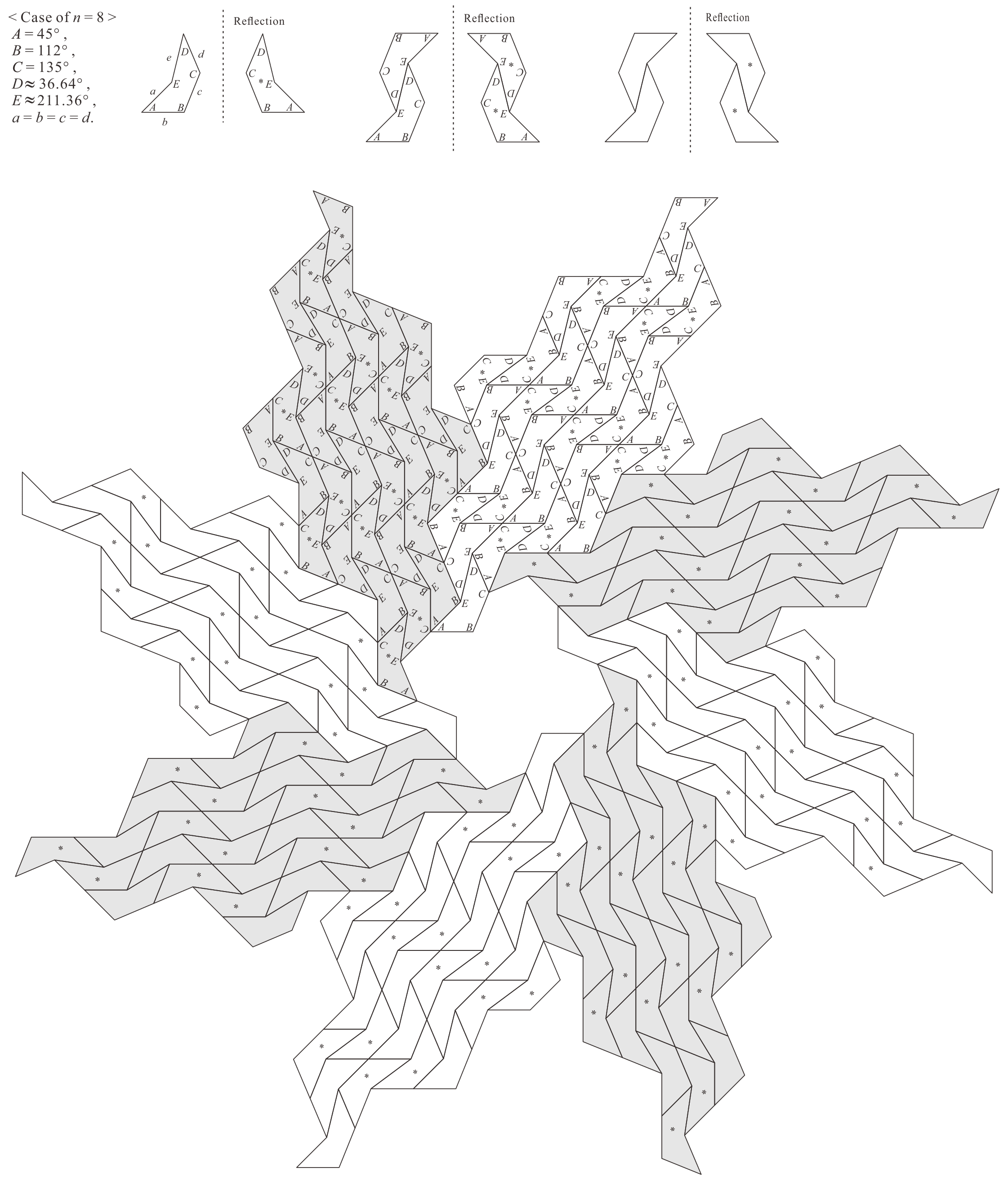} 
  \caption{{\small 
Rotationally symmetric tiling with $C_{4}$ symmetry, with an equilateral 
concave 16-gonal hole with $D_{4}$ symmetry at the center, by a concave 
pentagon of $n = 8$ in Table~\ref{tab3}  
 (The figure is solely a depiction of  the area around the rotationally symmetric
 center, and the tiling can be spread in all directions as well)
} 
\label{fig34}
}
\end{figure}

\renewcommand{\figurename}{{\small Figure.}}
\begin{figure}[htbp]
 \centering\includegraphics[width=15cm,clip]{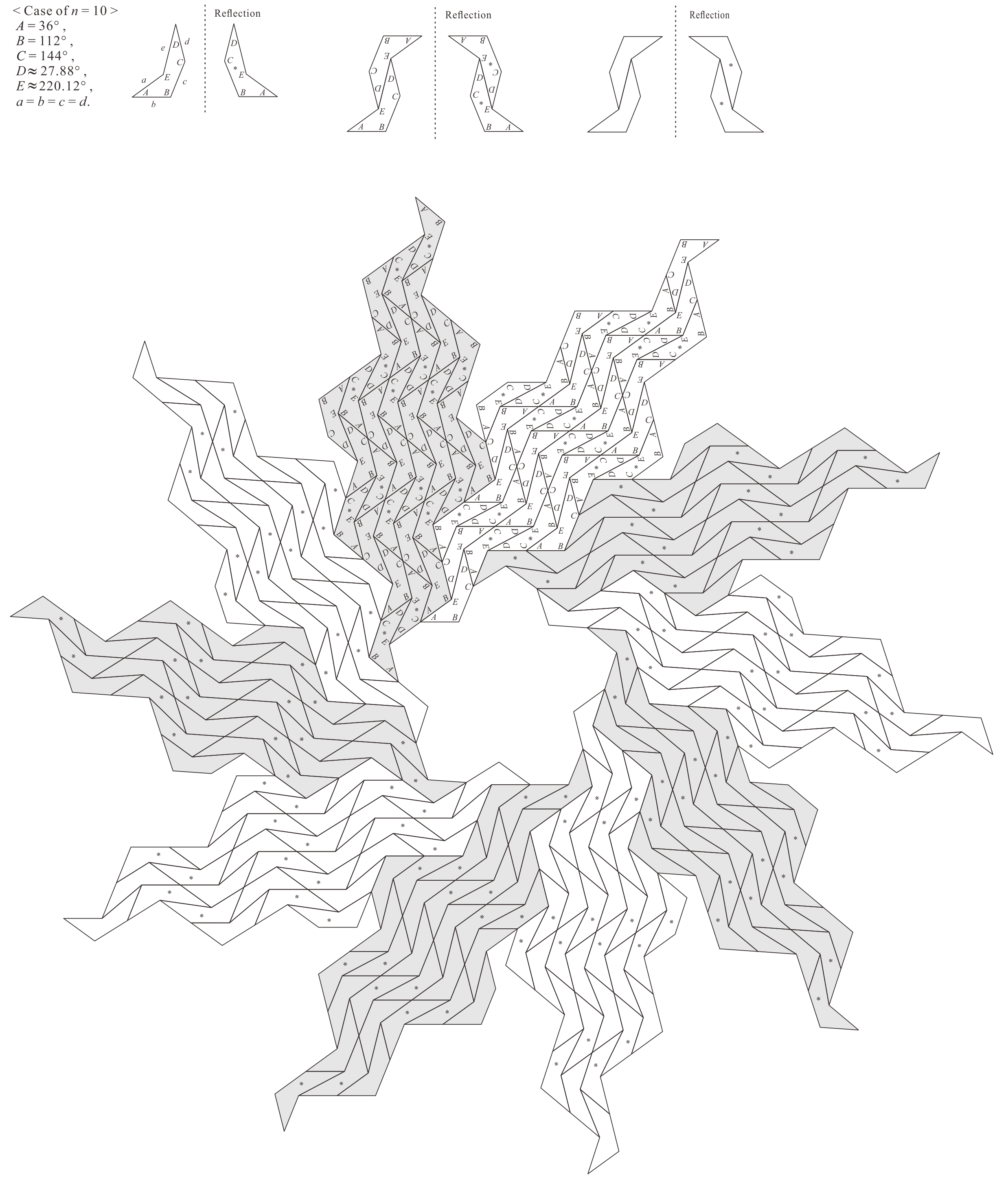} 
  \caption{{\small 
Rotationally symmetric tiling with $C_{5}$ symmetry, with an equilateral 
concave 20-gonal hole with $D_{5}$ symmetry at the center, by a concave 
pentagon of $n = 10$ in Table~\ref{tab3} 
 (The figure is solely a depiction of  the area around the rotationally symmetric
 center, and the tiling can be spread in all directions as well)
} 
\label{fig35}
}
\end{figure}

\renewcommand{\figurename}{{\small Figure.}}
\begin{figure}[htbp]
 \centering\includegraphics[width=15cm,clip]{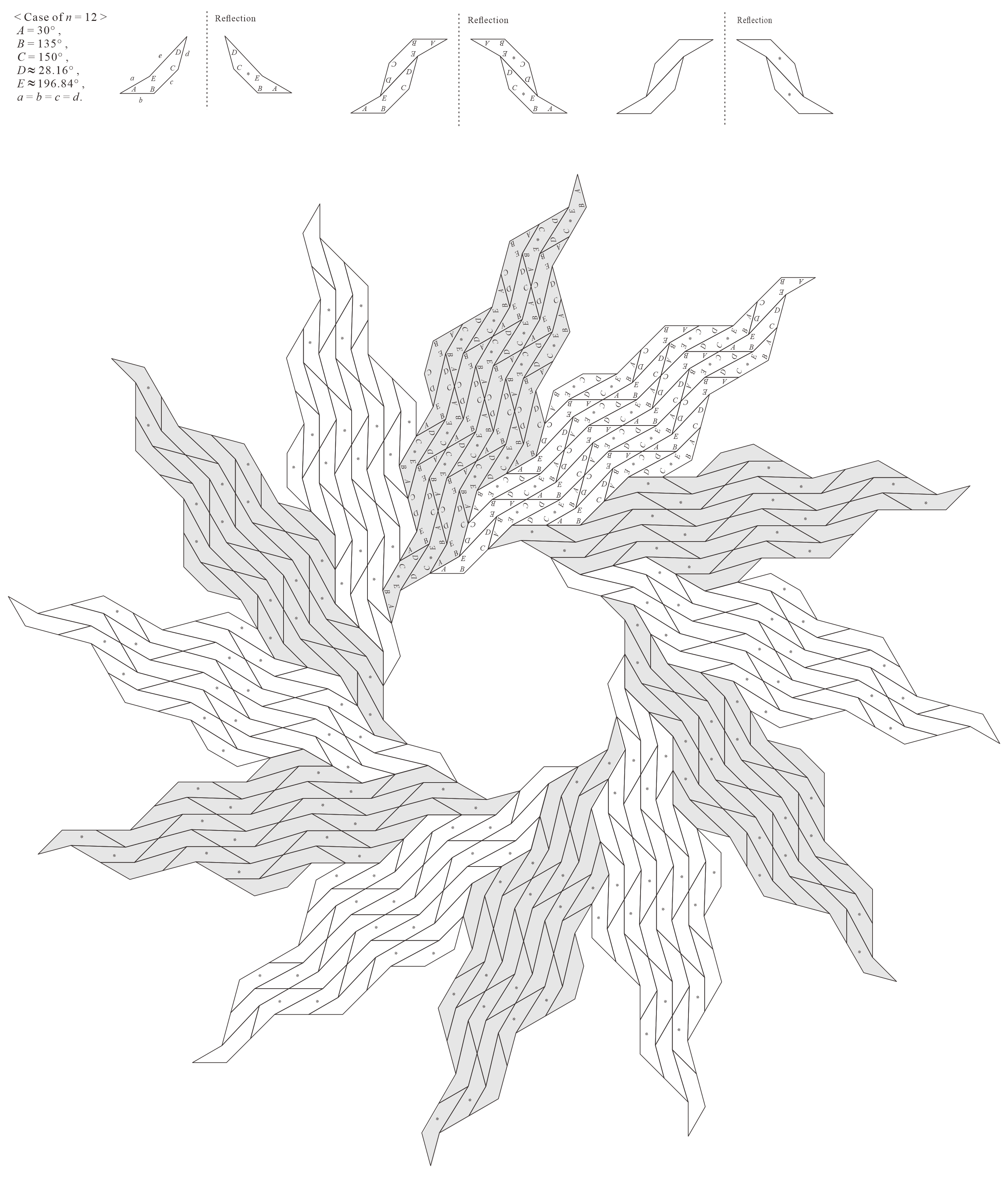} 
  \caption{{\small 
Rotationally symmetric tiling with $C_{6}$ symmetry, with an equilateral 
concave 24-gonal hole with $D_{6}$ symmetry at the center, by a concave 
pentagon of $n = 12$ in Table~\ref{tab3} 
 (The figure is solely a depiction of  the area around the rotationally symmetric
 center, and the tiling can be spread in all directions as well)
} 
\label{fig36}
}
\end{figure}

\renewcommand{\figurename}{{\small Figure.}}
\begin{figure}[htbp]
 \centering\includegraphics[width=15cm,clip]{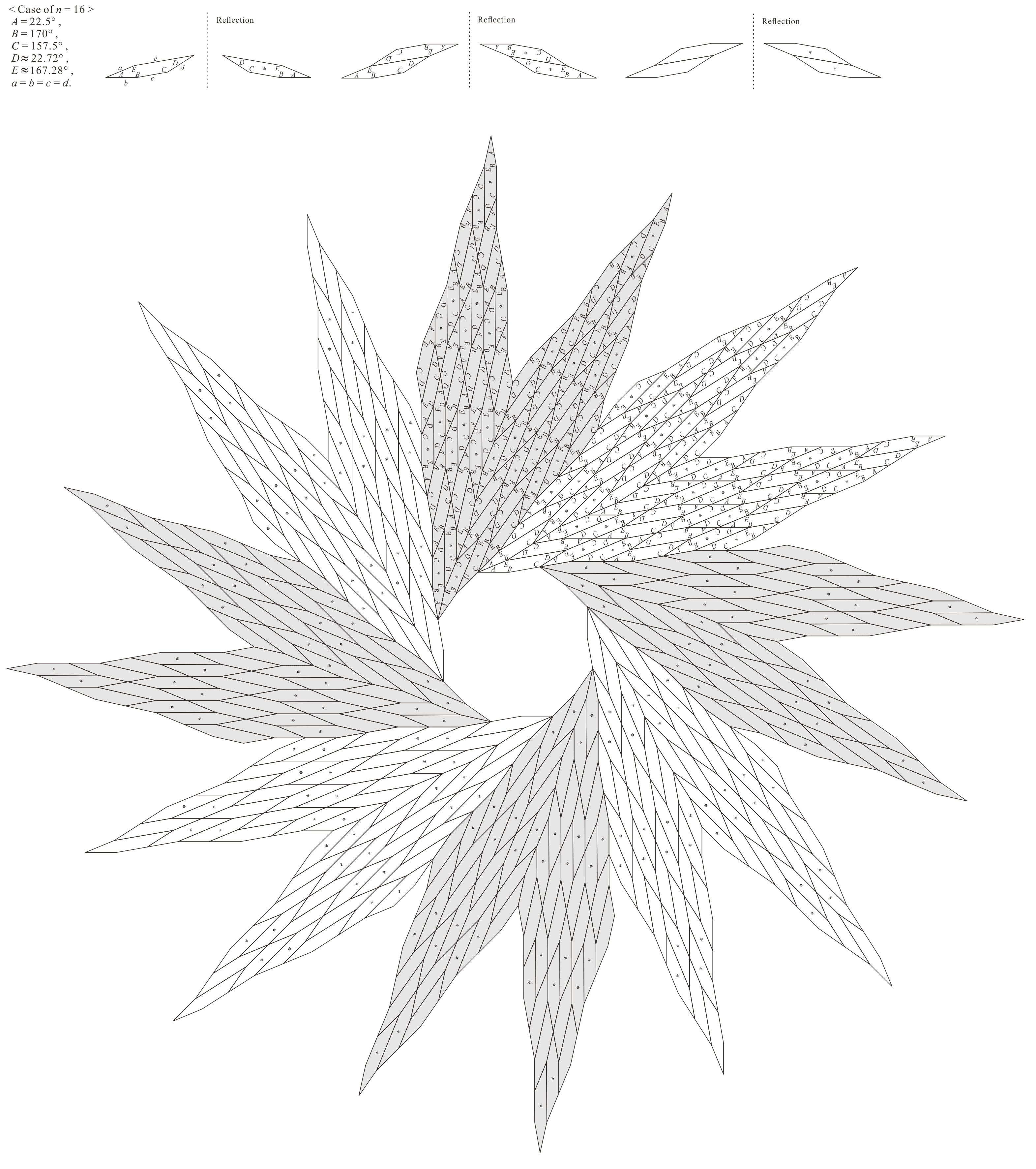} 
  \caption{{\small 
Rotationally symmetric tiling with $C_{4}$ symmetry, with an equilateral 
concave 16-gonal hole with $D_{4}$ symmetry at the center, by a convex 
pentagon of $n = 16$ in Table~\ref{tab1} 
 (The figure is solely a depiction of  the area around the rotationally symmetric
 center, and the tiling can be spread in all directions as well)
} 
\label{fig37}
}
\end{figure}

\renewcommand{\figurename}{{\small Figure.}}
\begin{figure}[htbp]
 \centering\includegraphics[width=15cm,clip]{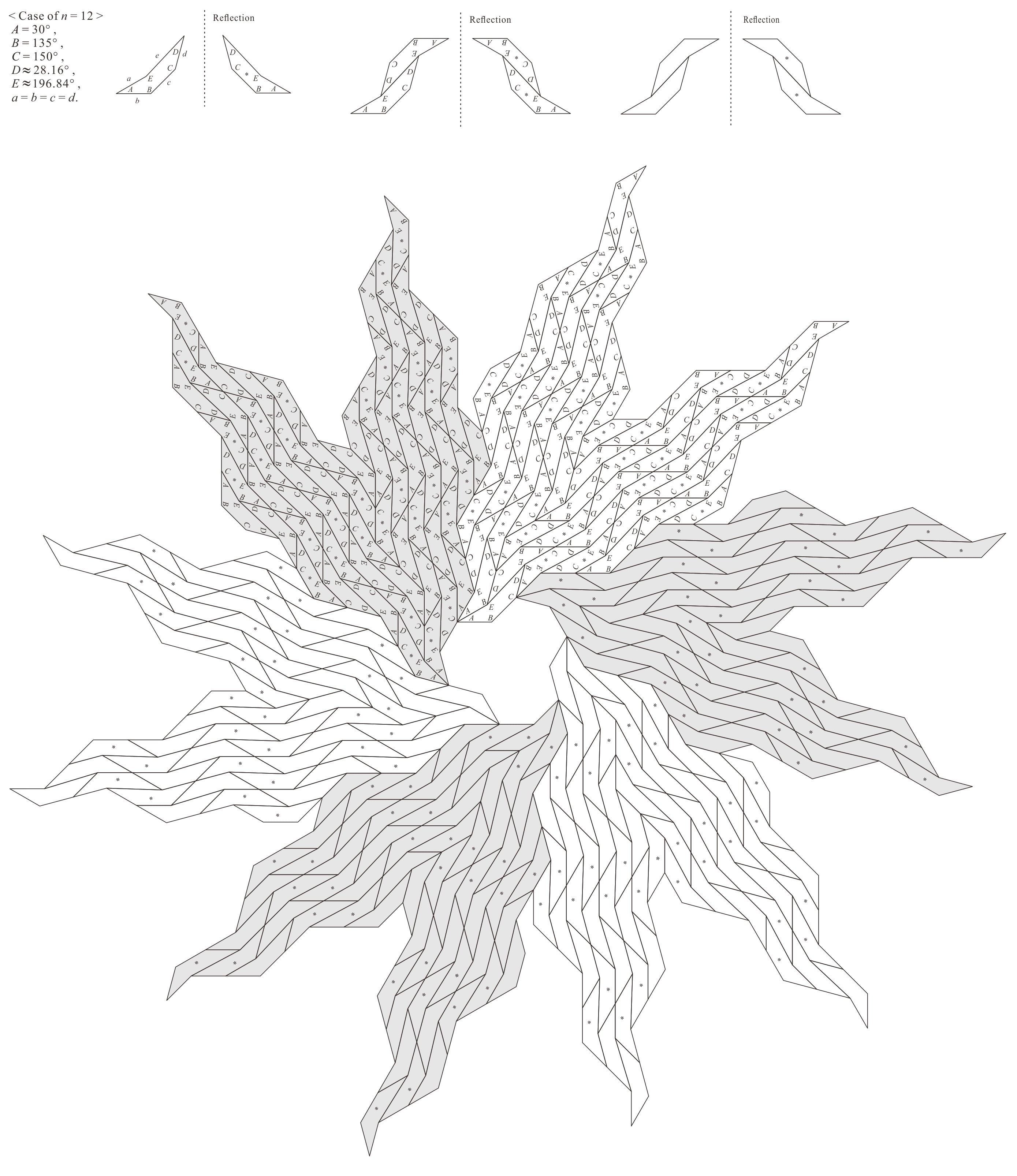} 
  \caption{{\small 
Rotationally symmetric tiling with $C_{3}$ symmetry, with an equilateral 
concave 12-gonal hole with $D_{3}$ symmetry at the center, by a concave 
pentagon of $n = 12$ in Table~\ref{tab3} 
 (The figure is solely a depiction of  the area around the rotationally symmetric
 center, and the tiling can be spread in all directions as well)
} 
\label{fig38}
}
\end{figure}

\renewcommand{\figurename}{{\small Figure.}}
\begin{figure}[htbp]
 \centering\includegraphics[width=15cm,clip]{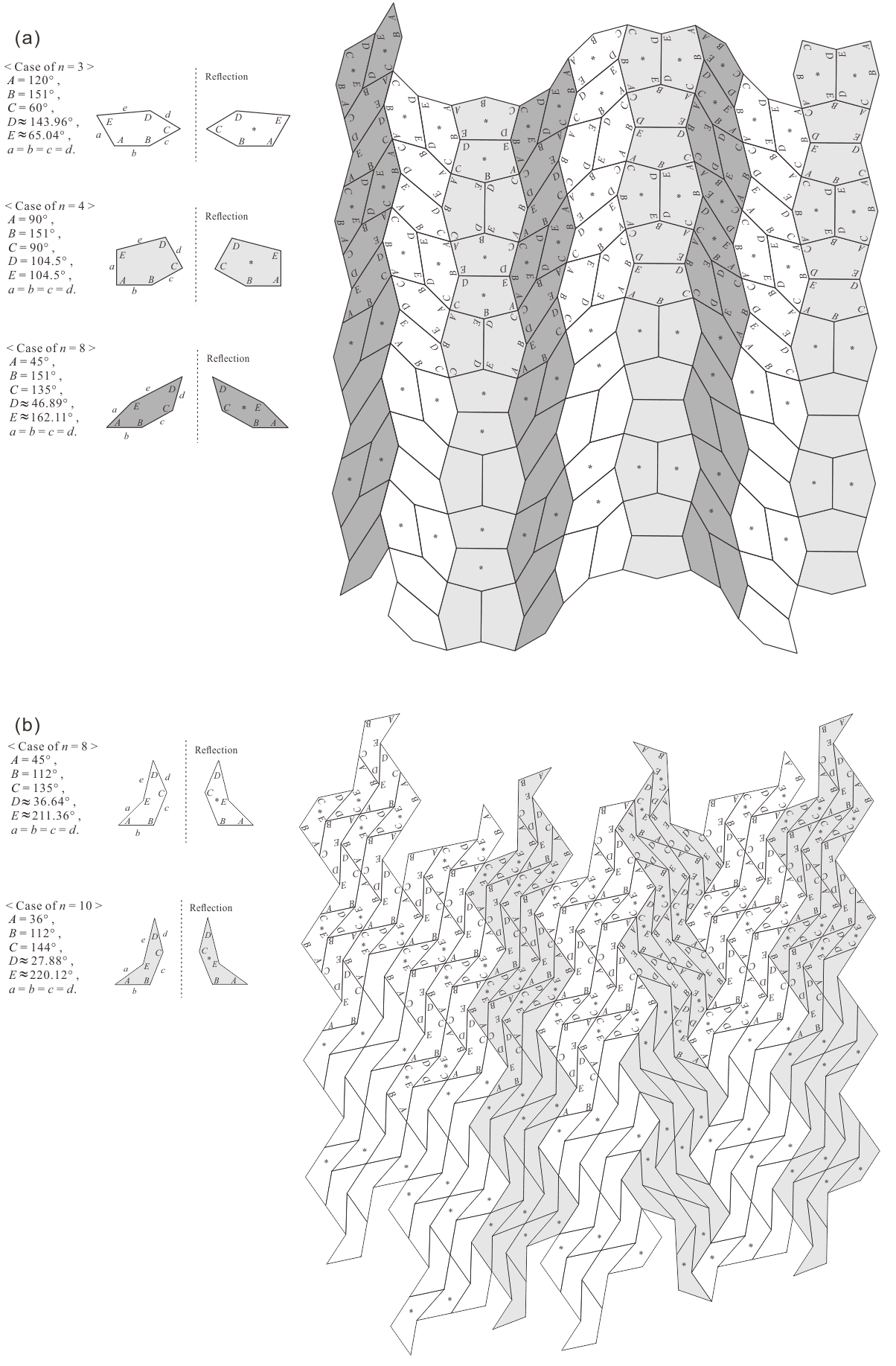} 
  \caption{{\small 
Example of tiling using convex pentagons with $n=3, 4, 8$ in Table~\ref{tab1}, 
and example of tiling by concave pentagons with $n=8, 10$ in Table~\ref{tab3}} 
\label{fig39}
}
\end{figure}

\renewcommand{\figurename}{{\small Figure.}}
\begin{figure}[htbp]
 \centering\includegraphics[width=15cm,clip]{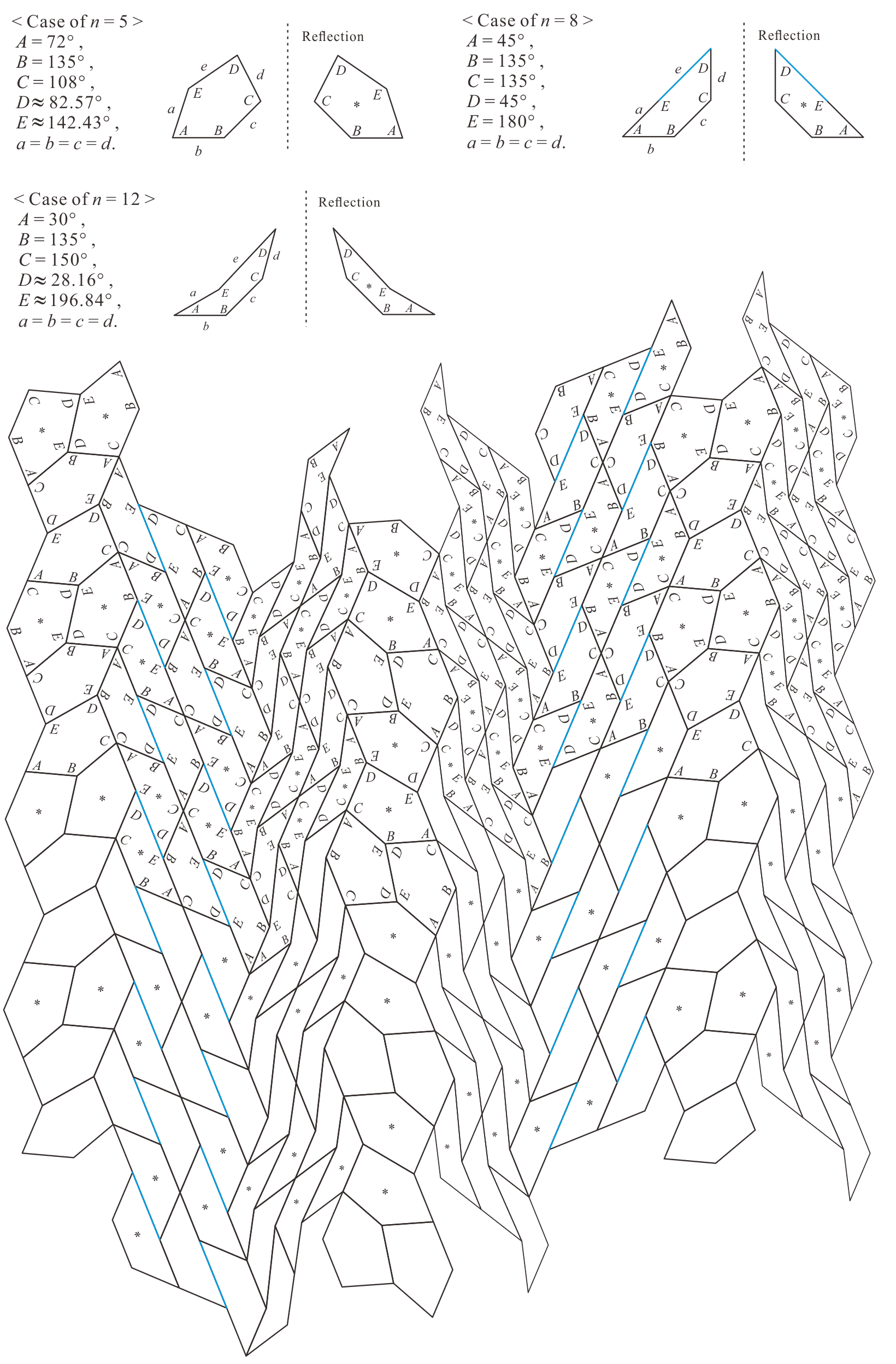} 
  \caption{{\small 
Example of tiling by convex pentagons satisfying (\ref{eq2}) with
 $\alpha = 54^ \circ $ and $\theta = 45^ \circ $, 
trapezoids satisfying (\ref{eq2}) with 
$\alpha = 67.5^ \circ $ and $\theta = 45^ \circ $, 
and concave pentagons satisfying (\ref{eq2}) with 
$\alpha = 75^ \circ $ and $\theta = 45^ \circ $} 
\label{fig40}
}
\end{figure}

\renewcommand{\figurename}{{\small Figure.}}
\begin{figure}[htbp]
 \centering\includegraphics[width=15cm,clip]{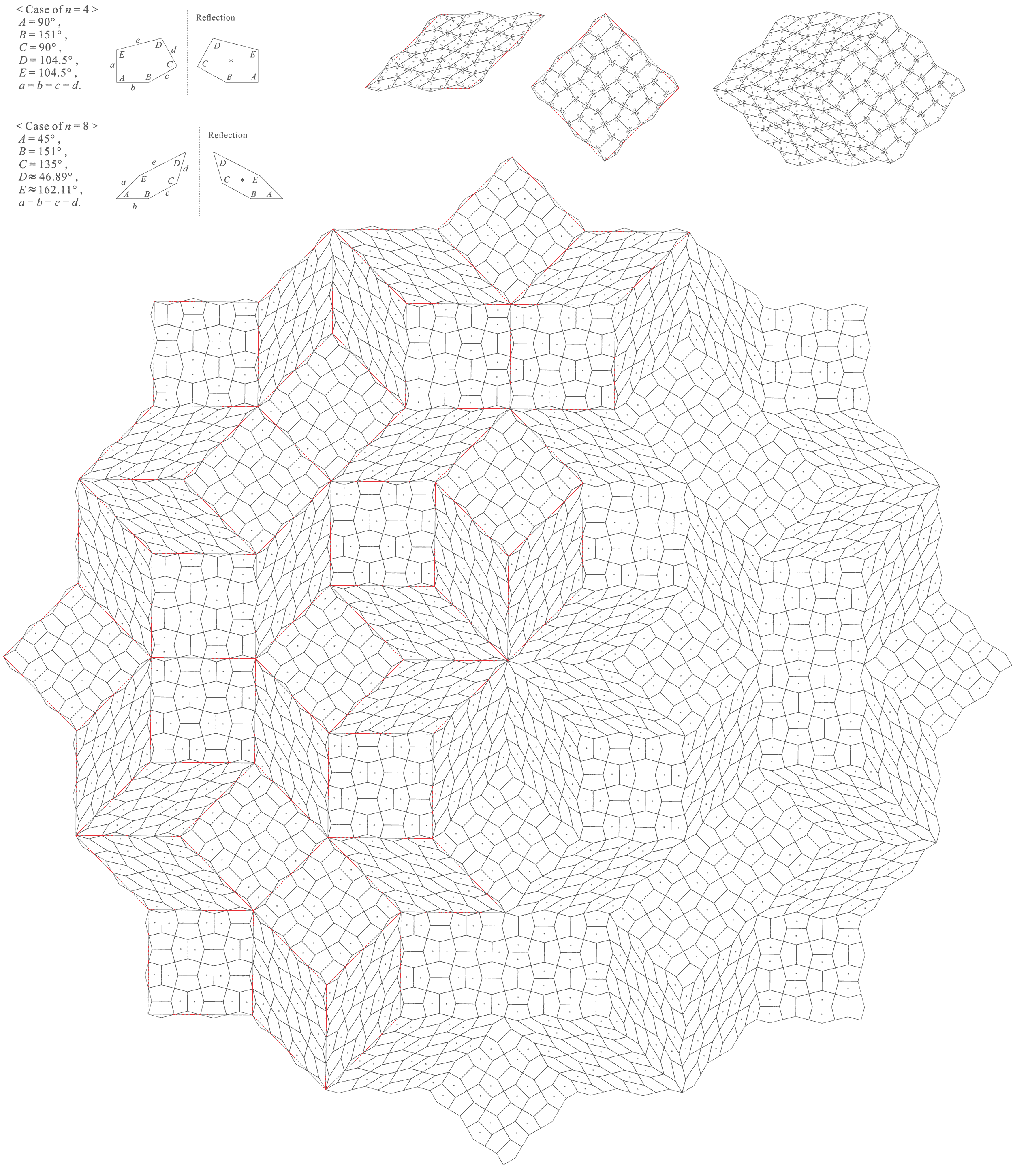} 
  \caption{{\small 
Eight-fold rotationally symmetric tiling by convex pentagons with 
$n=4, 8$ in Table~\ref{tab1} 
(The figure is solely a depiction of the area around the rotationally symmetric 
structure, and the tiling can be spread in all directions)
} 
\label{fig41}
}
\end{figure}

\renewcommand{\figurename}{{\small Figure.}}
\begin{figure}[htbp]
 \centering\includegraphics[width=15cm,clip]{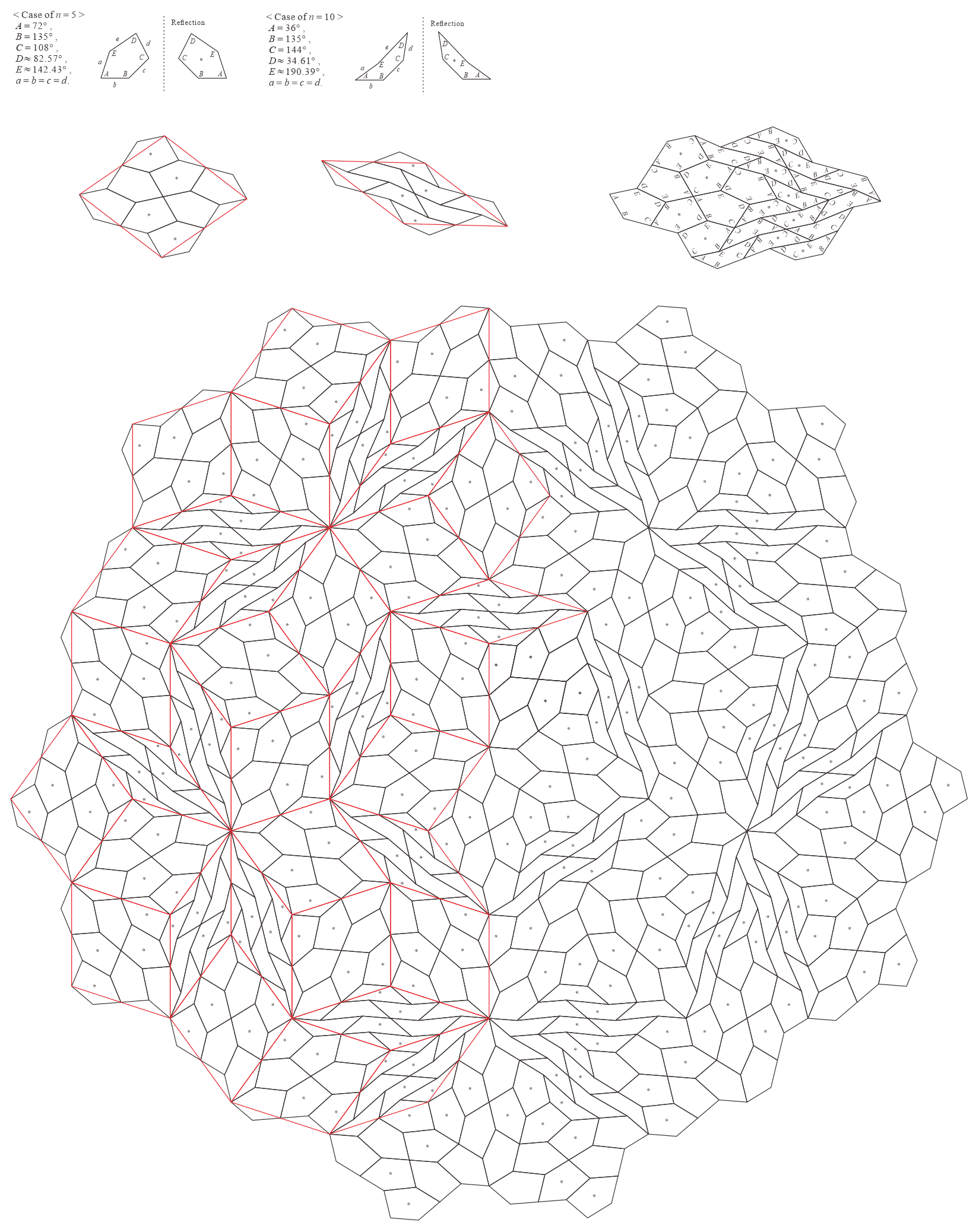} 
  \caption{{\small 
Five-fold rotationally symmetric tiling by convex 
pentagons satisfying (\ref{eq3}) with $n=5$ and $\theta = 45^ \circ $, 
and concave pentagons satisfying (\ref{eq3}) with $n=10$ and 
$\theta = 45^ \circ $ 
(The figure is solely a depiction of the area around the rotationally symmetric 
structure, and the tiling can be spread in all directions)
} 
\label{fig42}
}
\end{figure}

\renewcommand{\figurename}{{\small Figure.}}
\begin{figure}[htbp]
 \centering\includegraphics[width=15cm,clip]{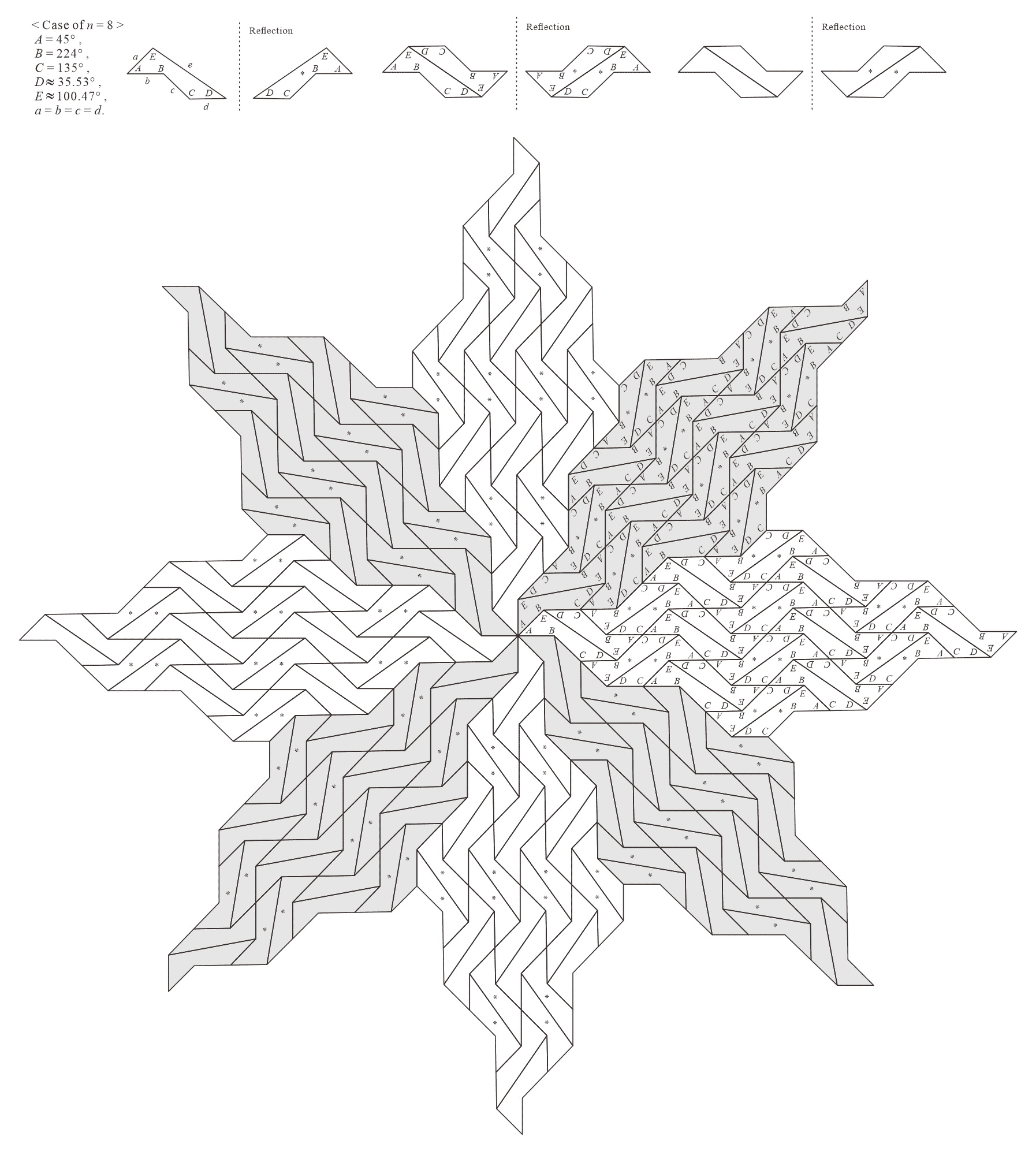} 
  \caption{{\small 
Eight-fold rotationally symmetric tiling by a concave 
pentagon  satisfying (\ref{eq3}) with $n=8$ and $B = 224^ \circ$
 (The figure is solely a depiction of  the area around the rotationally symmetric
 center, and the tiling can be spread in all directions as well)
}  
\label{fig43}
}
\end{figure}

\renewcommand{\figurename}{{\small Figure.}}
\begin{figure}[htbp]
 \centering\includegraphics[width=15cm,clip]{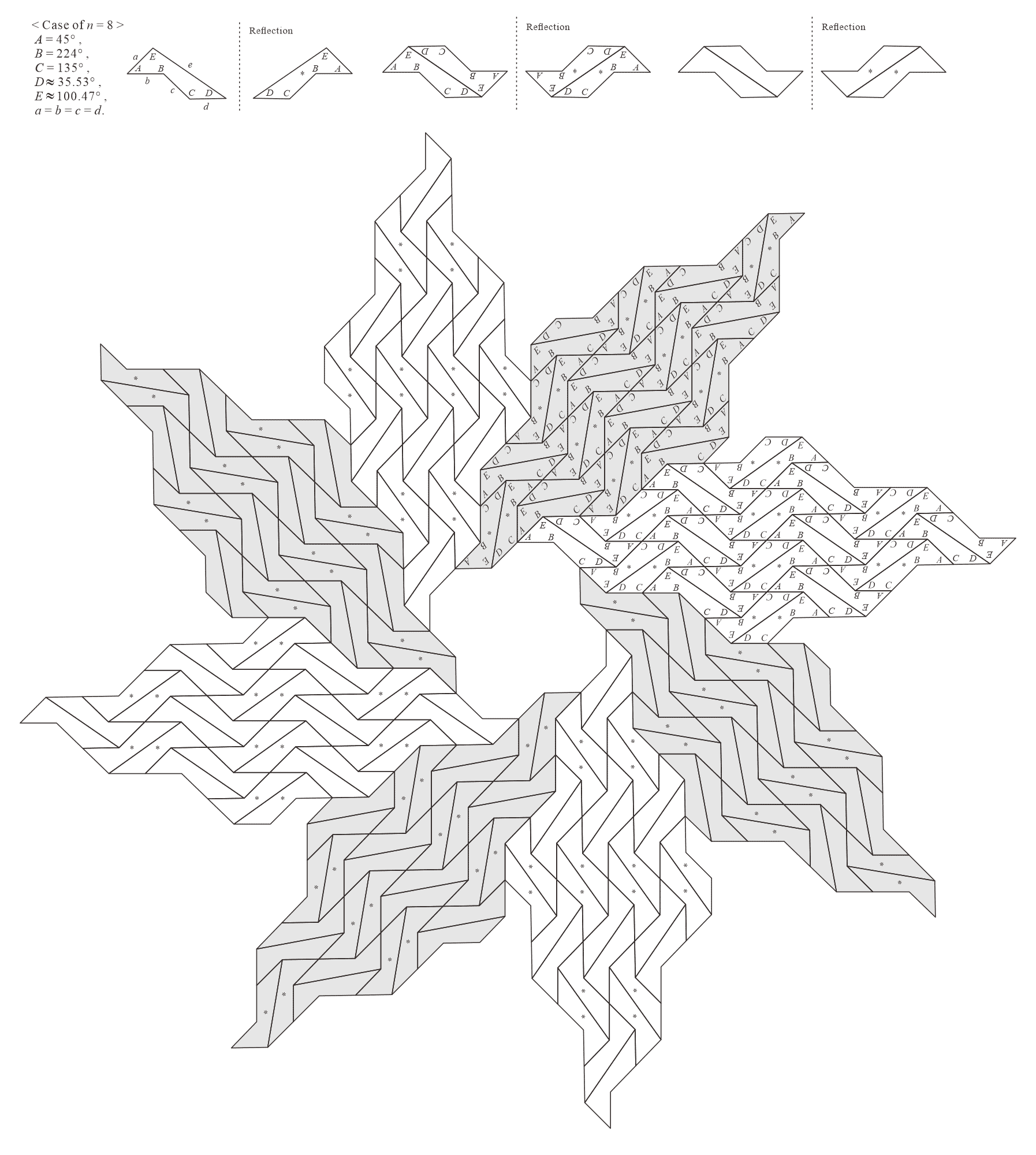} 
  \caption{{\small 
Rotationally symmetric tiling with $C_{4}$ symmetry, with an equilateral 
concave 16-gonal hole with $D_{4}$ symmetry at the center, by a concave 
pentagon  satisfying (\ref{eq3}) with $n = 8$ and $B = 224^ \circ$
 (The figure is solely a depiction of  the area around the rotationally symmetric
 center, and the tiling can be spread in all directions as well)
} 
\label{fig44}
}
\end{figure}


\begin{thebibliography}{99}
{\small

\bibitem{G_and_S_1987}
B.~Gr\"{u}nbaum, G.C.~Shephard, \textit{Tilings and Patterns}. 
W. H. Freeman and Company, New York, 1987, pp.15--35 (Chapter 1), 
pp.471--518 (Chapter 9).


\bibitem{H_and_H_1985}
M. D.~Hirschhorn, D. C.~Hunt, Equilateral convex pentagons which tile 
the plane, \textit{Journal of Combinatorial Theory}, Series \textbf{A 39}, 1--18


\bibitem{H_and_M_2015}
H.~Hosoya, K.~Miyazaki (ed.), \textit{Takakkeihyakka} (in Japanese), Maruzen, 
Tokyo, 2015, pp.172--173.


\bibitem{Klaassen_2016}
B.~Klaassen, Rotationally symmetric tilings with convex pentagons and hexagons, 
\textit{Elemente der Mathematik}, \textbf{71} (2016) 137--144.
doi:10.4171/em/310. 
Available online: \url{https://arxiv.org/abs/1509.06297} 
(accessed on 23 February 2020).


\bibitem{Sugimoto_NoteTP}
T.~Sugimoto, Tiling problem: Convex pentagons for edge-to-edge tiling and 
convex polygons for aperiodic tiling (2015).
\url{https://arxiv.org/abs/1508.01864} (accessed on 23 February 2020).


\bibitem{Sugimoto_2016}
\textemdash, Convex pentagons for edge-to-edge tiling, III, 
\textit{Graphs and Combinatorics}, \textbf{32} (2016) 785--799. 
doi:10.1007/s00373-015-1599-1. 


\bibitem{Sugimoto_2017_2}
\textemdash, Properties of strongly balanced tilings by convex polygons, 
\textit{Research and Communications in Mathematics and Mathematical Sciences}, 
\textbf{8} (2017) 95--114. Available online:
\url{https://www.jyotiacademicpress.net/properties_of_strongly.pdf}
(accessed on 27 February 2022).


\bibitem{Sugimoto_2020_1}
\textemdash, Convex pentagons and concave octagons that can form 
rotationally symmetric tilings (2020).
\url{https://arxiv.org/abs/2005.08470}
(accessed on 19 May 2020).


\bibitem{Sugimoto_2020_2}
\textemdash, Convex pentagons and convex hexagons that can form 
rotationally symmetric tilings (2020).
\url{https://arxiv.org/abs/2005.10639}
(accessed on 22 May 2020).


\bibitem{S_and_O_2009}
T.~Sugimoto, T.~Ogawa, Systematic study of convex pentagonal 
tilings, II: Tilings by convex pentagons with four equal-length edges,
\textit{Forma}, \textbf{24} (2009) 93--109. Available online:
\url{https://forma.katachi-jp.com/abstract/2403/24030093.html} 
(accessed on 4 March 2022).


\bibitem{wiki_pentagon_tiling}
Wikipedia contributors, Pentagonal tiling, 
Wikipedia, The Free Encyclopedia,
\url{https://en.wikipedia.org/wiki/Pentagonal_tiling}
(accessed on 31 March 2020).


\bibitem{wiki_point_group}
\textemdash, Point group, 
Wikipedia, The Free Encyclopedia,
\url{https://en.wikipedia.org/wiki/Point_group}
(accessed on 23 February 2020).


\bibitem{wiki_schoenflies_notation}
\textemdash, Schoenflies notation, 
Wikipedia, The Free Encyclopedia,
\url{https://en.wikipedia.org/wiki/Schoenflies_notation}
(accessed on 23 February 2020).


\bibitem{Zucca_2003}
L.~Zucca, Equilateral pentagons that tile the plane and their 
different tilings, 
\url{https://www.iread.it/lz/pentagons.html}
(accessed on 28 March 2020).


}

\end{thebibliography}
\end{document}